\documentclass{amsart}
\pdfoutput 1
\usepackage[all,hyperref,numberbysection]{bi-discrete}
\usepackage{amsmath}
\usepackage{upgreek}
\usepackage{microtype}
\numberwithin{equation}{section}
\usepackage[a4paper,margin=24ex]{geometry}
\allowdisplaybreaks

\usepackage[backend=biber,maxbibnames=5,maxalphanames=5,style=alphabetic,bibencoding=utf8,giveninits,url=false,isbn=false]{biblatex}
\AtBeginBibliography{\small}

\DeclareMathOperator*{\osc}{\mathrm{osc}}

\begin{document}

\author[Diening]{Lars Diening}
\address{Fakult\"at f\"ur Mathematik, Universit\"at Bielefeld, Postfach 100131, D-33501 Bielefeld, Germany}
\email{lars.diening@uni-bielefeld.de}
\author[Kim]{Kyeongbae Kim}
\address{Department of Mathematical Sciences, Seoul National University, Seoul 08826, Korea}
\email{kkba6611@snu.ac.kr}
\author[Lee]{Ho-Sik Lee}
\address{Fakult\"at f\"ur Mathematik, Universit\"at Bielefeld, Postfach 100131, D-33501 Bielefeld, Germany}
\email{ho-sik.lee@uni-bielefeld.de}
\author[Nowak]{Simon Nowak}
\address{Fakult\"at f\"ur Mathematik, Universit\"at Bielefeld, Postfach 100131, D-33501 Bielefeld, Germany}
\email{simon.nowak@uni-bielefeld.de}

\makeatletter
\@namedef{subjclassname@2020}{\textup{2020} Mathematics Subject Classification}
\makeatother


\thanks{Lars Diening and Simon Nowak gratefully acknowledge funding by the Deutsche Forschungsgemeinschaft (DFG, German Research Foundation) - SFB 1283/2 2021 - 317210226. Kyeongbae Kim thanks for funding of the National Research Foundation of Korea (NRF) through IRTG 2235/NRF-2016K2A9A2A13003815 at Seoul National University. Ho-Sik Lee thanks for funding of the Deutsche Forschungsgemeinschaft through GRK 2235/2 2021 - 282638148.
}

\title{Gradient estimates for parabolic nonlinear nonlocal equations}

\begin{abstract}
	The primary objective of this work is to establish 
	pointwise gradient estimates for solutions to a class of parabolic nonlinear nonlocal measure data problems, expressed in terms of caloric Riesz potentials of the data. As a consequence of our pointwise estimates, we obtain that the first-order regularity properties of solutions to such general parabolic nonlinear nonlocal equations, both in terms of size and oscillations of the spatial gradient, closely resemble the ones of the fractional heat equation even at highly refined scales. 
	Along the way, we show that solutions to homogeneous parabolic nonlinear nonlocal equations have H\"older continuous spatial gradients under optimal assumptions on the nonlocal tails.
\end{abstract}

\maketitle

\tableofcontents


\section{Introduction}
\subsection{Aim and scope}
This paper aims to explore the fine pointwise first-order properties as well as the gradient regularity of solutions to parabolic nonlinear nonlocal equations of the type
\begin{equation}\label{eq:nonlocaleq}
	\partial_t u+\mathcal{L}u=\mu\quad\text{in }\Omega_T \subset \Omega \times (0,T),
\end{equation}
where $T>0$, $\Omega$ is an open subset of $\mathbb{R}^n$ for some $n \geq 2$ and the nonlinear nonlocal operator $\mathcal{L}$ is formally defined by
\begin{equation}\label{eq:nonlocalop}
	\mathcal{L}u(x,t)=(1-s)\mathrm{P.V.}\int_{\bbR^n}\Phi\left(\frac{u(x,t)-u(y,t)}{|x-y|^s}\right)\frac{\,dy}{|x-y|^{n+s}}.
\end{equation}
Here $s \in (0,1)$ is a parameter that determines the order of the nonlocal operator $\mathcal{L}$ given by $2s$, while $\mu$ belongs to the class $\mathcal{M}(\mathbb{R}^{n+1})$ of signed Radon measures on $\mathbb{R}^{n+1}$ with finite total mass. In addition, the nonlinearity $\Phi$ is assumed to satisfy the following Lipschitz and monotonicity assumptions:
\begin{assumption} \label{assump}
	We assume that $\Phi:\bbR\to\bbR$ is an odd function such that for all $t,t^\prime \in \mathbb{R}$ and some $\Lambda \geq 1$, we have
	\begin{equation}
		\label{pt : assmp.phi}
		|{{\Phi}}(t)-\Phi(t^\prime)|\leq \Lambda|t-t^\prime|\quad\text{and}\quad ({\Phi}(t)-{\Phi}(t^\prime))(t-t^\prime)\geq\Lambda^{-1}|t-t^\prime|^{2}.
	\end{equation}
\end{assumption}

Studying the regularity of solutions to elliptic and parabolic nonlinear nonlocal equations has become a highly active research area in recent years, see e.g.\ \cite{KassCalcVar,CSannals,CCV,FKPar,KuuMinSir15,KuuMinSir15s,SerraParabolic,SchikorraMA,IanMosSqu16,DKP,BraLin17,BLS,DFPJDE,MSY,MeH,CKWCalcVar,BKJMA,MeV,GL23,KaWe23,DuzLiao,DuzLiao1} for a non-exhaustive list of fundamental contributions in this direction.
This rapid development of nonlocal regularity theory was largely driven by its wide-ranging applications in both pure and applied mathematics such as for instance stochastic processes of jump-type (see e.g.\ \cite{bertoin,Fukushima}), classical harmonic analysis (see e.g.\ \cite{landkof}), conformal geometry (see e.g.\ \cite{GZ,Case-Chang}), phase transitions (see e.g.\ \cite{fife}), relativistic models (see e.g.\ \cite{LiebYau2}), fluid dynamics (see e.g.\ \cite{NInvent,CaffVass}) and kinetic theory (see e.g.\ \cite{ImbSil}). Moreover, nonlocal operators of the particular type \eqref{eq:nonlocalop} arise in image processing (see e.g.\ \cite{GOsh}).

\subsubsection{Gradient potential estimates for local parabolic equations} \label{sec:111}
A local analogue of the nonlocal equation \eqref{eq:nonlocaleq} in the (formal) limit $s \to 1$ is given by nonlinear second-order parabolic equations of the type
\begin{equation}\label{eq:localeq}
	\partial_t u-\textnormal{div} (\mathbf{a}(\nabla u))=\mu\quad\text{in }\Omega_T \subset \Omega \times (0,T),
\end{equation}
where the vector field $\mathbf{a}$ satisfies suitable growth and ellipticity assumptions. In the case of assumptions on $\mathbf{a}$ that are similar to ours imposed on $\Phi$ in Assumption \ref{assump}, inspired by previous zero-order and first-order potential estimates in the elliptic case provided in \cite{KM,TWAJM,Min}, Duzaar and Mingione in \cite{DM2} managed to prove that solutions to nonlinear parabolic equations of the type \eqref{eq:localeq} satisfy gradient potential estimates of the form
\begin{equation}\label{ineq.potloc}
	\begin{aligned}
		|\nabla u(z_0)|\lesssim I_{1}^{|\mu|}(z_0,R) + \textnormal{lower-order terms}
	\end{aligned}
\end{equation}
for almost every $z_0=(x_0,t_0) \in \Omega_T$ and every $R>0$ such that $Q_R(z_0):=(t_0-R^2,t_0) \times B_R(x_0) \subset \Omega_T$. Here $I_{1}^{|\mu|}(z_0,R)$ denotes a truncated version of the classical caloric Riesz potential of order $1$. Indeed, more generally, for any $\beta \in (0,N_{\textnormal{par}})=(0,n+2)$, we define
$$
I_{\beta}^{|\mu|}(z_0,R)\coloneqq\int_{0}^{R}\frac{|\mu|(Q_r(z_0))}{r^{N_{\textnormal{par}}-\beta}}\frac{\,dr}{r} = \int_{0}^{R}\frac{|\mu|(Q_r(z_0))}{r^{n+2-\beta}}\frac{\,dr}{r},
$$
where $N_{\textnormal{par}}:=n+2$ denotes the standard parabolic dimension. Moreover, for any $z_0=(x_0,t_0) \in \mathbb{R}^{n+1}$ the classical caloric Riesz potential of order $\beta \in (0,N_{\textnormal{par}})$ mentioned above is defined by 
$$
I_{1}^{|\mu|}(z_0)\coloneqq\int_{\mathbb{R}^{n+1}}\frac{d|\mu|(z)}{d_{\textnormal{par}}(z,z_0)^{N_{\textnormal{par}}-\beta}} = \int_{\mathbb{R}^{n+1}}\frac{d|\mu|(z)}{d_{\textnormal{par}}(z,z_0)^{n+2-\beta}},
$$
where the standard parabolic distance $d_{\textnormal{par}}$ in $\mathbb{R}^{n+1}$ is defined by
$$d_{\textnormal{par}}(z,z_0)=d_{\textnormal{par}}((x,t),(x_0,t_0)):=\max \left \{|x-x_0|,|t-t_0|^\frac{1}{2} \right \}.$$
If $\Omega_T=\mathbb{R}^n \times (0,\infty)$, then the estimate \eqref{ineq.potloc} simplifies to one without lower-order terms and in terms of the classical, non-truncated potential, namely
\begin{equation}\label{ineq.potlocg}
	\begin{aligned}
		|\nabla u(z_0)|\lesssim I_{1}^{|\mu|}(z_0).
	\end{aligned}
\end{equation}

In addition to capturing the precise pointwise first-order behavior of solutions to \eqref{eq:localeq}, a major strength of potential estimates of the type \eqref{ineq.potloc}-\eqref{ineq.potlocg} is that they imply sharp Calder\'on-Zygmund-type gradient regularity estimates for solutions in a wide range of function spaces, including those that measure highly refined scales such as Lorentz spaces. 
While for linear parabolic equations such regularity estimates can often also be inferred by means of representation formulas in terms of fundamental solutions/heat kernels, for nonlinear equations of the type \eqref{eq:localeq} this is no longer feasible due to the lack of suitable representation formulas, highlighting the importance of the gradient potential estimate \eqref{ineq.potloc}. Furthermore, variations of the estimate \eqref{ineq.potloc} can be employed to establish sharp borderline criteria on the data in order to ensure control also of the oscillations of the gradient, for instance in the form of VMO regularity or continuity of the gradient of solutions (see \cite{KuMi}). 

Motivated by these powerful implications, in the local setting similar gradient potential estimates were later obtained also for more general nonlinear elliptic and parabolic equations and even systems of $p$-Laplacian-type, see for instance \cite{DM1,KuuMin14w,KuMiARMA1,CMJEMS,Baroni,KuMiV,BCDKS,BY,NNARMA,BCDS,DFJMPA,CKW23,DZJEMS,DZ22} for a non-exhaustive list of further contributions 
 direction.

\subsubsection{Gradient potential estimates for nonlocal parabolic equations} \label{sec:112}

Inspired by zero-order potential estimates for nonlinear nonlocal elliptic equations due to Kuusi, Mingione and Sire from \cite{KuuMinSir15} (see e.g.\ \cite{KLL,Minhyun1,DieNow23,NOS,MinMarv} for more results in this direction), and by the gradient potential estimates for linear elliptic equations with coefficients due to Kuusi, Sire and the last-named author (see \cite{KuuSimYan22}), in \cite{DieKimLeeNow24} we recently managed to establish gradient potential estimates for the elliptic counterpart of the parabolic nonlinear nonlocal equation \eqref{eq:nonlocaleq} in the range $s \in (1/2,1)$.

Moreover, in \cite{NguNowSirWei23} Nguyen, Sire, Weidner and the last-named author obtained zero-order potential estimates for a class of nonlocal drift-diffusion equations related to the dissipative quasi-geostrophic (SQG) equation from fluid dynamics, which in view of linearizing the nonlocal operator \eqref{eq:nonlocalop} in particular implies zero-order potential estimates for the class of parabolic nonlinear nonlocal equations we consider in the present paper.

In light of these recent development and the by now classical parabolic gradient potential estimates from \cite{DM2}, the intriguing question arises whether gradient potential estimates can also be obtained in our parabolic nonlinear nonlocal setting. And in fact, despite the already highly demanding technical nature of the proof of the elliptic case as carried out in \cite{DieKimLeeNow24} and the substantial additional difficulties arising due to the non-stationary nature of the equations we consider, in the present work we establish parabolic first-order potential estimates analogous to \eqref{ineq.potloc} also in the nonlocal case.

While due to the technical nature of their precise formalism, for our gradient potential estimates for solutions to equations posed in bounded domains we refer to Theorem \ref{thm.pt} below, our main result on the whole space can be stated more easily.
Indeed, for any $s \in (0,1)$, denote by $N_{\textnormal{par},s}:=n+2s$ the fractional parabolic dimension of order $s$. Then for any $z_0=(x_0,t_0) \in \mathbb{R}^{n+1}$ and any $\beta \in (0,N_{\textnormal{par},s})$, we denote by 
\begin{equation} \label{eq:GR}
	I_{\beta,s}^{|\mu|}(z_0)\coloneqq \int_{\mathbb{R}^{n+1}}\frac{d|\mu|(z)}{d_{\textnormal{par},s}(z,z_0)^{N_{\textnormal{par},s}-\beta}} = \int_{\mathbb{R}^{n+1}}\frac{d|\mu|(z)}{d_{\textnormal{par},s}(z,z_0)^{n+2s-\beta}}
\end{equation}
a version of the caloric Riesz potential of order $\beta$ suitable for our fractional setting, where the fractional parabolic distance $d_{\textnormal{par},s}$ of order $s$ in $\mathbb{R}^{n+1}$ is defined by
$$d_{\textnormal{par},s}(z,z_0)=d_{\textnormal{par},s}((x,t),(x_0,t_0)):=\max \left \{|x-x_0|,|t-t_0|^\frac{1}{2s} \right \}.$$
We then have the following result.
\begin{theorem}[Gradient potential estimates on the whole space] \label{thm:gradpotwhole}
	Let $s \in (1/2,1)$, $\mu \in \mathcal{M}(\mathbb{R}^{n+1})$ and let $u\in L^2(0,\infty;W^{s,2}(\bbR^n))\cap C(0,\infty;L^2(\bbR^n))$ be a weak solution of 
	\begin{align*}
		\partial_t u+\mathcal{L}u=\mu\quad\text{in }\bbR^n\times (0,\infty).
	\end{align*}
	Moreover, assume that $\Phi$ satisfies Assumption \ref{assump} for some $\Lambda \geq 1$.
	Then for almost every $z_0=(x_0,t_0) \in \mathbb{R}^n \times (0,\infty)$, we have the pointwise estimate
	\begin{equation} \label{eq:gpew}
		|\nabla u(z_0)|\leq cI_{2s-1,s}^{|\mu|}(z_0)
	\end{equation}
	for some constant $c=c(n,s,\Lambda)$. In addition, for any fixed $s_0 \in (1/2,1)$, the constant $c$ depends only on $n,s_0$ and $\Lambda$ whenever $s \in [s_0,1)$.
\end{theorem}

For the precise definition of weak solutions to \eqref{eq:nonlocaleq}, we refer to Definition \ref{def:weak} below. Moreover, in Theorem \ref{thm:gradpotwhole} and all other of our main results we provide estimates that are stable as $s \to 1$. Since at least formally nonlocal operators converge to local second-order ones as the order of the equation approaches two (see e.g.\ \cite{BBM1,FKV} for some rigorous results in this direction), our gradient potential estimates can indeed be considered to be nonlocal analogues of the ones obtained in the local parabolic setting in \cite{DM2}.

As indicated in the previous section, the gradient potential estimates we obtain imply fine regularity results in various function spaces, which in our nonlinear setting is no longer possible by means of estimates on fundamental solutions, see in particular Corollary \ref{cor:Lorentz} below for such regularity results in Lorentz spaces.

Moreover, in order to further strengthen the analogy of parabolic potential estimates playing a similar role for nonlinear parabolic equations as heat kernel estimates do in linear parabolic settings, let us observe that taking $\mu=\delta_{z_1}$ in Theorem \ref{thm:gradpotwhole}, where $\delta_{z_1}$ is the Dirac delta function concentrated at some fixed point $z_1=(x_1,t_1) \in \mathbb{R}^n \times (0,\infty)$, reveals that any solution $u$ to 
\begin{equation} \label{eq:dirac}
	\partial_t u+\mathcal{L}u=\delta_{z_1} \text{ in }\bbR^n\times (0,\infty)
\end{equation}
satisfies the pointwise estimate 
\begin{equation} \label{eq:timedecay}
	|\nabla u(z_0)|\lesssim |t_0-t_1|^{-\frac{n+1}{2s}}
\end{equation}
for any $z_0=(x_0,t_0) \in \mathbb{R}^n \times (0,\infty)$, see Remark \ref{rem:fhk} below. Indeed, the estimate \eqref{eq:timedecay} shows that at least in certain regimes, the gradient of solutions to parabolic nonlinear nonlocal measure data problems of the type \eqref{eq:dirac} possesses similar time-decay as provided by the well-known upper bounds for the gradient of the fractional heat kernel, that is, for the gradient of the fundamental solution of the fractional heat operator $\partial_t+ (-\Delta)^s$, see e.g.\ \cite{BodJak07}. 




\subsubsection{Gradient H\"older regularity for homogeneous parabolic nonlinear nonlocal equations}
A key step in the proof of the pointwise gradient estimates given by \eqref{thm:gradpotwhole} and \eqref{ineq.potl} below is to first prove suitable gradient estimates in the homogeneous case when $\mu \equiv 0$. In this case, H\"older regularity of the spatial gradient of weak solutions to equations similar to \eqref{eq:nonlocaleq} with $\mu \equiv 0$ posed on the whole space was first established by Caffarelli, Chan and Vasseur in \cite{CCV}. 

However, in order to be able to deduce our gradient potential estimates, both on the whole space and in domains, our approach requires localized H\"older estimates for the spatial gradient of solutions to homogeneous equations posed in bounded domains. Due to the nonlocal and nonlinear nature of the operator \eqref{eq:nonlocalop}, obtaining a sharp local analogue of the global gradient estimates from \cite{CCV} is already a nontrivial task. In particular, an additional difficulty present in the parabolic nonlocal setting is the appearance of time-dependent nonlocal tail terms. Handling these tail terms under optimal assumptions on their integrability in time turns out to be a delicate issue, which was recently resolved by Kassmann and Weidner in \cite{KaWe23} in the case of H\"older regularity of solutions to linear parabolic equations by means of certain localization arguments. 

In the present paper, on our way to proving our gradient potential estimates, we show that despite their nonlinear nature, similar localization arguments can be applied in the context of parabolic nonlocal equations of the type \eqref{eq:nonlocaleq}, enabling us in particular to prove localized gradient H\"older regularity for homogeneous parabolic nonlinear nonlocal equations under optimal assumptions on the nonlocal tails on the solution, thereby providing sharp localized counterparts of the global gradient estimates from \cite{CCV}.

Finally, while in the case of parabolic nonlinear nonlocal equations of the type \eqref{eq:nonlocaleq} posed in bounded domains our gradient H\"older estimates seem to be the first of their kind, we want to mention that H\"older estimates for the solution itself rather than their spatial gradients were studied in a substantial amount of previous works, see for instance \cite{FKPar,BLSt,KaWe23,ByuKimKim23,ByuKimKum23p,Tavakoli,NaianCVPDE,APT24,Nai24}.

\subsection{Setup and further main results}

Before being able to state our other main results, we need to introduce our setup more rigorously. In order to control the growth of solutions at infinity, we consider the tail space
$$L^1_{2s}(\mathbb{R}^n):= \left \{g \in L^1_{\loc}(\mathbb{R}^n) \mathrel{\Big|} \int_{\mathbb{R}^n} \frac{|g(y)|}{(1+|y|)^{n+2s}}\,dy < \infty \right \}$$
introduced in \cite{existence} and denote for any open interval $I \subset \mathbb{R}$ and any $q \in [1,\infty]$ by $L^q(I;L^1_{2s}(\bbR^n))$ and $L^q_{\mathrm{loc}}(I;L^1_{2s}(\bbR^n))$ the associated Bochner spaces (see Section \ref{sec2} for more details). We note that a function $g \in L^1_{\loc}(\mathbb{R}^n)$ belongs to the space $L^{1}_{2s}(\mathbb{R}^n)$ if and only if the \emph{nonlocal tails} of $g$ given by
$$ \textnormal{Tail}(g;B_R(x_0)):= (1-s) R^{2s} \int_{\mathbb{R}^n \setminus B_{R}(x_0)} \frac{|g(y)|}{|x_0-y|^{n+2s}}\,dy$$
are finite for all $R>0$ and $x_0 \in \mathbb{R}^n$.

Moreover, for $s \in (0,1)$, $z_0=(x_0,t_0) \in \mathbb{R}^{n+1}$ and $R>0$, we define the parabolic cylinder of order $s$ with radius $R$ and center $z_0$ by $Q_R^s(z_0):=B_R(x_0) \times I^s_R(t_0),$ where $I_R^s(t_0):=(t_0-R^{2s},t_0)$. 

\begin{definition}[Parabolic excess functionals]
	Fix $s \in (0,1)$, $z_0=(x_0,t_0) \in \mathbb{R}^{n+1}$, $R>0$ and $q\in[1,\infty)$. For any $u \in L^q_{\mathrm{loc}}(I_R^s;L^1(B_R(x_0)))$, we define the local parabolic $q$-excess by
	\begin{equation}\label{loc.exc}
		E^q_\loc(u;Q_R^s(z_0)) \coloneqq\left(\dashint_{Q_R(z_0)}|u-(u)_{Q_R^s(z_0)}|^q\,dz\right)^{\frac1q}.
	\end{equation}
	Moreover, for any $u \in L^q_{\mathrm{loc}}(I;L^1_{2s}(\bbR^n))$ we define the nonlocal parabolic $q$-excess by
	\begin{equation}\label{not.exc}
		\begin{aligned}
			E^q(u;Q_R^s(z_0))&\coloneqq\left(\dashint_{Q_R^s(z_0)}|u-(u)_{Q_R^s(z_0)}|^q\,dz\right)^{\frac1q}\\
			&\quad+\left(\dashint_{I_R^s(z_0)}\mathrm{Tail}(u-(u)_{Q_R^s(z_0)};B_R(x_0))^q\,dt\right)^{\frac1q}.
		\end{aligned}
	\end{equation}
Furthermore, when $q=1$, for convenience we write $E_\loc(u;Q_R^s(z_0)) \coloneqq E^1_\loc(u;Q_R^s(z_0))$ and $E(u;Q_R^s(z_0)) \coloneqq E^1(u;Q_R^s(z_0))$.
\end{definition}

Next, we define standard energy-type weak solutions to \eqref{eq:nonlocaleq} as follows.
\begin{definition} \label{def:weak}
	Let $\Omega \subset \mathbb{R}^n$ be an open set and let $\mu\in L^{p_1}_{\mathrm{loc}}(0,T;L^{p_2}_{\mathrm{loc}}(\Omega))$ with $\frac{n}{2p_2s}+\frac{1}{p_1}\leq 1+\frac{n}{4s}$. We say that 
	\begin{equation*}
		u\in L^2_{\mathrm{loc}}(0,T;W^{s,2}_{\mathrm{loc}}(\Omega))\cap C_{\mathrm{loc}}(0,T;L^2_{\mathrm{loc}}(\Omega))\cap L^1_{\mathrm{loc}}(0,T;L^1_{2s}(\bbR^n))
	\end{equation*}
	is a weak solution to \eqref{eq:nonlocaleq}, if 
	\begin{align*}
		&-\int_{t_1}^{t_2}\int_{\Omega}u\partial_t \phi\,dz\\
		&\quad+(1-s)\int_{t_1}^{t_2}\int_{\bbR^n}\int_{\bbR^n}\Phi\left(\frac{u(x,t)-u(y,t)}{|x-y|^s}\right)\frac{\phi(x,t)-\phi(y,t)}{|x-y|^{n+s}}\,dx\,dy\,dt\\
		&=\int_{t_1}^{t_2}\int_{\Omega}\mu\phi\,dz-\int_{\Omega}u\phi\,dx\Bigg\lvert_{t=t_1}^{t=t_2}
	\end{align*}
	holds for any function $\phi\in L^2(t_1,t_2;W^{s,2}(\Omega))\cap W^{1,2}(t_1,t_2;L^2(\Omega))$ with support in the spatial direction compactly contained in $\Omega$, whenever $(t_1,t_2)\Subset (0,T)$.
\end{definition}
\begin{remark} \label{rem:tail} \normalfont
	Under the condition $\frac{n}{2p_2s}+\frac{1}{p_1}\leq 1+\frac{n}{4s}$, the corresponding initial boundary value problem to \eqref{eq:nonlocaleq} is uniquely solvable (see \cite[Appendix A]{ByuKimKim23}).
\end{remark}

We are now in the position to state our next main result, which is concerned with H\"older estimates for the spatial gradient of weak solutions to parabolic nonlinear nonlocal equations of the type \eqref{eq:nonlocaleq} in the homogeneous case when $\mu \equiv 0$. It is noteworthy that in contrast to our gradient potential estimates under general measure data, our H\"older estimates below are valid in the whole range $s \in (0,1)$.

\begin{theorem}[Gradient H\"older regularity]\label{thm.hol}
	Let $s \in (0,1)$ and fix $q>1$. Let 
	\begin{equation*}
		u\in L^2_{\mathrm{loc}}(0,T;W^{s,2}_{\mathrm{loc}}(\Omega))\cap C_{\mathrm{loc}}(0,T;L^2_{\mathrm{loc}}(\Omega))\cap L^q_{\mathrm{loc}}(0,T;L^1_{2s}(\bbR^n))
	\end{equation*}
	be a weak solution of
	$$ \partial_t u+\mathcal{L}u=0 \quad \textnormal{in } \Omega_T.$$ Furthermore, assume that $\Phi$ satisfies Assumption \ref{assump} for some $\Lambda \geq 1$.
	Then there exists some $\alpha=\alpha(n,s,\Lambda,q) \in (0,1)$ such that $\nabla u\in C^{\alpha}_{\mathrm{loc}}(\Omega_T)$. Moreover, for any $z_0=(x_0,t_0) \in \Omega_T$ and any $R>0$ such that $Q_{R}^s(z_0)\Subset \Omega_T$, we have
	\begin{equation} \label{eq:C1a}
		\|\nabla u\|_{L^\infty(Q_{R/2}^s(z_0))}+ R^{\alpha}[\nabla u]_{C^{\alpha}(Q_{R/2}^s(z_0))}\leq E^q(u/R;Q_{R}^s(z_0)),
	\end{equation}
	where $c=c(n,s,\Lambda,q)$. In addition, for any fixed $s_0 \in (0,1)$, the constants $c$ and $\alpha$ depend only on $n,s_0,\Lambda$ and $q$ whenever $s \in [s_0,1)$.
\end{theorem}
\begin{remark}[Sharpness] \normalfont
	We note that in Theorem \ref{thm.hol} the tail assumption that $u \in L^q_{\mathrm{loc}}(0,T;L^1_{2s}(\bbR^n))$ for some $q>1$ is sharp, since under the slightly weaker assumption that $u \in L^1_{\mathrm{loc}}(0,T;L^1_{2s}(\bbR^n))$, already in the linear case of the fractional heat equation, that is, for $\Phi(t)=c_{n,s}t$ and $\mu \equiv 0$, where $c_{n,s}$ is some appropriate positive constant such that $\mathcal{L}=(-\Delta)^s$ is the fractional Laplacian, weak solutions to \eqref{eq:nonlocaleq} are in general not locally H\"older continuous in time, see \cite[Example 5.2]{KaWe23}.
\end{remark}

Since under general measure data weak solutions to \eqref{eq:nonlocaleq} in the sense of Definition \ref{def:weak} might not exist, we are going to state our gradient potential estimates in bounded domains and their consequences in terms of the following more general solution concept called SOLA (= solutions obtained by limiting approximations). 

\begin{definition} \label{def:SOLA}
	Let $\mu\in\mathcal{M}(\mathbb{R}^{n+1})$, $g\in L^{2}\left(0,T;W^{s,2}(\Omega)\right) \cap L^{1}\left(0,T;L^{1}_{2s}(\mathbb{R}^{n})\right)$ and $g_{0}\in L^{1}(\Omega)$. We say that a function $u\in L^{p}(0,T;W^{\sigma,p}(\Omega))\cap L^{\infty}(0,T;L^{1}(\Omega))\cap L^{1}\left(0,T;L^{1}_{2s}(\mathbb{R}^{n})\right)$, for $\sigma\in(0,s)$ and $p\in\left[1,\frac{n+2s}{n+s}\right)$ is a SOLA of the initial boundary-value problem
	\begin{equation}
		\label{eq: IVP}
		\left\{
		\begin{alignedat}{3}
			\partial_t u+\mathcal{L}u&= \mu&&\qquad \mbox{in  $\Omega_{T}$}, \\
			u&=g&&\qquad  \mbox{in $\big(\mathbb{R}^{n}\setminus\Omega\big)\times(0,T]$},\\
			u(\cdot,0)&=g_{0}&&\qquad  \mbox{in $\Omega$},
		\end{alignedat} \right.
	\end{equation}
	if $u$ satisfies 
	\begin{equation}\label{defn.sola}
		\begin{aligned}
			&-\int_{\Omega_T}u\phi_t\,dz\\
			&\quad+(1-s)\int_{0}^T\int_{\bbR^n}\int_{\bbR^n}\Phi\left(\frac{u(x,t)-u(y,t)}{|x-y|^s}\right)\frac{\phi(x,t)-\phi(y,t)}{|x-y|^{n+s}}\,dx\,dy\,dt=\int_{\Omega_T}\phi\,d\mu
		\end{aligned}
	\end{equation}
	for any $\phi\in C_c^\infty(\Omega_T)$,
	$u=g$ a.e. in $\left(\mathbb{R}^{n}\setminus\Omega\right)\times (0,T)$ and 
	\begin{equation}
		\label{defn sola : initial data}
		\lim_{h\to0}\frac{1}{h}\int_{0}^{h}\|u(\cdot,t)-g_{0}\|_{L^{1}(\Omega)}\,dt=0.
	\end{equation} Moreover, there exists a sequence of weak solutions $$\{u_{i}\}_{i\in\mathbb{N}}\subset C([0,T];L^{2}(\Omega))\cap L^{2}\left(0,T;W^{s,2}(\mathbb{R}^{n})\right)$$ to the regularized problems
	\begin{equation}
		\label{thm eq : regularized1}
		\left\{
		\begin{alignedat}{3}
			\partial_tu_i+\mathcal{L}u_{i}&= \mu_{i}&&\qquad \mbox{in  $\Omega_{T}$}, \\
			u_{i}&=g_{i}&&\qquad  \mbox{in $\big(\mathbb{R}^{n}\setminus\Omega\big)\times[0,T]$},\\
			u_{i}(\cdot,0)&=g_{0,i}&&\qquad  \mbox{in $\Omega$},
		\end{alignedat} \right.
	\end{equation}
	where $\mu_{i}\in C_{c}^{\infty}(\mathbb{R}^{n}\times(0,T))$, $g_i\in L^2(0,T;W^{s,2}(\bbR^n))$, and $g_{0,i}\in L^2(\Omega)$ satisfy
	\begin{equation}
		\label{regular limit}
		\left\{
		\begin{alignedat}{3}
  &u_j\to u&&\quad\mbox{a.e.\ in $\bbR^n\times(0,T)$ and in $L^1_{\mathrm{loc}}(\bbR^n\times(0,T))$}\\
	&\mu_{i}\rightharpoonup\mu&&\quad\text{in the sense of measures}\\
   &g_i\to g&&\quad\text{in }L^2(0,T;W^{s,2}(\Omega))\cap L^1(0,T;L^1_{2s}(\bbR^n))\\
   &g_{0,i}\to g_0&&\quad\text{in }L^1(\Omega)
		\end{alignedat} \right.
	\end{equation}
	and
	\begin{equation*}
		\limsup_{i\to\infty}|\mu_i|(Q)\leq |\mu|(\bar{Q})
	\end{equation*}
	for every $Q\subset\mathbb{R}^{n}\times(0,T)$.
\end{definition}

The main advantage of working with SOLA instead of standard weak solutions is that SOLA always exist even in the presence of general measure data, as our next result shows, which is a parabolic counterpart of \cite[Theorem 1.1]{KuuMinSir15}.
\begin{theorem}[Existence of SOLA]
	\label{thm: existence}
	Let $s \in (0,1)$, $\mu\in\mathcal{M}(\mathbb{R}^{n+1})$, $g\in L^{2}\left(0,T;W^{s,2}_{\mathrm{loc}}(\mathbb{R}^{n})\right) \cap L^{q}\left(0,T;L^{1}_{2s}(\mathbb{R}^{n})\right)$ for some $q\in[1,\infty]$ such that $\partial_t g\in \left(L^{2}\left(0,T;W^{s,2}(\Omega)\right)\right)^{*}$ and assume that $\Phi$ satisfies Assumption \ref{assump} for some $\Lambda \geq 1$. Then, for any $p\in\left[1,\frac{n+2s}{n+s}\right)$ and any $\sigma\in(0,s)$, there exists a SOLA 
	\begin{equation*}
		u\in L^{p}(0,T;W^{\sigma,p}(\Omega))\cap L^\infty(0,T;L^1(\Omega))\cap L^{q}(0,T;L^{1}_{2s}(\mathbb{R}^{n}))
	\end{equation*}to \eqref{eq: IVP} with $g_{0}(x)=g(x,0)$.
\end{theorem}

For any $s \in (0,1)$, $\beta\in(0,n+2s)$, $z_0 =(x_0,t_0) \in \mathbb{R}^{n+1}$, $R>0$ and $\mu \in \mathcal{M}(\mathbb{R}^{n+1})$, we define a truncated version of the caloric Riesz-type potential \eqref{eq:GR} by
\begin{equation*}
	I_{\beta,s}^{|\mu|}(z_0,R)\coloneqq\int_{0}^{R}\frac{|\mu|(Q_r^s(z_0))}{r^{n+N_{\textnormal{par},s}-\beta}}\frac{\,dr}{r} = \int_{0}^{R}\frac{|\mu|(Q_r^s(z_0))}{r^{n+2s-\beta}}\frac{\,dr}{r}.
\end{equation*}

We are now in the position to state our gradient potential estimates for SOLA of initial boundary-value problems.
\begin{theorem}[Gradient potential estimates for SOLA]\label{thm.pt}
	Let $s \in (1/2,1)$ and let $u$ be a SOLA to \eqref{eq: IVP} with $\mu$, $g$ and $g_0$ as in Definition \ref{def:SOLA}. Furthermore, assume that $\Phi$ satisfies Assumption \ref{assump} for some $\Lambda \geq 1$. Then for almost every $z_0=(x_0,t_0) \in \Omega_T$ and any $R>0$ such that $Q_R^s(z_0)\Subset \Omega_T$, we have  
	\begin{equation}\label{ineq.potl}
		\begin{aligned}
			|\nabla u(z_0)|&\leq cE(u/R;Q_R^s(z_0))+cI_{2s-1,s}^{|\mu|}(z_0,R)\\
			&\quad+c\int_{0}^{R}\left(\int_{I_r^s(t_0)}R^{-2s}\mathrm{Tail}(u-(u)_{Q_R^s(z_0)};B_R(x_0))\,dt\right)\frac{dr}{r^2},
		\end{aligned}
	\end{equation}
	where $c=c(n,s,\Lambda)$. In addition, for any fixed $s_0 \in (1/2,1)$, the constant $c$ depends only on $n,s_0$ and $\Lambda$ whenever $s \in [s_0,1)$.
\end{theorem}

As mentioned, a key strength of the potential estimate \eqref{ineq.potl} is that since the mapping properties of the caloric Riesz-type potential $I_{2s-1,s}$ can easily be inferred with respect to many function spaces, as an immediate corollary we obtain gradient regularity estimates even in function spaces measuring highly refined scales. For instance, recall that the Lorentz spaces $L^{p,q}$, $p \in [1,\infty]$, $q \in (0,\infty]$ refine the scale of $L^p$ spaces in the sense that $L^{p,p}(\Omega_T)=L^p(\Omega_T)$ and $L^{p,q_0}(\Omega_T) \subsetneq L^{p,q_1}(\Omega_T)$ whenever $q_0<q_1$. For a precise definition of Lorentz spaces and more relations between them, see e.g.\ \cite[Section 1.3]{KuuSimYan22}.

Theorem \ref{thm.pt} now yields the following gradient regularity estimates of Calder\'on-Zygmund-type in Lorentz spaces.
\begin{corollary}[Calder\'on-Zygmund estimates in Lorentz spaces]\label{cor:Lorentz}
	Let $s \in (1/2,1)$ and let $u$ be a SOLA to \eqref{eq: IVP} with $\mu$, $g$ and $g_0$ as in Definition \ref{def:SOLA}. Moreover, assume that $\Phi$ satisfies Assumption \ref{assump} for some $\Lambda \geq 1$.
	\begin{itemize}
		\item We have the implication
		$$
		\mu\in \mathcal{M}(\setR^{n+1})\Longrightarrow \nabla u\in L^{\frac{n+2s}{n+1},\infty}_{\mathrm{loc}}(\Omega_T).
		$$
		\item If $p \in \left (1,\frac{n+2s}{2s-1} \right )$ and $q \in (0,\infty]$, then we have the implication
		$$
		\mu\in L^{p,q}(\Omega_T), \, u\in L^{p,q}(0,T;L^1_{2s}(\bbR^n)) \Longrightarrow \nabla u\in L_\loc^{\frac{p(n+2s)}{n+2s-(2s-1)p},q}(\Omega_T).
		$$
	\end{itemize}
\end{corollary}
In particular, the second implication in Corollary \ref{cor:Lorentz} yields the slightly coarser implication in standard $L^p$ spaces: For any $p \in \left (1,\frac{n+2s}{2s-1} \right )$, we have
\begin{equation} \label{eq:Lpreg}
	\mu\in L^{p}(\Omega_T), \, u\in L^{p}(0,T;L^1_{2s}(\bbR^n)) \Longrightarrow \nabla u\in L_\loc^{\frac{p(n+2s)}{n+2s-(2s-1)p}}(\Omega_T).
\end{equation}


While the gradient potential estimate \eqref{ineq.potl} yields precise control of the size of $\nabla u$ in terms of the size of the data, we also provide fine control of the oscillations of $\nabla u$ if the data and the long-range interactions of $u$ are sufficiently well-behaved. In particular, we obtain the following gradient continuity criterion via potentials.
\begin{theorem}[Gradient continuity via potentials]\label{thm.conti}
	Let $s \in (1/2,1)$ and let $u$ be a SOLA to \eqref{eq: IVP} with $\mu$, $g$ and $g_0$ as in Definition \ref{def:SOLA}. Moreover, assume that $\Phi$ satisfies Assumption \ref{assump} for some $\Lambda \geq 1$. Let us fix $z_0=(x_0,t_0) \in \Omega_T$ and $R>0$ such that $Q_{R}^s(z_0)\Subset \Omega_T$. If 
	\begin{align}\label{ass.coni}
		\lim_{\varrho\to0} \sup_{z_1 \in Q_R^s(z_0)}\left[I^{|\mu|}_{2s-1,s}(z_1,\varrho)+\int_{0}^{\varrho}\int_{I_r^s(t_1)}\mathrm{Tail}(u-(u)_{Q_R^s(z_1)};B_R(x_1))\,dt\frac{\,dr}{r^{2}}\right]=0,
	\end{align}
	then $\nabla u$ is continuous in $Q_{R/2}^s(z_0)$.
\end{theorem}
Combining Theorem \ref{thm.conti} with Lemma \ref{lem.lorentz} below, we directly obtain the following borderline criterion for gradient continuity in terms of Lorentz spaces.
\begin{corollary}[Lorentz spaces criterion for gradient continuity] \label{cor:L}
	Let $s \in (1/2,1)$ and let $u$ be a SOLA to \eqref{eq: IVP} with $\mu$, $g$ and $g_0$ as in Definition \ref{def:SOLA}. Moreover, assume that $\Phi$ satisfies Assumption \ref{assump} for some $\Lambda \geq 1$. If 
	\begin{equation*}
		\mu\in L^{\frac{n+2s}{2s-1},1}(\Omega_T)\quad\text{and}\quad u\in L^{\frac{2s}{2s-1},1}(0,T;L^1_{2s}(\bbR^n)),
	\end{equation*}
	then $\nabla u$ is continuous in $\Omega_T$.
\end{corollary}

\begin{remark} \normalfont
	The assumptions $u\in L^{p,q}(0,T;L^1_{2s}(\bbR^n))$ in Corollary \ref{cor:Lorentz} and $u\in L^{\frac{2s}{2s-1},1}(0,T;L^1_{2s}(\bbR^n))$ from Corollary \ref{cor:L} only restrict the global behavior of $u$, since under the assumptions on $\mu$ in the mentioned corollaries, \cite[Corollary 1.1]{NguNowSirWei23} already implies that $u\in L^{p,q}_{\mathrm{loc}}(0,T;L^1_\loc(\Omega))$ in the second part of Corollary \ref{cor:Lorentz} and that $u\in L^{\frac{2s}{2s-1},1}_{\mathrm{loc}}(0,T;L^1_\loc(\Omega))$ in Corollary \ref{cor:L}. 
\end{remark}
\begin{remark} \normalfont
	Note that for $\mu \equiv 0$, Theorem \ref{thm.hol} yields gradient continuity under the weaker assumption that $u\in L^{q}(0,T;L^1_{2s}(\bbR^n))$ for some $q>1$ in comparison to the assumption $u\in L^{\frac{2s}{2s-1},1}(0,T;L^1_{2s}(\bbR^n))$ made in Corollary \ref{cor:L}. The reason for this discrepancy might be technical: While in the homogeneous setting of Theorem \ref{thm.hol} we are able to reduce the tail assumption to an optimal one by means of a bootstrap argument that involves differentiating the equation, this is no longer possible in the setting of Corollary \ref{cor:L} due to the presence of the non-differentiable data $\mu$. Thus, investigating whether it is possible to reduce the assumptions on the tails of $u$ made in Corollary \ref{cor:L} and also in the other consequences of Theorem \ref{thm.pt} represents an interesting open problem.
	
	Nevertheless, to the best of our knowledge, the tail assumptions in all of our results are already considerably weaker than in the previous literature concerned with the higher regularity of solutions to parabolic nonlocal problems. Indeed, while the recent works \cite{KaWe23,BKJEE} provide H\"older regularity of $u$ for some in general very small H\"older exponent under sharp tail assumptions, all contributions that provide more regularity seem to assume that at least $u\in L^{\infty}(0,T;L^1_{2s}(\bbR^n))$ (see e.g.\ \cite{BLSt,ByuKimKim23,ByuKimKum23p,Tavakoli}), which is a significantly stronger assumption than any of our assumptions on the tails of $u$.
\end{remark}

Finally, we also provide the following criterion for VMO gradient regularity, which yields slightly weaker control on the oscillations of $\nabla u$ than gradient continuity under slightly weaker assumptions on the data. 
\begin{theorem}[VMO gradient regularity via potentials]\label{thm.vmo} 	Let $s \in (1/2,1)$ and let $u$ be a SOLA to \eqref{eq: IVP} with $\mu$, $g$ and $g_0$ as in Definition \ref{def:SOLA}. Moreover, assume that $\Phi$ satisfies Assumption \ref{assump} for some $\Lambda \geq 1$. If for some $z_0=(x_0,t_0) \in \Omega_T$ and $R>0$ such that $Q_{R}^s(z_0)\Subset \Omega_T$ we have
	\begin{equation}\label{ass1.vmo}
		\sup_{z_1\in Q_R^s(z_0)}\left[I^{|\mu|}_{2s-1,s}(z_1,\varrho)+\int_{0}^{\varrho}\int_{I_r^s(t_1)} \mathrm{Tail}(u-(u)_{Q_R^s(z_1)};B_R(x_1))\,dt\frac{\,dr}{r^{2}}\right]<\infty
	\end{equation}
	and
	\begin{equation}\label{ass2.vmo}
		\lim_{\varrho\to0}\sup_{z_1\in Q_R^s(z_0)}\left[\frac{|\mu|(Q_\varrho^s(z_1))}{\varrho^{n+1}}+\frac{1}{\varrho}\int_{I_\varrho(t_1)} \mathrm{Tail}(u-(u)_{Q_R^s(z_1)};B_R(x_1))\,dt\right]=0,
	\end{equation}
	then $\nabla u \in \textnormal{VMO}(Q^s_{R/2}(z_0))$.
\end{theorem}

Before concluding this section, let us shed some light on the connection of parabolic potential estimates and heat kernel estimates that we already mentioned in Section \ref{sec:112}.
\begin{remark}[Potential estimates and heat kernel estimates] \label{rem:fhk} \normalfont
	In view of the gradient potential estimate \eqref{eq:gpew} from Theorem \ref{thm:gradpotwhole}, any solution $u$ to
	\begin{equation} \label{eq:dirac1}
		\partial_t u+\mathcal{L}u=\delta_{z_1} \text{ in }\bbR^n\times (0,\infty),
	\end{equation}
	where $\delta_{z_1}$ is the Dirac delta function concentrated at some fixed point $z_1=(x_1,t_1) \in \mathbb{R}^n \times (0,\infty)$,
	satisfies the pointwise estimate 
	\begin{equation} \label{eq:timedecay1}
		\begin{aligned}
			|\nabla u(z_0)|\lesssim I_{2s-1,s}^{|\mu|}(z_0) & \lesssim \min\left\{|x_0-x_1|^{-(n+1)},|t_0-t_1|^{-\frac{n+1}{2s}}\right\} \\ & \lesssim |t_0-t_1|^{-\frac{n+1}{2s}}
		\end{aligned}
	\end{equation}
	for any $z_0=(x_0,t_0) \in \mathbb{R}^n \times (0,\infty)$.
	
	In the linear case when $\Phi(t)=c_{n,s} t$ for some appropriate positive constant $c_{n,s}$ so that $\mathcal{L}=(-\Delta)^s$ is the fractional Laplacian, the fractional heat kernel $p(t,x,y)$, that is, the fundamental solution of $\partial_t u + (-\Delta)^s u$ in $\mathbb{R}^n$, has the property that for any $z_1=(x_1,t_1) \in \mathbb{R}^n \times (0,\infty)$ with $t>t_1$,
	\begin{equation} \label{eq:heatrep}
		v(x,t)\coloneqq \int_{0}^\infty\int_{\bbR^{n}}P(t-\tau,x,y) \, d\delta_{z_1}(y,\tau)=P(t-t_1,x,x_1)
	\end{equation}
	is a solution to
	\begin{equation} \label{eq:diracl}
		\partial_t v+(-\Delta)^s v=\delta_{z_1} \text{ in }\bbR^n\times (0,\infty),
	\end{equation}
	so that \eqref{eq:timedecay1} implies the heat kernel gradient estimate
	\begin{align*}
		|\nabla P(t_0-t_1,x_0,x_1)|=|\nabla v(z_0)|\lesssim |t_0-t_1|^{-\frac{n+1}{2s}}
	\end{align*}
	for all $z_0=(x_0,t_0),z_1=(x_1,t_1) \in \mathbb{R}^n \times (0,\infty)$ with $t_0>t_1$. At least in the regime $|x_0-x_1| \eqsim |t_0-t_1|^\frac{1}{2s}$, this corresponds to the known optimal time-decay for the gradient of the fractional heat kernel, since in this regime the sharp gradient estimates for the fractional heat kernel given e.g.\ in \cite[Lemma 5]{BodJak07} yield
	$$
		|\nabla P(t_0-t_1,x_0,x_1)|\eqsim |x_0-x_1|\min\left\{\frac{t_0-t_1}{|x_0-x_1|^{n+2+2s}},(t_0-t_1)^{-\frac{n+2}{2s}}\right\} \eqsim (t_0-t_1)^{-\frac{n+1}{2s}}
	$$
	whenever $|x_0-x_1|\eqsim |t_0-t_1|^{\frac1{2s}}$. 
	
	Therefore, in the above sense our potential estimates can indeed be considered to be nonlinear analogues of corresponding heat kernel estimates in the linear setting.
\end{remark}

\subsection{Technical approach}
Let us now briefly outline the approach we take to prove our main results in a heuristic manner, with a particular focus on the novelties compared to previous approaches for obtaining gradient potential estimates.

In the local parabolic setting considered in \cite{DM2}, the authors establish the gradient potential estimate \eqref{ineq.potloc} by means of a potential-theoretic Campanato-type iteration below the natural duality exponent in terms of the local excess functional given by \eqref{loc.exc} with the choice $q=s=1$ and with $u$ replaced by $\nabla u$. Roughly speaking, this approach is predicated on establishing decay estimates for the gradient of solutions to the corresponding equations with zero right-hand side by means of De Giorgi-Nash-Moser theory along with differentiating the equation. These gradient excess decay estimates are then transferred to solutions of the associated measure data problem through gradient comparison estimates, which eventually leads to the desired gradient potential estimates.

In our nonlocal setting, there arise several additional difficulties in comparison to the local second-order case:
\begin{itemize}
	\item The presence of non-differentiable tail terms due to the nonlocality of the equation.
	\item The lack of obvious energy estimates at the gradient level due to the lower order of the equation.
\end{itemize}

In the elliptic nonlocal case treated in \cite{DieKimLeeNow24}, we overcame the mentioned difficulties by combining the Campanato-type methods introduced in \cite{DM2,KuuMinSir15} with certain localization arguments and difference quotient techniques based on a nonlinear atomic decomposition originally introduced in \cite{KrMin} in the study of variational problems, which was later utilized as a tool to differentiate measure data problems in \cite{MinCZMD}, see also \cite{KrMin1,Min,AKMARMA,DFMInvent,DFM2} for further applications of such methods.

While this general philosophy established in \cite{DieKimLeeNow24} can also be put into practice in the parabolic nonlinear nonlocal setting studied in the present work, a number of additional intricacies present in our parabolic setting nevertheless lead to severe complications in comparison to the elliptic case studied in \cite{DieKimLeeNow24}, which require new ideas to be dealt with. Indeed, as already indicated, in \cite{DieKimLeeNow24} the appearance of non-differentiable tail terms was surmounted by means of cutoff arguments that essentially enabled us to treat the nonlocal tails as a right-hand side, which turns out to be suitably regular. In our parabolic setting, not only differentiability issues, but also integrability issues of the nonlocal tails with respect to the time variable arise. In the case of our H\"older regularity result for equations with zero right-hand side given by Theorem \ref{thm.hol}, we address this issue by combining the cutoff arguments developed in the elliptic nonlinear nonlocal setting in \cite{DieKimLeeNow24} with further localization arguments recently introduced in the linear parabolic nonlocal setting in \cite{KaWe23}, which allows to establish gradient H\"older regularity under optimal tail assumptions, as noted in Remark \ref{rem:tail}.

In addition to being interesting for their own sake, the estimates in the homogeneous case given by Theorem \ref{thm.hol} form the base in order to establish our gradient potential estimates under general measure data. However, in the process of transferring the information obtained in the homogeneous case to the setting of measure data, more severe difficulties in comparison to the elliptic case arise. Most of these complications find their root in the absence of suitable Poincar\'e-type inequalities on the parabolic cylinders that respect the natural space-time scaling exhibited by the type of equations we study.

Let us make this point more precise.
Indeed, as indicated above, in the elliptic case treated in \cite{DieKimLeeNow24} we apply nonlinear atomic decomposition methods for the following multitude of purposes. First of all, in \cite{DieKimLeeNow24} we use such difference quotient methods to prove that the gradient of solutions to elliptic nonlinear nonlocal measure data problems is locally integrable and belongs to a certain range of fractional Sobolev spaces, which in contrast to the local setting is a priori not known in the nonlocal case. Secondly, similar methods are applied to provide suitable first-order excess decay estimates that respect both the nonlocality of the equation as well as the lack of boundary regularity estimates for solutions. Thirdly, the obtained estimates are then combined in order to upgrade known zero-order comparison estimates to suitable first-order comparison estimates, enabling us to conclude the proof  in the elliptic case by nonlocal adaptations of the Campanato-type methods from \cite{DM2}. In contrast, in our parabolic setting, due to the mentioned lack of suitable Poincar\'e-type inequalities on parabolic cylinders, we cannot use such nonlinear atomic decomposition methods in order to prove suitable higher differentiability estimates for the spatial gradient in space and time simultaneously. Instead, in sharp contrast to the local setting, space and time differentiability need to be treated \emph{separately}, leading to a propagation of technical issues throughout most of the remaining paper. 

In particular, the separate differentiability estimates for the spatial gradient of the solution in space and time need to be interpolated in a suitable fashion, which requires the additional use of \emph{affine functions} in order to be able to utilize time differentiability estimates of order strictly smaller than one in our first-order setting. Such affine function techniques were not used in the elliptic nonlinear nonlocal case treated in \cite{DieKimLeeNow24}, however, they were successfully used in order to establish gradient potential estimates for linear elliptic nonlocal equations in \cite{KuuSimYan22}.

The different treatments of space and time eventually lead to the necessity of carrying out the concluding Campanato-type iteration with respect to a more complicated nonlocal excess functional than in the previous literature. Indeed, while in the local parabolic case simply a local excess functional of the type \eqref{loc.exc} was used in \cite{DM2} and in the nonlocal elliptic setting a stationary version of the nonlocal excess functional \eqref{not.exc} was utilized, we are forced to work with the following modified parabolic nonlocal excess functional
\begin{equation} \label{modex}
	{E}(u,\nabla; Q^s_R(z_0))\coloneqq E^p(\nabla u;Q^s_R(z_0))+\overline{E}(u;Q^s_R(z_0)),
\end{equation}
where $p \in \left ( 1, \frac{n+2s}{n+1}\right )$ is fixed and for some sufficiently large $q>1$,
\begin{align*}
	\overline{E}(u;Q^s_R(z_0))&\coloneqq\left(\dashint_{I_R(t_0)}\left(\dashint_{B_R(x_0)}\frac{|u-(\nabla u)_{Q^s_R(z_0)}\cdot (x-x_0)- (u)_{B_R(x_0)}(t)|}{R}\,dx\right)^q\,dt\right)^{\frac1q}\\
	&\,+\left(\dashint_{I_R(t_0)}\mathrm{Tail}\left(\frac{u-(\nabla u)_{Q^s_R(z_0)}\cdot (y-x_0)- (u)_{B_R(x_0)}(t)}R;B_R(x_0)\right)^q\,dt\right)^{\frac1q}.
\end{align*}

Here the first term on the right-hand side of \eqref{modex} originates from our differentiability estimates in space analogous to the elliptic setting, while the additional second term originates from our differentiability estimates in time involving affine functions. In view of gradient estimates in space and time for the corresponding homogeneous problem and suitable first-order comparison estimates that follow by interpolation, we are then able to prove excess decay estimates in the presence of measure data, which take the form
\begin{equation}\label{exc.ineqintro}
	\begin{aligned}
		{E}(u,\nabla;Q^s_{\rho R}(z_0))
		&\lesssim \rho^{\alpha}{E}(u,\nabla;Q^s_{R}(z_0))\\
		&\quad +\rho^{-(n+2s+1)}\left(\frac{|\mu|(Q^s_R(z_0))}{R^{n+1}}\right)^{1-\theta}{E}(u,\nabla;Q^s_{R}(z_0))^\theta\\
		&\quad+\rho^{-(n+2s+1)}\frac{|\mu|(Q^s_R(z_0))}{R^{n+1}}
	\end{aligned}
\end{equation}
for any $\rho \in (0,1]$ and some exponents $\alpha,\theta \in(0,1)$.

These excess decay estimates can then be iterated in a similar way as in the local parabolic setting treated in \cite{DM2} in order to obtain our pointwise gradient potential estimates in bounded domains given by Theorem \ref{thm.pt}. Finally, the corresponding gradient potential estimates on the whole space given by Theorem \ref{thm:gradpotwhole} then simply follow by letting $R \to \infty$ in Theorem \ref{thm.pt}.

\subsection{Outline}
The paper is structured as follows. In Section \ref{sec2}, we gather some basic notation as well as some definitions, embeddings and functional inequalities that will be used frequently throughout the paper.
In Section \ref{sec4} we then proceed to establish our localization lemma for parabolic nonlinear nonlocal equations and apply it in order to deduce our gradient H\"older regularity result in the homogeneous case given by Theorem \ref{thm.hol} along with some further useful decay estimates that turn out to be crucial in the proof of our gradient potential estimates. In Section \ref{sec3}, we establish comparison estimates that enable us to transfer information from the homogeneous case to the case of general measure data. As a first application, we then use these comparison estimates to establish Theorem \ref{thm: existence}, that is, the existence of SOLA to parabolic nonlinear nonlocal initial boundary value problems.
In Section \ref{sec5} we combine the results obtained in the homogeneous setting in Section \ref{sec4} with the comparison estimates obtained in Section \ref{sec3} to prove higher differentiability results in fractional Sobolev spaces under general measure data as well as suitable decay estimates for homogeneous problems and first-order comparison estimates. Finally, in Section \ref{sec6} we then utilize the estimates inferred in Section \ref{sec5} to obtain our pointwise gradient potential estimates and the associated fine regularity results.

\section{Preliminaries} \label{sec2}
\subsection{Some notation}
First of all, throughout this paper by $c$ we denote general positive
constants which could vary line by line. In addition, we use a parentheses to highlight relevant dependencies on parameters, i.e., $c=c(n,s,\Lambda)$ indicates that the constant $c$ depends only on $n,s$ and $\Lambda$.

For $U\subset\setR^n$, we define the indicator function of $U$ as
\begin{align*}
	\bfchi_{U}(x):=
	\begin{cases}
		1\quad\text{if }x\in U\\
		0\quad\text{if }x\in \setR^n\setminus U.
	\end{cases}
\end{align*}
Given a measurable function $g:\setR^{n+1}\rightarrow\setR$, we use the notation
\begin{align*}
	g_{\pm}(x):= \max\{\pm g(x),0\}.
\end{align*}
If $g$ is integrable over a measurable set $U\subset\setR^{n+1}$ with positive measure, i.e., $0<|U|<\infty$, then we denote by the integral average of $g$ over $U$ 
\begin{align*}
	(g)_U:=\dashint_{U}g\,dx=\dfrac{1}{|U|}\int_{U}g\,dx.
\end{align*}

In addition, given a signed Radon measure $\mu$ on $\mathbb{R}^{n+1}$, as usual we define the variation of $\mu$ as the measure defined by
$$ |\mu|(E):=\mu^+(E) + \mu^-(E), \quad E \subset \mathbb{R}^{n+1} \text{ measurable},$$
where $\mu^+$ and $\mu^-$ are the positive and negative parts of $\mu$, respectively. In the case when $|\mu|(\mathbb{R}^{n+1})<\infty$, then we say that $\mu$ has finite total mass.

Finally, given a domain $\Omega \subset \mathbb{R}^{n+1}$, throughout the paper we conceptualize functions $g \in L^1(\Omega)$ as signed Radon measures of class $\mathcal{M}(\mathbb{R}^{n+1})$ by extending $g$ by $0$ to $\mathbb{R}^{n+1}$ if necessary and denoting
$$ g(E):=\int_E g dx, \quad  E \subset \mathbb{R}^{n+1} \text{ measurable}.$$
Note that in this case for any measurable set $E \subset \mathbb{R}^{n+1}$, we have $$|g|(E)=\int_E |g| dx.$$

Next, we turn to some geometric notation and function spaces which will be used in this paper. First of all, occasionally we are going to write the spatial variable as $x$, the time variable as $t$, the space-time variable as $z$ and its variants, i.e., $z=(x,t),z_0=(x_0,t_0),z_1=(x_1,t_1),z_2=(x_2,t_2)\in\setR^n\times\setR=\setR^{n+1}$ and so on. Also, as indicated in the introduction, for any radius $R>0$, any $z_0=(x_0,t_0) \in \mathbb{R}^{n+1}$ and any $s \in (0,1]$, let us define the parabolic cylinder
\begin{equation*}
    Q^s_R(z_0)\coloneqq B_R(x_0)\times I^s_R(t_0),
\end{equation*}
where we write
\begin{equation*}
      I^s_R(t_0)\coloneqq(t_0-R^{2s},t_0).
\end{equation*}

\subsection{Function spaces}
In this subsection, let us define various function spaces tailored to our parabolic fractional setting. First of all, let us define a notion of parabolic H\"older spaces that is appropriate for our setting.
\begin{definition}
    For any $\alpha\in(0,1]$, we say that $g=g(x,t)\in C^{0,\alpha}(Q^s_R(z_0))$ if $g$ is continuous in $Q^s_R(z_0)$ with 
    \begin{align*}
        \sup_{z_1,z_2\in Q^s_R(z_0)}\frac{|g(z_1)-g(z_2)|}{\left(|x_1-x_2|+|t_1-t_2|^{\frac1{2s}}\right)^{\alpha}}<\infty.
    \end{align*}
\end{definition}

Next, we introduce fractional Sobolev spaces suitable for our parabolic setting, which are defined via the Bochner integral.
Let us denote by $X$ a Banach space and $I$ a time interval.
We say that if $g\in L^q(I;X)$ for some $q\geq1$, then 
\begin{align*}
    \|g\|_{L^q(I;X)}\coloneqq\left(\int_{I}\|g(\cdot,t)\|^q_{X}\,dt\right)^{\frac1q}<\infty
\end{align*}
and if $g\in W^{\sigma,q}(I;X)$ for some $\sigma\in(0,1)$, then
\begin{align}\label{def.vvsp}
    \|g\|_{W^{\sigma,q}(I;X)}\coloneqq\|g\|_{L^q(I;X)}+\left(\int_{I}\int_{I}\frac{\|g(\cdot,t)-g(\cdot,\tau)\|_X^q}{|t-\tau|^{1+\sigma q}}\,dt\,d\tau\right)^{\frac1q}<\infty.
\end{align}
In addition, we denote
\begin{align*}
    [g]_{L^q(I;W^{\sigma,q}(B))}\coloneqq \left(\int_{I}\int_B\int_B\frac{|g(x,t)-g(y,t)|^q}{|x-y|^{n+\sigma q}}\,dx\,dy\,dt\right)^{\frac1q}
\end{align*}
and
\begin{align*}
    [g]_{W^{\sigma,q}(I;L^{q}(B))}\coloneqq \left(\int_{I}\int_I\int_B\frac{|g(x,t)-g(x,\tau)|^q}{|t-\tau|^{1+\sigma q}}\,dx\,dt\,d\tau\right)^{\frac1q}.
\end{align*}
Observe that 
\begin{align*}
    &\|g\|_{W^{\sigma,q}(I;W^{\gamma,q}(B))}\\
    &\eqsim\left(\int_I\int_I\int_B\int_B\frac{|(g(x,t)-g(x,\tau))-(g(y,t)-g(y,\tau))|^q}{|x-y|^{n+\gamma q}|t-\tau|^{1+\sigma q}}\,dx\,dy\,dt\,d\tau\right)^{\frac1q}\\
    &\quad+ [g]_{L^q(I;W^{\sigma,q}(B))}+[g]_{W^{\sigma,q}(I;L^{q}(B))}+\|g\|_{L^q(I;L^q(B))}
\end{align*}
and 
\begin{align*}
    \|g\|_{W^{\sigma,q}(I;W^{1+\gamma,q}(B))}\eqsim\|\nabla g\|_{W^{\sigma,q}(I;W^{\gamma,q}(B))}+\|g\|_{W^{\sigma,q}(I;W^{\gamma,q}(B))}.
\end{align*}
We refer \cite{Sim90, Sim87} for various embedding and interpolation results concerning these function spaces. 
\subsection{Parabolic fractional Sobolev spaces and difference quotients}

Let us denote
\begin{equation*}
    \delta_h g(x,t)\coloneqq g(x+h,t)-g(x,t),\quad\delta_h^2g(x,t)\coloneqq\delta_h(\delta_h g)(x,t)
\end{equation*}
for any $h\in \bbR^n\setminus \{0\}$, and 
\begin{equation}\label{defn.diff.t}
    {\delta}^t_hg(x,t)\coloneqq g(x,t-h^{2s})-g(x,t)
\end{equation}
for any $h>0$.
We now provide various embedding inequalities in terms of difference quotients.
\begin{lemma}\label{lem.firemb}
    Let $q\in[1,\infty)$ and $\gamma\in(0,1)$. If $g\in L^q(Q^s_{R+h_0}(z_0))$ satisfies
    \begin{equation*}
        \sup_{0<|h|<h_0}\left\|\frac{\delta_hg}{|h|^{\gamma}}\right\|_{L^{q}(Q^s_{R}(z_0))}< \infty
    \end{equation*}
    for some constant $h_0\in(0,R/4)$, then for every $\widetilde{\gamma}\in(0,\gamma)$, 
    \begin{equation*}
        g\in L^q(I^s_{R/2}(t_0);W^{\widetilde{\gamma},q}({B_{R/2}(x_0)})).
    \end{equation*} Moreover, we have
    \begin{align*}
        [g]^q_{L^q(I^s_{R/2}(t_0);W^{\widetilde{\gamma},q}({B_{R/2}(x_0)}))}&\leq \frac{ch_0^{q(\gamma-\widetilde{\gamma})}}{(\gamma-\widetilde{\gamma})^q}\sup_{0<|h|<h_0}\left\|\frac{\delta_h g}{|h|^\gamma}\right\|^q_{L^q(Q^s_R(z_0))}\\
        &\quad+c\left(\frac{h_0^{q(1-\widetilde{\gamma})}}{R^q(\gamma-\widetilde{\gamma})^q}+ \frac{h_0^{-\widetilde{\gamma} q}}{\widetilde{\gamma}}\right)\|g-k\|^q_{L^{q}(Q^s_R(z_0))}
    \end{align*}
    for any $k\in\bbR$ with some constant $c=c(n,q)$.
\end{lemma}
\begin{proof}
    By \cite[Lemma 2.3]{DieKimLeeNow24d}, we observe that
    \begin{align*}
        [g(\cdot,t)]^q_{W^{\widetilde{\gamma},q}(B_{R/2}(x_0))}&\leq \frac{ch_0^{q(\gamma-\widetilde{\gamma})}}{(\gamma-\widetilde{\gamma})^q}\sup_{0<|h|<h_0}\left\|\frac{\delta_h g(\cdot,t)}{|h|^\gamma}\right\|^q_{L^q(B_R(x_0))}\\
        &\quad+c\left(\frac{h_0^{q(1-\widetilde{\gamma})}}{R^q(\gamma-\widetilde{\gamma})^q}+ \frac{h_0^{-\widetilde{\gamma} q}}{\widetilde{\gamma}}\right)\|(g-k)(\cdot,t)\|^q_{L^{q}(B_R(x_0))}
    \end{align*}
    holds a.e. $t\in I^s_R(t_0)$, where $c=c(n,q)$. By integrating both sides of the above inequality with respect to the time variable, we deduce the desired estimate.
\end{proof}
\begin{lemma}\label{lem.secemb}
    Let $q\in[1,\infty)$ and $g\in L^q(Q^s_{R+6h_0}(z_0))$ for some $h_0\in(0,R/4)$. If 
    \begin{align*}
        h_0^{-q}\sup_{0<|h|<h_0}\int_{Q^s_{R+4h_0}(z_0)}\frac{|\delta_h g|^q}{|h|^{q\gamma}}\,dz+\sup_{0<|h|<h_0}\int_{Q^s_{R+4h_0}(z_0)}\frac{|\delta^2_h g|^q}{|h|^{q+q\gamma}}\,dz\leq M^q
    \end{align*}
    for some constants $M>0$ and $\gamma\in(0,1)$, then we have $g\in L^q(I^s_R(t_0);W^{1,q}(B_R(x_0)))$ with the estimate 
    \begin{align*}
        \|\nabla g\|_{L^q(Q^s_R(z_0))}\leq cM+c(h_0^{-(1+\gamma)}+1)\|g\|_{L^q(Q^s_{R+4h_0}(z_0))}
    \end{align*}
    for some constant $c=c(n,q,\gamma)$.
\end{lemma}
\begin{proof}
    We note that the corresponding elliptic version of this lemma with $q=1$ is given in \cite[Lemma 2.10]{DieKimLeeNow24}. In addition, a careful inspection of Lemma 2.10 in \cite{DieKimLeeNow24} yields the corresponding elliptic version of this lemma for every $q\geq1$.
    The desired estimate in the parabolic setting can now be deduced from its elliptic counterpart in a similar fashion as in the proof of Lemma \ref{lem.firemb}.
\end{proof}
In addition, the following embedding result is a direct consequence of \cite[Lemma 2.9]{DieKimLeeNow24}.
\begin{lemma}\label{lem.thiemb}
    Let $q\in[1,\infty)$ and $g\in L^{q}(I^s_{R+6h_0}(t_0);W^{1,q}(B_{R+6h_0}(x_0))$ with $h_0\in(0,R/4)$. If 
    \begin{align*}
    \sup_{0<|h|<h_0}\int_{Q^s_{R+4h_0}(z_0)}\frac{|\delta_h^2g|^q}{|h|^{q(1+\gamma)}}\,dz\leq M^q
    \end{align*}
    for some constants $M>0$ and $\gamma\in(0,1)$, then we have 
    \begin{align*}
        &[\nabla g]^q_{L^q(I^s_R(t_0);W^{\widetilde{\gamma},q}(B_R(x_0)))}\\
        &\leq \frac{ch_0^{q(\gamma-\widetilde{\gamma})}M^q}{(\gamma-\widetilde{\gamma})\gamma^q(1-\gamma)^q}+\frac{ch_0^{q(\gamma-\widetilde{\gamma})}}{(\gamma-\widetilde{\gamma})\gamma^q(1-\gamma)^q}\frac{(R+4h_0)^{q+n}}{h_0^{q(1+\gamma)}}E^q_{\mathrm{loc}}(\nabla g;Q^s_{R+4h_0}(z_0))^q
    \end{align*}
    for some constant $c=c(n,q)$, where $\widetilde{\gamma}\in(0,\gamma)$.
\end{lemma}
Next, we establish an embedding result in terms of the function $\delta^t_h g$. Since $h$ is only allowed to be a positive number, in this case we cannot directly use the results given in the above lemmas. However, with aid of the even extension with respect to the time direction, we are able to prove the following.
\begin{lemma}\label{lem.fouemb}
    Let $g\in L^q(Q^s_{R+h_0}(z_0))$ satisfy
    \begin{align*}
        \sup_{0<h<h_0}\left\|\frac{\delta^t_h g}{h^{\gamma}}\right\|_{L^{q}(Q^s_{R}(z_0))}\leq M
    \end{align*}
    for some constant $h_0\in(0,s^{\frac1{2s}}R/100)$.
    Then we have for any $\widetilde{\gamma}\in(0,\gamma)$,
    \begin{align*}
        [g]^q_{W^{\widetilde{\gamma}/2s,q}(I^s_{R/2}(t_0);L^{q}({B_{R/2}(x_0)}))}&\leq \frac{ch_0^{q({\gamma}-\widetilde{\gamma})}}{(\gamma-\widetilde{\gamma})^q}M^q\\
       &\quad+c\left(\frac{h_0^{2sq(1-\widetilde{\gamma}/2s)}}{R^{2sq}(\gamma-\widetilde{\gamma})^q}+\frac{h_0^{-\widetilde{\gamma}q}}{\widetilde{\gamma}}\right)\|g-k\|^q_{L^q(Q^s_R(z_0))}
    \end{align*}
    for any $k\in\bbR$ with some constant $c=c(n,s_0,q)$. 
\end{lemma}
\begin{proof}We may assume $z_0=0$.
Let us fix $s_0\in(0,1)$ and $s\in[s_0,1)$. Define
    \begin{align*}
        G(x,t)\coloneqq\begin{cases}
            g(x,t)&\quad\text{if }t\in(-(R+h_0)^{2s},0],\\
            g(x,-t)&\quad\text{if }t\in [0,(R+h_0)^{2s}).
        \end{cases}
    \end{align*}
    We are going to prove that
    \begin{align}\label{emb.time.goalf}
        \sup_{0<|h|<h_0^{2s}/3}\left\|\frac{G(x,t-h)-G(x,t)}{|h|^{\gamma/2s}}\right\|_{L^q(B_R\times (-R^{2s},R^{2s}))}\leq cM
    \end{align}
    for some constant $c=c(n,s_0,\gamma)$.
    We first observe that
   \begin{align}\label{fir.ob}
        \sup_{0<h<h_0^{2s}/3}\left\|\frac{G(x,t-h)-G(x,t)}{|h|^{\gamma/2s}}\right\|_{L^q(Q^s_R)}\leq \sup_{0<h<h_0}\left\|\frac{\delta_h^{t}g}{|h|^{\gamma}}\right\|_{L^q(Q^s_R)},
    \end{align}
    and note that 
    \begin{align*}
        &\sup_{-h_0^{2s}/3<h<0}\left\|\frac{G(x,t-h)-G(x,t)}{|h|^{\gamma/2s}}\right\|_{L^q(Q^s_R)}\\
        &\leq\sup_{t<h<0}\left\|\frac{G(x,t-h)-G(x,t)}{|h|^{\gamma/2s}}\right\|_{L^q(Q^s_R)}\\
        &\quad+\sup_{-h_0^{2s}/3<h<t}\left\|\frac{G(x,t-h)-G(x,t)}{|h|^{\gamma/2s}}\right\|_{L^q(Q^s_R)}\eqqcolon J_1+J_2.
    \end{align*}
    Since $t-h\in (-R^{2s},0)$ if $t<h<0$, it follows that
    \begin{align*}
        J_1\leq \sup_{0<h<h_0}\left\|\frac{\delta_h^{t}g}{|h|^{\gamma}}\right\|_{L^q(Q^s_R)}.
    \end{align*}
    On the other hand, we observe 
    \begin{align*}
        |2t-h|\leq 3|h|\quad\text{if }t\in(-R^{2s},0]\quad\text{and}\quad h\in(-h_0^{2s}/3,t).
    \end{align*}
    Therefore, we deduce 
    \begin{align*}
        J_2&=\sup_{-h_0^{2s}/3<h<t}\left\|\frac{g(x,h-t)-g(x,t)}{|h|^{\gamma/2s}}\right\|_{L^q(Q^s_R)}\\
        &\leq c\sup_{0<2t-h<h_0^{2s}}\left\|\frac{g(x,t)-g(x,t-(2t-h))}{|2t-h|^{\gamma/2s}}\right\|_{L^q(Q^s_R)}
        \leq c\sup_{0<h<h_0}\left\|\frac{\delta_h^{t}g}{|h|^{\gamma}}\right\|_{L^q(Q^s_R)}
    \end{align*}
    for some constant $c=c(n,s_0)$. Combining \eqref{fir.ob} and the estimates $J_1$ and $J_2$, we obtain 
    \begin{align*}
       \sup_{0<|h|<h_0^{2s}/3}\left\|\frac{G(x,t-h)-G(x,t)}{|h|^{\gamma/2s}}\right\|_{L^q(B_R\times (-R^{2s},0))}\leq cM.
   \end{align*}
   Similarly, we also obtain 
   \begin{align*}
       \sup_{0<|h|<h_0^{2s}/3}\left\|\frac{G(x,t-h)-G(x,t)}{|h|^{\gamma/2s}}\right\|_{L^q(B_R\times (0,R^{2s}))}\leq cM.
   \end{align*}
   The above estimates imply \eqref{emb.time.goalf}. By following the same lines as in the proof \cite[Lemma 2.3]{DieKimLeeNow24d} with $n,R,p$ and $\psi(x)$ replaced by $1,R^{2s},q$ and $\eta(t)$, respectively, where $\eta(t)\in C_{c}^{\infty}(-(3R/4)^{2s},(3R/4)^{2s})$ with $\eta\equiv 1$ on $(-(R/2)^{2s},(R/2)^{2s})$, we observe 
   \begin{align*}
       &[G(x,\cdot)]^q_{W^{\widetilde{\gamma}/2s,q}(-(R/2)^{2s},(R/2)^{2s})}\\
       &\leq \frac{ch_0^{{{2s}}q({\gamma}-\widetilde{\gamma})}}{(\gamma-\widetilde{\gamma})^q}\sup_{0<h<h_0^{2s}/3}\left\|\frac{\delta_h G(x,\cdot)}{|h|^{\gamma/2s}}\right\|^q_{L^q(-R^{2s},R^{2s})}\\
       &\quad+c\left(\frac{ch_0^{2sq(1-\widetilde{\gamma}/2s)}}{R^{2sq}(\gamma-\widetilde{\gamma})^q}+\frac{h_0^{-\widetilde{\gamma}q}}{\widetilde{\gamma}}\right)\|G(x,\cdot)-k\|^q_{L^q(-R^{2s},R^{2s})}
   \end{align*}
   for some constant $c=c(n,s_0,q)$. By integrating both sides of the above inequality on the spatial direction and by using \eqref{emb.time.goalf}, we get 
   \begin{align*}
       [G]^q_{W^{\widetilde{\gamma}/2s,q}(I^s_{R/2};L^q(B_{R/2}))}&\leq \frac{ch_0^{q({\gamma}-\widetilde{\gamma})}}{(\gamma-\widetilde{\gamma})^q}M^q\\
       &\quad+c\left(\frac{h_0^{2sq(1-\widetilde{\gamma}/2s)}}{R^{2sq}(\gamma-\widetilde{\gamma})^q}+\frac{h_0^{-\widetilde{\gamma}q}}{\widetilde{\gamma}}\right)\|g-k\|^q_{L^q(Q^s_R)}
   \end{align*}
   for some constant $c=c(n,s_0,q)$. Since $G(\cdot,t)=g(\cdot,t)$ on $I^s_R$, the desired result follows from the above estimate.
\end{proof}

\subsection{Fractional Poincar\'e-type inequalities}
The following two lemmas yield certain Poincar\'e-type inequalities, which can be considered to be replacements for standard fractional (Sobolev-)Poincar\'e inequalities in our parabolic setting.
\begin{lemma}\label{lem.spwa}
    Let $g\in W^{\gamma/2s,q}(I^s_R;L^q(B_R))\cap L^q(I^s_R;W^{\gamma,q}(B_R))$ for some constant $\gamma\in(0,1)$ and $q\in[1,\infty)$. Then we have 
    \begin{align*}
        \left(\dashint_{Q^s_R}|g-(g)_{Q^s_R}|^q\,dz\right)^{\frac1q}\leq cR^{-\frac{n+2s}{q}+\gamma }\left([g]_{W^{\gamma/2s,q}(I^s_R;L^q(B_R))}+[g]_{L^q(I^s_R;W^{\gamma,q}(B_R))}\right)
    \end{align*}
    for some constant $c=c(n)$.
\end{lemma}
\begin{proof}
    We first note that
    \begin{align*}
        \left(\dashint_{Q^s_R}|g-(g)_{Q^s_R}|^q\,dz\right)^{\frac1q}&\leq \left(\dashint_{Q^s_R}|g-(g)_{B_R}(t)|^q\,dz\right)^{\frac1q}\\
        &\quad+\left(\dashint_{Q^s_R}|(g)_{B_R}(t)-(g)_{Q^s_R}|^q\,dz\right)^{\frac1q}\eqqcolon J_1+J_2.
    \end{align*}
    Next, we estimate $J_1$ as 
    \begin{align*}
        J_1\leq \left(\dashint_{Q^s_R}\dashint_{B_R}|g(x,t)-g(y,t)|^q\,dy\,dz\right)^{\frac1q}\leq cR^{\gamma}\left(\dashint_{Q^s_R}\dashint_{B_R}\frac{|g(x,t)-g(y,t)|^q}{|x-y|^{n+\gamma q}}\,dy\,dz\right)^{\frac1q}
    \end{align*}
    for some constant $c=c(n)$, where we have used H\"older's inequality and a few algebraic inequalities. Similarly, we deduce
    \begin{align*}
        J_2\leq \left(\dashint_{Q^s_R}\dashint_{I^s_R}|g(x,t)-g(x,\tau)|^q\,d\tau\,dz\right)^{\frac1q} \leq cR^{\gamma}\left(\dashint_{Q^s_R}\int_{I^s_R}\frac{|g(x,t)-g(x,\tau)|^q}{|t-\tau|^{1+\gamma q/2s}}\,d\tau\,dz\right)^{\frac1q}
    \end{align*}
    for some constant $c=c(n)$. Combining the above estimates for $J_1$ and $J_2$ yields the desired estimate.
\end{proof}

\begin{lemma}\label{lem.inter}
    Fix $p\in [1,2]$ and $\sigma\in(0,1)$. If $g\in L^p(I^s_R;W^{\sigma,p}(B_{2R}))\cap L^\infty(I^s_R;L^1(B_{2R}))$ with $g(x,\cdot)\equiv0$ on $B_{2R}\setminus B_R$, then we have 
\begin{align*}
    \|g\|_{L^{\frac{p(n+\sigma)}{n}}(B_{2R}\times I^s_R)}&\leq c \left((1-\sigma)\int_{I^s_R}\int_{B_{2R}}\int_{B_{2R}}\frac{|g(x,t)-g(y,t)|^{p}}{|x-y|^{n+ \sigma p}}\,dy\,dz\right)^{\frac{n}{p(n+\sigma)}}\\
    &\qquad\times\left(\sup_{t\in I^s_{R}}\|g(\cdot,t)\|_{L^{1}(B_{2R})}\right)^{\frac{\sigma}{n+\sigma}}
\end{align*}
for some constant $c=c(n,\sigma,p)$, where the constant $c$ depends only on $n,\sigma_0$ and $p$  whenever $\sigma\in[\sigma_0,1)$.
\end{lemma}
\begin{proof}
Using an interpolation argument, we get
\begin{align*}
    \|g\|_{L^{\frac{p(n+\sigma)}{n}}(B_{2R}\times I^s_R)}&\leq c\|g\|_{L^{p}\left(I^s_{R};L^{\frac{np}{n-\sigma p}}(B_{2R})\right)}^{\frac{n}{n+\sigma}}\|g\|_{L^{\infty}(I^s_{R};L^{1}(B_{2R}))}^{\frac{\sigma}{n+\sigma}}.
\end{align*}
We note from \cite[Corollary 4.9]{Coz17} that
\begin{align*}
    \|g\|^2_{L^{p}\left(I^s_{R};L^{\frac{pn}{n-p\sigma}}(B_{2R})\right)}\leq c(1-\sigma)\int_{I^s_{R}}\int_{B_{2R}}\int_{B_{2R}}\frac{|g(x,t)-g(y,t)|^{p}}{|x-y|^{n+\sigma p}}\,dy\,dz
\end{align*}
for some constant $c=c(n,s_0,p)$. The desired result follows by combining the above two inequalities.
\end{proof}

\subsection{Embeddings and interpolation in parabolic fractional Sobolev spaces}

We now prove a simple embedding result in the parabolic fractional Sobolev spaces defined in \eqref{def.vvsp}. 
\begin{lemma}\label{lem.holspI}
	Assume $s\in(1/2,1)$. Let $g\in C^{\gamma}(Q^s_1)$ for some $\gamma\in(0,1)$. Then 
	\begin{align*}
		g\in W^{\gamma/(8s),q}({I^s_1};W^{\gamma/4,q}(B_1))
	\end{align*} with the estimate
	\begin{align*}
		\|g\|_{W^{\gamma/(8s),q}({I^s_1};W^{\gamma/4,q}(B_1))}\leq cR^{\gamma}[g]_{C^{\gamma}(Q^s_1)}+c\|g\|_{L^q(Q^s_1)}
	\end{align*}
	for some constant $c=c(n,s_0,q,\gamma)$.
\end{lemma}
\begin{proof}
	We observe that for any $x,y\in B_1$ and $t,\tau\in I^s_1$,
	\begin{align*}
		&\frac{|(g(x,t)-g(x,\tau))-(g(y,t)-g(y,\tau))|^q}{|x-y|^{n+\gamma q/4}|t-\tau|^{1+\gamma q/(8s)}}\\
		&\leq \frac{(|g(x,t)-g(x,\tau)|+|g(y,t)-g(y,\tau)|)^{q/2}}{|t-\tau|^{1+\gamma q/(8s)}}
	\frac{(|g(x,t)-g(y,t)|+|g(x,\tau)-g(y,\tau)|)^{q/2}}{|x-y|^{n+\gamma q/4}}\\
		&\leq \frac{[g]_{C^{\gamma}(Q^s_{1})}^{q/2}}{|t-\tau|^{1-\gamma q/(8s)}}\frac{[g]_{C^{\gamma}(Q^s_{1})}^{q/2}}{|x-y|^{n-\gamma q/4}}.
	\end{align*}
	Using this, we obtain
	\begin{align*}
		&\int_{I^s_1}\int_{I^s_1}\int_{B_1}\int_{B_1}\frac{|(g(x,t)-g(x,\tau))-(g(y,t)-g(y,\tau))|^q}{|x-y|^{n+\gamma q/4}|t-\tau|^{1+\gamma q/(8s)}}\,dx\,dy\,dt\,d\tau\\  
		&\leq \int_{I^s_1}\int_{I^s_1}\int_{B_1}\int_{B_1}\frac{[g]_{C^{\gamma}(Q^s_{1})}^{q}}{|t-\tau|^{1-\gamma q/(8s)}|x-y|^{n-\gamma q/4}}\,dx\,dy\,dt\,d\tau\leq c[g]_{C^{\gamma}(Q^s_{1})}^{q}
	\end{align*}
	for some constant $c=c(n,s_0,q,\gamma)$. Similarly, we have 
	\begin{align*}
		[g]_{L^q(I^s_1;W^{\gamma/4,q}(B_1))}+[g]_{W^{\gamma/(8s),q}(I^s_1;L^q(B_1))}\leq c[g]_{C^{\gamma}(Q^s_{1})}^{q}.
	\end{align*}
	Therefore, combining all the estimates yields the desired result.
\end{proof}
We will also utilize the following interpolation lemma between parabolic fractional Sobolev spaces.
\begin{lemma}\label{lem.interpovec}
	Let $g\in W^{s_0,q}(I^s_1;W^{\overline{s}_0,q}(B_1))\cap W^{s_1,q}(I^s_1;W^{1+\overline{s}_1,q}(B_1))$
	for some constants $s_i,\overline{s}_i\in(0,1)$, where $i=1,2$. Let us fix $\Theta\in(0,1)$ and choose
	\begin{align*}
		s_{\Theta}=s_0\Theta+s_1(1-\Theta)\quad\text{and}\quad \overline{s}_{\Theta}=\overline{s}_0\Theta+(1+\overline{s}_1)(1-\Theta).
	\end{align*}
	Then if $\overline{s}_{\Theta}\neq 1$, then we have
	\begin{align}\label{ineq.inter}
		\|g\|_{W^{s_{\Theta},q}(I^s_1;W^{\overline{s}_{\Theta},q}(B_1))}\leq c\|g\|^{\Theta}_{W^{s_0,q}(I^s_1;W^{\overline{s}_0,q}(B_1))}\|g\|^{1-\Theta}_{W^{s_1,q}(I^s_1;W^{1+\overline{s}_1,q}(B_1))}
	\end{align}
	for some constant $c=c(n,s_0,\overline{s}_0,s_1,\overline{s}_1,\Theta,q)$.
\end{lemma}
\begin{proof}
	By the standard extension lemma, it suffices to prove \eqref{ineq.inter} with $B_1$ replaced by $\bbR^n$.
	We first note from \cite[Proposition in 2.5.1]{RunSic96} that
	\begin{align*}
		[W^{\overline{s}_0,q}(\bbR^n),W^{\overline{s}_1,q}(\bbR^n)]_{\Theta,q}=W^{\overline{s}_{\Theta},q}(\bbR^n).
	\end{align*}
	Using this and \cite[Lemma 7]{Sim87}, we have 
	\begin{align*}
		\|g\|_{W^{s_{\Theta},q}(I^s_1;W^{\overline{s}_{\Theta},q}(\bbR^n))}\leq c\|g\|^{\Theta}_{W^{s_0,q}(I^s_1;W^{\overline{s}_0,q}(\bbR^n))}\|g\|^{1-\Theta}_{W^{s_1,q}(I^s_1;W^{1+\overline{s}_1,q}(\bbR^n))}
	\end{align*}
	for some constant $c=c(n,s_0,s_1,\overline{s}_0,\overline{s}_1,\Theta)$, which completes the proof.
\end{proof}

\subsection{A parabolic Campanato-type embedding}
The following Campanato-type embedding can be deduced in a standard way, for instance by following the proof of \cite[Lemma 2.1]{ByuKimKim23}.
\begin{lemma}\label{lem.camp}
    Let $g\in L^1(Q^s_{2R}(z_0))$. Suppose that there exist constants $M>0$ and $\gamma\in(0,1)$ such that
    \begin{align*}
        \dashint_{Q^s_{r}(z_1)}|g-(g)_{Q^s_r(z_1)}|\,dz\leq M\rho^{\gamma}
    \end{align*}
    whenever $Q^s_r(z_1)\Subset Q^s_{R}(z_0)$.
    Then $g\in C^{\gamma}(Q^s_{R/2}(z_0))$ with the estimate 
    \begin{align}\label{ineq.camp1}
        [g]_{C^{0,\gamma}(Q^s_{R/4}(z_0))}\leq cR^{-\gamma}\left[M+\dashint_{Q^s_{3R/4}}|g-(g)_{Q^s_{3R/4}(z_0)}|\,dz\right]
    \end{align}
    for some constant $c=c(n,s_0,\gamma)$.
\end{lemma}

\subsection{A Riesz potential estimate via Lorentz spaces}
We next prove the following Riesz potential estimates in terms of Lorentz spaces that are crucial to deduce local boundedness and gradient continuity later in the paper. The proof follows the arguments given in \cite{Cianchi}.
\begin{lemma}\label{lem.lorentz}
   Let $s>1/2$ and let $g\in L^{\frac{n+2s}{2s-1},1}(Q^s_{2R})$. Then we have 
    \begin{align*}
        I_{2s-1}^{|g|}(z_0,R)\leq c \|g\|_{L^{\frac{n+2s}{2s-1},1}(Q^s_R(z_0))}
    \end{align*}
    for any $z_0\in Q^s_R$, where $c=c(n,s)$.
\end{lemma}
\begin{proof}
    We first note that $Q^s_R(z_0)\Subset Q^s_{2R}$, as $s>1/2$. Let us define 
    \begin{equation*}
        g^{*}(r)=\sup\{\tau\geq0\,:\,|\{z\in Q^s_R(z_0):|g(z)|>\tau\}>r|\},
    \end{equation*}
    which is called the decreasing rearrangement of $g$ on $Q^s_R(z_0)$. Next, we define $g^{**}(\tau)=\frac{1}{\tau}\int_{0}^{\tau}g^{*}(\xi)\,d\xi$ (see \cite{Cianchi} for more details of the function $g^*$).
    Then we note from \cite[Equation (2.19)]{Cianchi} that 
    \begin{align}\label{ineq.ris}
        \int_{0}^{R^{n+2s}}{g^{**}(r)}r^{\frac{2s-1}{n+2s}-1}\,dr\eqsim \|g\|_{L^{\frac{n+2s}{2s-1},1}(Q^s_R(z_0))}.
    \end{align}
    Therefore, using \cite[Equation (2.18)]{Cianchi} with $\sigma=1$ and $\rho=1$, and Fubini's theorem, we get
    \begin{equation}\label{ineq.ris2}
    \begin{aligned}
        I_{2s-1}^{|g|}(z_0,R)&=\int_{0}^R{|g|(Q^s_r(z_0))}{r^{-(n+1)}}\,dr\\
        &=\int_{0}^{R}r^{-(n+1)}\int_{0}^{r^{n+2s}}(g\chi_{Q^s_r(z_0)})^{*}(\tau)\,d\tau \,dr\\
        &=\int_{0}^{R}r^{2s-1}(g\chi_{Q^s_r(z_0)})^{**} (r^{n+2s})\,dr\\
        &\leq \int_{0}^{R}r^{2s-1}g^{**}(r^{n+2s})\,dr\\
        &=(n+2s)\int_{0}^{R^{n+2s}}{g^{**}(r)}r^{\frac{2s-1}{n+2s}-1}\,dr,
    \end{aligned}
    \end{equation}
    where we have also used $(g\chi_{Q^s_r(z_0)})^{*}(\tau)=0$ if $\tau>r^{n+2s}$ and the fact that 
    \begin{equation*}
        g_1^{**}(r)\leq g_2^{**}(r)\quad\text{if }|g_1(z)|\leq |g_2(z)|.
    \end{equation*}
    Plugging \eqref{ineq.ris} into the last line in \eqref{ineq.ris2}, we obtain the desired result.
\end{proof}

\subsection{A covering lemma}
Next, we prove the following simple covering lemma which will be a crucial tool to establish higher differentiability results for parabolic nonlinear nonlocal measure data problems.
\begin{lemma}\label{lem.cov}
	Let $R>0$ be fixed. Let us choose $r\in(0,s^{\frac1{2s}}R/2)$. Then there is a constant $c=c(n)$, finite index sets $\mathcal{I}$ and $\mathcal{J}$ and sequences $\{x_i\}_{i \in \mathcal{I}} \subset B_R$ and $\{t_j\}_{j \in \mathcal{J}} \subset I^s_R$ such that for any $k\in\setN$,
	\begin{equation}\label{cov.ineq1}
		B_{R}\subset \bigcup_{i\in \mathcal{I}}B_{r}(z_{i})\subset B_{2R}, \quad \sup_{x_i \in \mathbb{R}^{n}} \sum_{i\in \mathcal{I}}\bfchi_{B_{2^{k}r}(x_{i})}(x)\leq c2^{nk}, \quad |\mathcal{I}| \leq c \frac{R^{n}}{r^{n}}
	\end{equation}
 and 
 \begin{equation}\label{cov.ineq2}
		I^s_{R}\subset \bigcup_{j\in \mathcal{J}}I^s_{r}(t_{j})\subset I^s_{2R}, \quad \sup_{t_j \in \mathbb{R}} \sum_{j\in \mathcal{J}}\bfchi_{I^s_{2^{k}r}(t_{j})}(t)\leq c2^{2sk}, \quad |\mathcal{J}| \leq c \frac{R^{2s}}{r^{2s}}
	\end{equation}
 where we denote by $|\mathcal{I}|$ and $|\mathcal{J}|$ the number of elements in the sets $\mathcal{I}$ and $\mathcal{J}$, respectively.
\end{lemma}
\begin{proof}
For the proof of \eqref{cov.ineq1}, we refer \cite[Lemma 2.11]{DieKimLeeNow24}. Now, we are going to prove \eqref{cov.ineq2}.
	We note that there is a mutually disjoint covering $\{I^s_{r}(t_j)\}_{i\in \mathcal{J}}$ of $I^s_{R}$ such that
 \begin{equation*}
     t_j\in I^s_R\quad\text{and}\quad I^s_{r}(t_j)\subset I^s_{2R},
 \end{equation*} as $r<s^{\frac1{2s}}R/2$. Therefore, we have
	\begin{align}
		\label{el : ineq1.besi}
		|\mathcal{J}||I^s_{r}|=\sum_{j\in \mathcal{J}}|I^s_{r}|=\sum_{j\in \mathcal{J}}|I^s_{r}(t_j)|\leq |I^s_{2R}|
	\end{align}
 which gives the third inequality given in \eqref{cov.ineq2}.
	We are now in the position to prove 
	\begin{align}
		\label{el : ineq00.bdry}
		 \sup_{t \in \mathbb{R}}\sum_{j\in \mathcal{J}}\bfchi_{I^s_{2^{k}r}(t_{j})}(t)\leq 2^{2sk+1}.
	\end{align}
	Suppose there is a point $\overline{t}\in \bbR$ such that $\sum_{j\in \mathcal{J}}\bfchi_{I^s_{2^{k}r}(t_{j})}(\overline{t})> 2^{2sk+1}$.
	We now denote by $\mathcal{J}_{0}$ the set $\{j\in \mathcal{J}\,:\, \bfchi_{I^s_{2^{k}r}(t_j)}(\overline{t})=1\}$.
	Then we observe 
	\begin{equation*}
		\bigcup_{j\in\mathcal{J}_{0}} I^s_{2^{k}r}(t_{j})  \subset I^s_{2^{\frac1{2s}+k}r}(\overline{t}+(2^{k}r)^{2s}),
	\end{equation*}
	which implies 
	\begin{align*}
		(2^{2sk+1}+\!1)|I^s_{2^{k}r}|\leq\sum_{j\in \mathcal{J}_0}|I^s_{2^{k}r}(t_j)|\leq 2^{2sk}\sum_{j\in \mathcal{J}_{0}}|I^s_{r}(t_j)|&\leq 2^{2sk}|I^s_{2^{\frac1{2s}+k}r}(\overline{t}+(2^kr)^{2s})|\\
  &\leq 2^{2sk+1}|I^s_{2^kr}|,
	\end{align*}
	where for the third inequality we have used the fact that $\{I^s_{r}(t_{j})\}_{j\in \mathcal{J}}$ is a mutually disjoint set. This is a contradiction. Thus we show \eqref{cov.ineq2}, which completes the proof.
\end{proof}

\subsection{Parabolic tails and affine function}

We now conclude Section \ref{sec2} with a tail estimate that will be play a crucial role for obtaining our gradient potential estimates (see Lemma \ref{lem.exc.dec.fin} below).
\begin{lemma}\label{lem.tail1} Suppose $s>s_0\geq 1/2$.
   Let 
   \begin{align*}
       u\in L^q(I^s_{2^iR}(t_0);W^{1,q}(B_{2^iR}(x_0)))
   \end{align*}for some constant $q\in[1,\infty)$ and for some positive integer $i$. 
    Then we have
    \begin{align*}
      &\left(\dashint_{I^s_{R}(t_0)}\mathrm{Tail}\left(\frac{u-(\nabla u)_{Q^s_{R}(z_0)}\cdot (y-x_0)-(u)_{B_{R}(x_0)}(t)}{R};B_{R}(x_0)\right)^q\,dt\right)^{\frac1q}\\
      &\leq c\left(\dashint_{I^s_{R}(t_0)}\left(\sum_{j=0}^{i}2^{j(1-2s)}\dashint_{B_{2^{j}R}(x_0)}|\nabla u-(\nabla u)_{Q^s_{2^{j}R}(z_0)}|\,dx\right)^q\,dt\right)^{\frac1q}\\
      &\quad +c\sum_{j=0}^{i}2^{j(1-2s)}\dashint_{Q^s_{2^jR}(z_0)}|\nabla u-(\nabla u)_{Q^s_{2^{j}R}(z_0)}|\,dz\\
      &\quad+ 2^{i(1-2s)}\left(\dashint_{I^s_{R}(t_0)}\mathrm{Tail}\left(\frac{u-(\nabla u)_{Q^s_{2^iR}(z_0)}\cdot (y-x_0)-(u)_{B_{2^iR}(x_0)}(t)}{2^iR};B_{2^iR}(x_0)\right)^q\,dt\right)^{\frac1q}
    \end{align*}
    for some constant $c=c(n,s_0,q)$.
\end{lemma}
\begin{proof}
    We may assume $R=1$ and $z_0=0$ in view of the scaling invariance of the desired estimate. We denote 
    \begin{align*}
        l_{k}(y,t)=(\nabla u)_{Q^s_{2^{k}}}\cdot y+(u)_{B_{2^{k}}}(t)
    \end{align*}
    for any $k\geq0$. We first note from  \begin{align*}
        \dashint_{B}A\cdot y\,dy=0
    \end{align*}
    for any $A\in\bbR^n$ and any ball $B \subset \mathbb{R}^n$ centered at the origin that
    \begin{equation}\label{ineq.solmlin}
    \begin{aligned}
        \dashint_{B_{2^{k}}}|u-l_k|\,dx&=\dashint_{B_{2^{-k}}}|u-l_k-(u-l_k)_{B_{2^{-k}}}(t)|\,dx\\
        &\leq c2^{-k}\dashint_{B_{2^{-k}}}|\nabla u(x,t)-(\nabla u)_{Q^s_{2^{-k}}}|\,dx
    \end{aligned}
    \end{equation}
    and
    \begin{align*}
        \mathrm{Tail}\left(u-l_0;B_{1}\right)=\mathrm{Tail}\left(u-l_0-(u-l_0)_{B_{1}}(t);B_{1}\right)\eqqcolon J,
    \end{align*}
    where we have also used the Poincar\'e inequality in \eqref{ineq.solmlin}.
    We next observe from \cite[Lemma 2.1]{DieKimLeeNow24} with $g$ replaced by $u-l_0$ that
    \begin{align*}
        J&\leq c\sum_{k=0}^{i}\dashint_{B_{2^{k}}}{|u-l_0-(u-l_0)_{B_{2^{k}}}(t)|}\,dy\\
        &\quad+c2^{-2si}\mathrm{Tail}\left(u-l_0-(u-l_0)_{B_{2^{2si}}}(t);B_{2^i}\right)\eqqcolon J_1+J_2,
    \end{align*}
where $c=c(n,s_0)$. We first estimate $J_1$ as 
\begin{equation}\label{est.tailj1}
\begin{aligned}
    J_1&\leq c\sum_{k=0}^{i}2^{-2sk}\dashint_{B_{2^{k}}}|u-l_{k}-(u-l_k)_{B_{2^{k}}}(t)|\,dy\\
    &\quad+c\sum_{k=0}^{i}2^{-2sk}\sum_{j=1}^{k}\dashint_{B_{2^k}}|(\nabla u)_{Q^s_{2^{j}}}\cdot y-(\nabla u)_{Q^s_{2^{j-1}}}\cdot y|\,dy\eqqcolon J_{1,1}+J_{1,2}.
\end{aligned}
\end{equation}
We now estimate $J_{1,2}$ as
\begin{align*}
    J_{1,2}&\leq c\sum_{k=0}^{i}2^{(-2s+1)k}\sum_{j=0}^{k}\dashint_{Q^s_{2^j}}|\nabla u-(\nabla u)_{Q^s_{2^j}}|\,dz\\
    &\leq c\sum_{j=0}^{i}\dashint_{Q^s_{2^j}}|\nabla u-(\nabla u)_{Q^s_{2^j}}|\,dz\sum_{k=j}^{i}2^{(-2s+1)k}\\
    &\leq c\sum_{j=0}^{i}2^{(-2s+1)j}\dashint_{Q^s_{2^j}}|\nabla u-(\nabla u)_{Q^s_{2^j}}|\,dz,
\end{align*}
where we have used Fubini's theorem along with a few simple calculations.
Plugging the estimate $J_{1,2}$ into \eqref{est.tailj1}, we get 
\begin{align*}
    J_1&\leq c\sum_{j=0}^{i}2^{-2sj}\dashint_{B_{2^{j}}}|u-(\nabla u)_{Q^s_{2^j}}\cdot y -(u)_{B_{2^{j}}}(t)|\,dy\\
    &\quad+c\sum_{j=0}^{i}2^{j(1-2s)}\dashint_{Q^s_{2^{j}}}|\nabla u-(\nabla u)_{Q^s_{2^{j}}}|\,dz
\end{align*}
for some constant $c=c(n,s_0)$. 
We next estimate $J_{2}$ as 
\begin{align*}
    J_{2}&\leq 2^{-2si}\mathrm{Tail}(u-(\nabla u)_{Q^s_1}\cdot y-(u)_{B_{2^i}}(t);B_{2^i})\\
    &\leq 2^{-2si}\mathrm{Tail}(u-(\nabla u)_{Q^s_{2^i}}\cdot y-(u)_{B_{2^i}}(t);B_{2^i})\\
    &\quad+\sum_{j=0}^i2^{(1-2s)j}\dashint_{Q^s_{2^j}}|\nabla u-(\nabla u)_{Q^s_{2^j}}|\,dz,
\end{align*}
where we have used the fact that $2^{(1-2s)i}\leq 2^{(1-2s)j}$ for any $j\leq i$ for the last inequality.
Combining all the estimates $J_1$ and $J_2$, we get
\begin{equation}\label{tail.est}
\begin{aligned}
    &\mathrm{Tail}\left(u-(\nabla u)_{Q^s_1}\cdot y-(u)_{B_1}(t);B_{1}\right)\\
    &\leq c\sum_{j=0}^{i}2^{-2sj}\dashint_{B_{2^{j}}}|u-(\nabla u)_{Q^s_{2^j}}\cdot y -(u)_{B_{2^{j}}}(t)|\,dy\\
    &\quad+c\sum_{j=0}^{i}2^{j(1-2s)}\dashint_{Q^s_{2^{j}}}|\nabla u-(\nabla u)_{Q^s_{2^{j}}}|\,dz\\
    &\quad+c2^{-2si}\mathrm{Tail}(u-(\nabla u)_{Q^s_{2^i}}\cdot y-(u)_{B_{2^i}}(t);B_{2^i}).
\end{aligned}
\end{equation}
Applying \eqref{ineq.solmlin} into the first term in the right-hand side of \eqref{tail.est} and integrating both sides of \eqref{tail.est} with respect to the time variable, we obtain the desired estimate.
\end{proof}

\section{Localization and gradient H\"older regularity} \label{sec4}

Throughout this section, we fix a parameter $s_0\in(0,1)$ and some
\begin{equation}\label{s0.cond1}
	s\in[s_0,1)
\end{equation}
to describe estimates that are stable as $s\to1$.

\subsection{Localization}
In this subsection, we prove the parabolic version of the localization argument given in \cite[Lemma 3.2]{DieKimLeeNow24}. Before providing this localization argument, we first observe the following straightforward scaling invariance of our equation \eqref{eq:nonlocaleq}.
\begin{lemma}\label{lem:scaling}
    Let 
\begin{equation*}
    u\in L^{2}(I^s_{R}(t_0); W^{s,2}(B_{R}(x_0))) \cap C\left(I^s_{R}(t_0);L^{2}(B_{R}(x_0))\right) \cap L^{1}\left(I^s_{R}(t_0);L^{1}_{2s}(\mathbb{R}^n)\right)
\end{equation*}
be a weak solution to 
\begin{align*}
    \partial_t u+\mathcal{L}u=\mu\quad\text{in }Q^s_R(z_0),
\end{align*}
with $\mu$ as in Definition \ref{def:weak}.
Then $u_R(x,t)\coloneqq u(Rx+x_0,R^{2s}t+t_0)/R^s$ is a weak solution to 
\begin{align*}
    \partial_t u_R+\mathcal{L}u_R=\mu_R\quad\text{in }Q^s_1,
\end{align*}
where $\mu_R(x,t)=R^s\mu(Rx+x_0,R^{2s}t+t_0)$.
\end{lemma}
We now prove the following lemma.
\begin{lemma}[Localization lemma]\label{lem.loc}
For some $q\geq 1$, let
\begin{equation*}
    u\in L^{2}(I^s_{5R}(t_0); W^{s,2}(B_{5R}(x_0))) \cap C\left(I^s_{5R}(t_0);L^{2}(B_{5R}(x_0))\right) \cap L^{q}\left(I^s_{5R}(t_0);L^{1}_{2s}(\mathbb{R}^n)\right)
\end{equation*}
be a weak solution to 
\begin{equation}\label{eq.loca}
    \partial_t u+\mathcal{L}u=\mu\quad\text{in }Q^s_{5R}(z_0).
\end{equation}
Let us fix $\xi\in C_c^{\infty}(B_{4R}(x_0))$ with $\xi\equiv 1$ on $B_{3R}(x_0)$ and $|\nabla \xi|\leq c/R$ for some constant $c$. Then we have that
\begin{equation*}
    w:=u\xi\in C\left(I^s_{4R}(t_0);L^{2}(\bbR^n)\right)\cap L^{2}(I^s_{4R}(t_0); W^{s,2}(\bbR^n)) 
\end{equation*}is a weak solution to 
\begin{equation}\label{eq.localized}
    \partial_t w+\mathcal{L}w=\mu+f\quad\text{in }Q^s_{2R}(z_0),
\end{equation}
where $f\in L^{q}\left(I^s_{5R/2}(t_0);L^{\infty}(B_{5R/2}(x_0))\right)$ satisfies that for any $r\in[0,5R/2]$
\begin{equation}\label{estb.loc.tail}
\begin{aligned}
    \left(\dashint_{I^s_{r}(t_0)}\|f(\cdot,t)\|^q_{ L^{\infty}(B_{5R/2}(x_0))}\,dt\right)^{\frac{1}{q}}
    \leq cR^{-2s}\left(\dashint_{I^s_{r}(t_0)}\mathrm{Tail}(u;B_{3R}(x_0))^q\,dt\right)^{\frac{1}{q}},
\end{aligned}
\end{equation}
where $c=c(n,\Lambda)$. In addition, if $u\in L^q(I^s_{5R}(t_0);C^{0,\beta}(B_{5R}(x_0)))$ for some $\beta\in(0,1]$, then 
\begin{equation*}
    w:=u\xi\in C\left(I^s_{4R}(t_0);L^{2}(\bbR^n)\right)\cap L^{2}(I^s_{4R}(t_0); W^{s,2}(\bbR^n)) \cap L^\infty\left(I^s_{4R}(t_0);C^{0,\beta}(\bbR^n)\right)
\end{equation*}is a weak solution to \eqref{eq.localized}, where $f\in L^{q}\left(I^s_{5R/2}(t_0);C^{0,\beta}(B_{5R/2}(x_0))\right)$ with the estimate
\begin{equation}\label{esth.loc.tail}
\begin{aligned}
    &\left(\dashint_{I^s_{\frac{5R}2}(t_0)}[f(\cdot,t)]^q_{ C^{0,\beta}(B_{\frac{5R}2}(x_0))}\,dt\right)^{\frac{1}{q}}\\
    &\leq cR^{-2s}\left(\dashint_{I^s_{\frac{5R}2}(t_0)}[u(\cdot,t)]^q_{C^{0,\beta}(B_{3R}(x_0))}\,dt\right)^{\frac{1}{q}}\\
    &\quad+cR^{-(2s+\beta)} \left(\dashint_{I^s_{\frac{5R}2}(t_0)}\|u(\cdot,t)\|^q_{ L^{\infty}(B_{3R}(x_0))}+\mathrm{Tail}(u;B_{3R}(x_0))^q\,dt\right)^{\frac{1}{q}}
\end{aligned}
\end{equation}
for some constant $c=c(n,s_0,\Lambda)$.
\end{lemma}
\begin{proof}
We may assume $z_0=(x_0,t_0)=0$ by Lemma \ref{lem:scaling}. We first prove that $w=u\xi$ is a weak solution to 
\eqref{eq.localized} with 
\begin{align*}
    f(x,t)&=2(1-s)\int_{\RRn\setminus B_{3R}}{{\Phi}\left(\frac{u(x,t)-(u\xi)(y,t)}{|x-y|^s}\right)}\frac{\,dy}{|x-y|^{n+s}}\\
    &\quad-2(1-s)\int_{\RRn\setminus B_{3R}}{{\Phi}\left(\frac{u(x,t)-u(y,t)}{|x-y|^s}\right)}\frac{\,dy}{|x-y|^{n+s}}
\end{align*}
for any $z=(x,t)\in Q^s_{5R/2}$.
Let us fix $\psi\in L^{2}(I^s_{2R};W^{s,2}(B_{2R}))\cap W^{1,2}(I^s_{2R};L^2(B_{2R}))$ where the support of $\psi$ is compactly contained in $B_{2R}$. By testing the function $\psi$ to \eqref{eq.loca}, after a few simple manipulations, we obtain
\begin{equation*}
\begin{aligned}
    &-\int_{I}\int_{B_{2R}}u\partial_{t}\psi\,dz\\
    &\quad+(1-s)\int_{I}\int_{\RRn}\int_{\RRn}{\Phi}\left(\frac{(u\xi)(x,t)-(u\xi)(y,t)}{|x-y|^{s}}\right)\frac{\psi(x,t)-\psi(y,t)}{|x-y|^{n+s}}\,dy\,dz\\
    &=\int_{I}\int_{B_{2R}}\mu\psi\,dz-\int_{B_{2R}}u\psi\,dx\Bigg\rvert_{t=t_1}^{t=t_2}\\
    &\quad-(1-s)\int_{I}\int_{\RRn}\int_{\RRn}{\Phi}\left(\frac{u(x,t)-u(y,t)}{|x-y|^{s}}\right)\frac{\psi(x,t)-\psi(y,t)}{|x-y|^{n+s}}\,dy\,dz\\
    &\quad+(1-s)\int_{I}\int_{\RRn}\int_{\RRn}{\Phi}\left(\frac{(u\xi)(x,t)-(u\xi)(y,t)}{|x-y|^{s}}\right)\frac{\psi(x,t)-\psi(y,t)}{|x-y|^{n+s}}\,dy\,dz
    \eqqcolon\sum_{i=1}^{4}J_i
\end{aligned}
\end{equation*}
for any $I=[t_1,t_2]\Subset I^s_{2R}$.
Using 
the facts $\xi\equiv 1$ on $B_{3R}$ and $\psi(x,\cdot)\equiv 0$ in $\RRn\setminus B_{2R}$, we observe
\begin{equation*}
\begin{aligned}
    J_3+J_4&=2(1-s)\int_{I}\int_{B_{2R}}\int_{\RRn\setminus B_{3R}}{\Phi}\left(\frac{u(x,t)-(u\xi)(y,t)}{|x-y|^{s}}\right)\frac{\psi(x,t)}{|x-y|^{n+s}}\,dy\,dz\\
    &\quad-2(1-s)\int_{I}\int_{B_{2R}}\int_{\RRn\setminus B_{3R}}{\Phi}\left(\frac{u(x,t)-u(y,t)}{|x-y|^{s}}\right)\frac{\psi(x,t)}{|x-y|^{n+s}}\,dy\,dz
    =\int_{I}\int_{B_{2R}}f\psi\,dz.
\end{aligned}
\end{equation*}
Using this along with the facts that 
\begin{equation*}
    \int_{I^s_{2R}}\int_{B_{2R}}u\partial_{t}\psi\,dz=\int_{I^s_{2R}}\int_{B_{2R}}(u\xi)\partial_{t}\psi\,dz
\end{equation*}
and
\begin{align*}
    J_1+J_2=\int_{I}\int_{B_{2R}}\mu\psi\,dz-\int_{B_{2R}}(u\xi)\psi\,dx\Bigg\rvert_{t=t_1}^{t=t_2},
\end{align*}
we verify that $w=u\xi$ is a weak solution to \eqref{eq.localized}.
Now, using \eqref{pt : assmp.phi} and the relation
\begin{align*}
    |x-y|\geq \frac{|y|}{6}\quad\text{for any }x\in B_{5R/2}\text{ and }y\in B_{3R}^c,
\end{align*}  
we obtain that for any $x\in B_{5R/2}$,
\begin{equation*}
\begin{aligned}
    |f(x,t)|\leq c(1-s)\int_{\bbR^n\setminus B_{3R}}\frac{|u(y,t)|}{|y|^{n+2s}}\,dy
\end{aligned}
\end{equation*}
holds for some constant $c=c(n,\Lambda)$, which implies \eqref{estb.loc.tail}. We now suppose 
\begin{equation*}
    u\in L^\infty(I^s_{5R};C^{0,\beta}(B_{5R}))
\end{equation*}
for some $\beta\in(0,1]$. As in the estimate of \cite[Equation (3.4) in Lemma 3.2]{DieKimLeeNow24}, we get 
\begin{align*}
    [f(\cdot,t)]_{C^{0,\beta}(B_{5R/2})}&\leq cR^{-2s}[u(\cdot,t)]_{C^{0,\beta}(B_{3R})}\\
    &\quad +cR^{-(2s+\beta)}\left[\|u(\cdot,t)\|_{ L^{\infty}(B_{3R})}+\mathrm{Tail}(u;B_{3R})\right]
\end{align*}
for some constant $c=c(n,s_0,\Lambda)$. From this, we deduce \eqref{esth.loc.tail}.
\end{proof}

\subsection{H\"older regularity}

In this subsection, we discuss some H\"older regularity results for nonlocal parabolic equations that concern the regularity of the solution itself rather than its gradient. 

Before that, for convenience of notation,
for any $g\in L^q(Q^s_R(z_0))\cap L^q(I^s_R(t_0);L^1_{2s}(\bbR^n))$ with $q\geq1$, we write
\begin{align}\label{defn.exctilde}
	\widetilde{E}^q_{\mathrm{loc}}(g;Q^s_R(z_0))\coloneqq\left(\dashint_{Q^s_R(z_0)}|g|^q\,dz\right)^{\frac1q},
\end{align}
\begin{align*}
	\widetilde{E}^q(g;Q^s_R(z_0))\coloneqq\widetilde{E}^q_{\mathrm{loc}}(g;Q^s_R(z_0))+\left(\dashint_{I^s_R(t_0)}\mathrm{Tail}(g;B_R(x_0))^q\,dt\right)^{\frac1q},
\end{align*}
By recalling the notation \eqref{not.exc}, we observe that
\begin{equation}\label{rel.exc}
	\begin{aligned}
		E^q(g;Q^s_R(z_0))&\leq cE^q_{\mathrm{loc}}(g;Q^s_R(z_0))\\
		&\quad+\left(\dashint_{I^s_R(t_0)}\mathrm{Tail}(g-(g)_{B_{R}(x_0)}(t);B_R(x_0))^q\,dt\right)^{\frac1q},
	\end{aligned}
\end{equation}
where $c=c(n,s_0)$. When $q=1$, for convenience we simply write $\widetilde{E}_{\mathrm{loc}}(g;Q^s_R(z_0)) \coloneqq \widetilde{E}^1_{\mathrm{loc}}(g;Q^s_R(z_0))$ and $\widetilde{E}(g;Q^s_R(z_0)) \coloneqq \widetilde{E}^1(g;Q^s_R(z_0))$.

The first result we state in this section is a known H\"older regularity result for solutions to nonlocal linear parabolic equations, see \cite[Theorem 1.5]{KaWe23}.
\begin{lemma}\label{lem.lin.hol}
	Let $g\in L^1(I^s_R(t_0);L^\infty(B_R(x_0)))$ and
	\begin{equation*}
		v\in L^2(I^s_R(t_0);W^{s,2}(B_R(x_0)))\cap C(I^s_R(t_0);L^2(B_R(x_0)))\cap L^1(I^s_R(t_0);L^1_{2s}(\bbR^n))
	\end{equation*}
	be a weak solution to
	\begin{align}\label{eq.linear}
		\partial_t v+\widetilde{\mathcal{L}}_{K}v=g\quad\text{in }Q^s_R(z_0),
	\end{align}
	where 
	\begin{align}\label{defn.nonloc}
		\widetilde{\mathcal{L}}_{K}v(x,t)=\mathrm{P.V.}\int_{\bbR^n}\left({v(x,t)-v(y,t)}\right)K(x,y,t)\,dy
	\end{align}
	for some measurable function $K:\bbR^n\times\bbR^n\times I^s_R(t_0)$ such that 
	\begin{align}\label{lin.kern.cond}
		K(x,y,t)=K(y,x,t)
	\end{align}
	and
	\begin{align}\label{lin.kern.cond2}
		\frac{(1-s)\Lambda^{-1}}{|x-y|^{n+2s}}\leq K(x,y,t)\leq \frac{(1-s)\Lambda}{|x-y|^{n+2s}}
	\end{align}
	for any $x,y\in\bbR^n$, $x\neq y$ and $t\in I^s_R(t_0)$. Then we have 
	\begin{align*}
		\|v\|_{L^\infty(Q^s_{R/2}(z_0))}\leq c\widetilde{E}(v;Q^s_R(z_0))+cR^{2s}\|g\|_{L^1(I^s_R(t_0);L^\infty(B_R(x_0)))}
	\end{align*}
	for some constant $c=c(n,s_0,\Lambda)$. In addition, if $v\in L^q(I^s_R(t_0);L^\infty(B_R(x_0)))$ and $g=0$ with $q>1$, then there is a constant $\gamma=\gamma(n,s_0,\Lambda,q)$ such that
	\begin{align*}
		R^\gamma[v]_{C^{0,\gamma}(Q^s_{R/2}(z_0))}\leq cE^q(v;Q^s_R(z_0))
	\end{align*}
	for some constant $c=c(n,s_0,\Lambda,q)$.
\end{lemma}

With slight modifications, we obtain a corresponding result for nonlinear nonlocal parabolic equations.
\begin{lemma}\label{lem.lo.hol}
	Let 
	\begin{equation*}
		v\in L^2(I^s_R(t_0);W^{s,2}(B_R(x_0)))\cap C(I^s_R(t_0);L^2(B_R(x_0)))\cap L^1(I^s_R(t_0);L^1_{2s}(\bbR^n))
	\end{equation*}
	be a weak solution to 
	\begin{align}\label{eq.nonhomo}
		\partial_t v+\mathcal{L}v=g\quad\text{in }Q^s_R(z_0),
	\end{align}
	where $g\in L^1(I^s_R(t_0);L^\infty(B_R(x_0)))$.
	Then we have 
	\begin{align}\label{bdd.est.zero}
		\|v\|_{L^\infty(Q^s_{R/2}(z_0))}\leq c\widetilde{E}(v;Q^s_R(z_0))+cR^{2s}\|g\|_{L^1(I^s_R(t_0);L^\infty(B_R(x_0)))}
	\end{align}
	for some constant $c=c(n,s_0,\Lambda)$. In addition, if $v\in L^q(I^s_R(t_0);L^1_{2s}(\bbR^n))$ and $g\in L^q(I^s_R(t_0);L^\infty(B_R(x_0)))$ for some constant $q>1$, then there is a constant $\gamma=\gamma(n,s_0,\Lambda,q)$ such that
	\begin{align}\label{hol.est.zero}
		R^\gamma[v]_{C^{0,\gamma}(Q^s_{R/2}(z_0))}\leq cE^q(v;Q^s_R(z_0))+cR^{2s}\left(\dashint_{I^s_R(t_0)}\|g(\cdot,t)\|_{L^\infty(B_R(x_0))}^q\,dt\right)^{\frac1q}
	\end{align}
	for some constant $c=c(n,s_0,\Lambda,q)$.
\end{lemma}
\begin{proof}
	We first remark that if $v$ is a weak solution to \eqref{eq.nonhomo}, then $v$ is also a weak solution to \eqref{eq.linear} with a measurable kernel $K$ satisfying \eqref{lin.kern.cond} and \eqref{lin.kern.cond2}. Indeed, by following the same lines as in \cite[Remark 4.1]{DieKimLeeNow24}, we deduce that if 
	\begin{align*}
		K(x,y,t)\coloneqq (1-s)\Phi\left(\frac{v(x,t)-v(y,t)}{|x-y|^s}\right)(v(x,t)-v(y,t))^{-1}|x-y|^{-(n+s)},
	\end{align*}
	then $v$ is a weak solution to \eqref{eq.linear} with \eqref{lin.kern.cond} and \eqref{lin.kern.cond2}.
	By Lemma \ref{lem.lin.hol}, we obtain \eqref{bdd.est.zero}. 
	
	We now prove \eqref{hol.est.zero}. Let us fix $Q^s_{10r}(z_1)\Subset Q^s_R(z_0)$.
	By the localization argument given in Lemma \ref{lem.loc} below, we observe that
	\begin{align*}
		w\coloneqq v\xi\in C(I^s_{4r}(t_1);L^2(\bbR^n))\cap L^2(I^s_{4r}(t_0);W^{s,2}(\bbR^n))
	\end{align*}
	is a weak solution to 
	\begin{align}\label{eq.trans}
		\partial_t w+\mathcal{L}w=g+f\quad\text{in }Q^s_{2r}(z_1)
	\end{align}
	for some $f\in L^q(I^s_{5r/2}(t_1);L^\infty(B_{5r/2}(x_1)))$, where $\xi$ is the function determined in Lemma \ref{lem.loc} with $R=r$ and $z_0=z_1$. Therefore, using  perturbation arguments as in \cite[Lemma 3.3--3.5]{ByuKimKim23} along with Lemma \ref{lem.lin.hol} and carefully tracking the factor $1-s$ in the constants, we get
	\begin{align*}
		r^\gamma[w]_{C^{0,\gamma}(Q^s_{r}(z_1))}&\leq c\|w\|_{L^\infty(Q^s_{2r}(z_1))}+c\sup_{t\in I^s_{2r}(t_1)}\mathrm{Tail}(w(\cdot,t);B_{2r}(x_1))\\
		&\quad+cr^{2s}\left(\dashint_{I^s_r(t_1)}\|(g+f)(\cdot,t)\|_{L^\infty(B_r(x_1))}^q\,dt\right)^{\frac1q}
	\end{align*}
	for some constant $c=c(n,s_0,\Lambda,q)$, where $\gamma=\gamma(n,s_0,\Lambda,q)\in(0,1)$. Using the fact that $w=v\xi$, \eqref{estb.loc.tail} and \eqref{bdd.est.zero}, we have
	\begin{align*}
		r^\gamma[v]_{C^{0,\gamma}(Q^s_{r}(z_1))}&\leq c\|v\|_{L^\infty(Q^s_{5r}(z_1))}+c\left(\dashint_{I^s_{5r}(t_1)}\mathrm{Tail}(v(\cdot,t);B_{5r}(x_1))^q\,dt\right)^{\frac1q}\\
		&\quad+cr^{2s}\left(\dashint_{I^s_{5r}(t_1)}\|g(\cdot,t)\|_{L^\infty(B_{5r}(x_1))}^q\,dt\right)^{\frac1q}\\
		&\leq c\widetilde{E}^q(v;Q^s_{10r}(z_1))+cr^{2s}\left(\dashint_{I^s_{10r}(t_1)}\|g(\cdot,t)\|_{L^\infty(B_{10r}(x_1))}^q\,dt\right)^{\frac1q},
	\end{align*}
	where $c=c(n,s_0,\Lambda,q)$. By standard covering arguments along with the fact that $v-(v)_{Q^s_R(z_0)}$ is a also weak solution to \eqref{eq.nonhomo}, we finally obtain \eqref{hol.est.zero}.
\end{proof}

\begin{remark}\label{rmk.hollin} \normalfont
	We point out that if $v$ is a weak solution to \eqref{eq.linear}, then we observe from the proof of Lemma \ref{lem.loc} that $w=v\xi$ is a weak solution to 
	\begin{align*}
		\partial_t w+\widetilde{\mathcal{L}}_Kw=\mu+f\quad\text{in }Q^s_{2R/5}(z_0)
	\end{align*}
	with $f$ satisfying \eqref{estb.loc.tail} and with $R$ replaced by $R/5$. Therefore, the weak solution $v$ also satisfies \eqref{hol.est.zero} by following the same lines as in the proof of Lemma \ref{lem.lo.hol}.
\end{remark}

\subsection{A parabolic Poincar\'e-type inequality}

We next give a version of the classical parabolic Poincar\'e inequality that holds for solutions to linear nonlocal parabolic equations.
\begin{lemma}\label{lem.glu}Let
	\begin{equation*}
		v\in L^2(I^s_R(t_0);W^{s,2}(B_R(x_0)))\cap C(I^s_R(t_0);L^2(B_R(x_0)))\cap L^1(I^s_R(t_0);L^1_{2s}(\bbR^n))
	\end{equation*}
	be a weak solution to 
	\begin{align}\label{eq.glu}
		\partial_t v+\widetilde{\mathcal{L}}_Kv=0\quad\text{in }Q^s_R(z_0),
	\end{align}
	where the nonlocal operator $\widetilde{\mathcal{L}}_K$ is defined in \eqref{defn.nonloc}. If $\nabla v\in L^1(Q^s_R(z_0))$, then we have 
	\begin{align*}
		E^q_{\mathrm{loc}}(v;Q^s_{R/2}(z_0))&\leq cR\dashint_{Q^s_{3R/4}(z_0)}|\nabla v|\,dz\\
		&\quad+c\dashint_{I^s_{3R/4}(t_0)}\mathrm{Tail}(v-(v)_{B_{3R/4}(x_0)}(t);B_{3R/4}(x_0))\,dt
	\end{align*}
	for some constant $c=c(n,s_0,\Lambda)$ and for any $q>1$.
\end{lemma}
\begin{proof}
	We may assume $z_0=0$ and $R=1$ in view of the scaling invariance of the desired estimate.
	Using \eqref{bdd.est.zero} together with the fact that $v-(v)_{Q^s_{1/2}}$ is also a weak solution to \eqref{eq.glu} with $R=1$ and $z_0=0$, we have 
	\begin{align*}
		E^q_{\mathrm{loc}}(v;Q^s_{1/2})&\leq cE(v;Q^s_{3/4})\\
		&\leq cE_{\mathrm{loc}}(v;Q^s_{3/4})+c\dashint_{I^s_{3/4}}\mathrm{Tail}(v-(v)_{B_{3/4}}(t);B_{3/4})\,dt
	\end{align*}
	for some constant $c=c(n,s_0,\Lambda)$. We now observe
	\begin{align*}
		\dashint_{Q^s_{3/4}}|v-(v)_{Q^s_{3/4}}|^q\,dz
		&\leq c\dashint_{Q^s_{3/4}}|v-(v)_{B_{3/4}}(t)|^q\,dz\\
		&\quad +c\dashint_{I^s_{3/4}}|(v)_{Q^s_{3/4}}-(v)_{B_{3/4}}(t)|^q\,dt\eqqcolon J_1+J_2
	\end{align*}
	for some constant $c=c(n,s_0,\Lambda)$. Let us choose a cutoff function $\psi \in C^\infty_c(B_{7/8})$ with $\psi\equiv 1$ on $B_{3/4}$ and $\abs{\nabla\psi}\leq c$. We next observe from the proof of \cite[Lemma A.1]{ByuKimKum23p} that 
	\begin{align*}
		J_2\leq c\dashint_{I^s_{3/4}}|(v)_{B_{3/4}(t)}-(v)^{\psi}_{B_{3/4}}(t)|^q\,dt+c\sup_{t,\tau\in I^s_{3/4}}|(v)^{\psi}_{B_{3/4}}(t)-(v)^{\psi}_{B_{3/4}}(\tau)|^q
	\end{align*}
	for some constant $c=c(n,s_0,\Lambda,q)$, where we denote 
	\begin{align*}
		(v)^{\psi}_{B_{3/4}}(t)=\frac{1}{\|\psi\|_{L^1}}\int_{B_{3/4}}(v\psi)(x,t)\,dx.
	\end{align*} 
	By following the same lines as in the proof of the gluing lemma given by \cite[Lemma 4.5]{ByuKimKim23} and taking into account the factor $1-s$ in front of the nonlocal operator, we get 
	\begin{equation}\label{est.j.glu}
		\begin{aligned}
			J_2&\leq cJ_1+c\left((1-s)\int_{B_{3/4}}\dashint_{Q^s_{3/4}}\frac{|v(x,t)-v(y,t)|}{|x-y|^{n+2s-1}}\,dz\,dy\right)^q\\
			&\quad+c\left((1-s)\int_{\bbR^n\setminus B_{3/4}}\dashint_{Q^s_{7/8}}\frac{|v(x,t)-v(y,t)|}{|x-y|^{n+2s}}\,dz\,dy\right)^q\\
			&\leq cJ_1+c\left((1-s)\int_{B_{3/4}}\dashint_{Q^s_{3/4}}\frac{|v(x,t)-v(y,t)|}{|x-y|^{n+2s-1}}\,dz\,dy\right)^q\\
			&\quad +c\left(\dashint_{I^s_{3/4}}\mathrm{Tail}(v-(v)_{B_{3/4}}(t);B_{3/4})\,dt\right)^q,
		\end{aligned}
	\end{equation}
	where $c=c(n,s_0,\Lambda)$. As in the proof of \cite[Proposition 2.2]{DinPalVal12} with $u$ replaced by $v-(v)_{B_{3/4}}(t)$, we deduce
	\begin{align*}
		\dashint_{Q^s_{3/4}}\int_{B_{3/4}}\frac{|v(x,t)-v(y,t)|}{|x-y|^{n+2s-1}}\,dy\,dz&\leq \frac{c}{1-s}\dashint_{Q^s_{3/4}}|\nabla v|\,dz\\
		&\quad+ c\dashint_{Q^s_{3/4}}|v-(v)_{B_{3/4}}(t)|\,dz
	\end{align*}
	for some constant $c=c(n,s_0)$. We now combine the above two inequalities and the estimate $J_1$ with H\"older's inequality and Poincar\'{e}'s inequality to conclude that
	\begin{align*}
		J_1+J_2 &\leq c\dashint_{Q^s_{3/4}}|\nabla v|\,dz+c\dashint_{I^s_{3/4}}\mathrm{Tail}(v-(v)_{B_{3/4}}(t);B_{3/4})\,dt
	\end{align*}
	for some constant $c=c(n,s_0,\Lambda)$, which completes the proof.
\end{proof}

\subsection{Gradient H\"older regularity}
In this subsection, we establish the H\"older continuity of the gradient of weak solutions to \eqref{eq:nonlocaleq} with $s\in(0,1)$ and $\mu=0$.

Let us fix $\beta\in(0,1]$. We first prove that a given regular and localized solution $w$, $\frac{\delta_h w}{|h|^\beta}$ is a weak solution of a nonhomogeneous weighted fractional heat equation.
\begin{lemma}\label{lem.linsol}
Let 
\begin{equation*}
    w\in C\left(I^s_{2R}(t_0);L^{2}(\bbR^n)\right)\cap L^{2}(I^s_{2R}(t_0); W^{s,2}(B_{2R}(x_0)))\cap L^q(I^s_{2R}(t_0); L^1_{2s}(\bbR^n))
\end{equation*}be a weak solution to 
\begin{equation}\label{eq.gra.frde}
    \partial_{t}w+\mathcal{L}w=f\quad\text{in }{Q}_{2R}(z_0),
\end{equation}  
where $f\in L^{q}\left(I^s_{2R}(t_0);L^\infty\left(B_{2R}(x_0)\right)\right)$ for some $q>1$.
Let us fix $h\in B_{R/100}\setminus \{0\}$ and $\beta\in(0,1]$. Then $\widetilde{w}\coloneqq\frac{\delta_h w}{|h|^\beta}$ is a weak solution to 
\begin{align}\label{eq.gra.dfq}
    \partial_t \widetilde{w}+\widetilde{\mathcal{L}}_{K}\widetilde{w}=\frac{\delta_h f}{|h|^\beta}\quad\text{in }Q^s_R(z_0)
\end{align}
for some symmetric measurable kernel $K$ satisfying
\eqref{lin.kern.cond} and \eqref{lin.kern.cond2}.
\end{lemma}
\begin{proof}
    We now fix $|h|<R/100$. We next choose $$\phi\in L^2(t_1,t_2;W^{s,2}(B_R(x_0)))\cap W^{1,2}(t_1,t_2;L^2(B_R(x_0)))$$ with the support compactly contained in the spatial direction, where $[t_1,t_2]\Subset I^s_R(t_0)$. By testing $\delta_{-h}\phi$ to \eqref{eq.gra.frde}, we get
    \begin{equation}\label{eq.test.lin}
    \begin{aligned}
        &-\int_{t_1}^{t_2}\int_{B_{2R}(x_0)}\delta_h w\partial_t \phi\,dz\\
        &\quad+(1-s)\int_{t_1}^{t_2}\int_{\bbR^n}\int_{\bbR^n}\left[\Phi\left(\frac{w_h(x,t)-w_h(y,t)}{|x-y|^s}\right)-\Phi\left(\frac{w(x,t)-w(y,t)}{|x-y|^s}\right)\right]\\
        &\qquad\qquad\qquad\qquad\qquad\quad\times \frac{\phi(x,t)-\phi(y,t)}{|x-y|^{n+s}}\,dx\,dy\,dt\\
        &=\int_{t_1}^{t_2}\int_{B_{2R}(x_0)}\delta_h f \phi\,dz-\int_{B_{2R}(x_0)}\delta_hw\, \phi\,dx\Bigg\lvert_{t=t_1}^{t=t_2}.
    \end{aligned}
    \end{equation}
We now write 
\begin{align*}
    K(x,y,t)=(1-s)\frac{\Phi\left(\frac{w_h(x,t)-w_h(y,t)}{|x-y|^s}\right)-\Phi\left(\frac{w(x,t)-w(y,t)}{|x-y|^s}\right)}{\frac{\delta_h w(x,t)-\delta_h w(y,t)}{|x-y|^s}} |x-y|^{-n-2s}
\end{align*}
to see that the coefficient function $K$ satisfies \eqref{lin.kern.cond} and \eqref{lin.kern.cond2}. Dividing both sides of \eqref{eq.test.lin} by $|h|^\beta$ along with the fact that $\phi(\cdot,x)=0$ on $B_{2R}(x_0)\setminus B_R(x_0)$ yields 
\begin{align*}
        &-\int_{t_1}^{t_2}\int_{B_R(x_0)}\frac{\delta_h w}{|h|^\beta}\partial_t \phi\,dz\\
        &\quad+(1-s)\int_{t_1}^{t_2}\int_{\bbR^n}\int_{\bbR^n}\frac{\delta_hw(x,t)-\delta_hw(y,t)}{|h|^\beta|x-y|^s}\frac{\phi(x,t)-\phi(y,t)}{|x-y|^{n+s}}K(x,y,t)\,dx\,dy\,dt\\
        &=\int_{t_1}^{t_2}\int_{B_R(x_0)}\frac{\delta_h f}{|h|^\beta} \phi\,dz-\int_{B_{R}(x_0)}\frac{\delta_h w}{|h|^\beta}\, \phi\,dx\Bigg\lvert_{t=t_1}^{t=t_2},
    \end{align*}
which implies $\frac{\delta_h w}{|h|^\beta}$ is a weak solution to \eqref{eq.gra.dfq}.
\end{proof}
Using Lemma \ref{lem.linsol}, we obtain the following H\"older estimates for the difference quotients of solutions to \eqref{eq.gra.frde}.
\begin{lemma}\label{lem.hol}
   Let $w$ be a weak solution to \eqref{eq.gra.frde}
with $f\in L^q(I^s_{2R}(t_0);C^{0,\beta}(B_{2R}(x_0)))$ with $q\geq 1$ and $\beta\in(0,1]$. Let us fix $h\in B_{R/100}\setminus \{0\}$. Then we have 
\begin{align*}
    \left\|\frac{\delta_h w}{|h|^\beta}\right\|_{L^\infty(Q^s_{R/2}(z_0))}&\leq c\widetilde{E}\left(\frac{\delta_h w}{|h|^\beta};Q^s_R(z_0)\right)\\
    &\quad+cR^{2s}\dashint_{I^s_R(t_0)}[f(\cdot,t)]_{C^{0,\beta}(B_{R}(x_0))}\,dt
\end{align*}
for some constant $c=c(n,s_0,\Lambda)$. Morevoer, if $q>1$, then there is a constant $\gamma=\gamma(n,s_0,\Lambda,q)\in(0,1)$ such that
\begin{align*}
    R^{\gamma}\left[\frac{\delta_h w}{|h|^\beta}\right]_{C^{0,\gamma}(Q^s_{R/2}(z_0))}&\leq c{E}^q\left(\frac{\delta_h w}{|h|^\beta};Q^s_R(z_0)\right)\\
    &\quad+cR^{2s}\left(\dashint_{I^s_R(t_0)}[f(\cdot,t)]^q_{C^{0,\beta}(B_{R}(x_0))}\,dt\right)^{\frac1q},
\end{align*}
where $c=c(n,s_0,\Lambda,q)$.
\end{lemma}
\begin{proof}
    By Lemma \ref{lem.linsol}, $\frac{\delta_hw}{|h|^\beta}$ is a weak solution to \eqref{eq.gra.dfq}. Then by Lemma \ref{lem.lin.hol} and Remark \ref{rmk.hollin}, we obtain 
\begin{equation}\label{hol.ineq1}
\begin{aligned}
    \left\|\frac{\delta_h w}{|h|^\beta}\right\|_{L^\infty(Q^s_{R/2}(z_0))}&\leq c\widetilde{E}\left(\frac{\delta_h w}{|h|^\beta};Q^s_{3R/4}(z_0)\right)\\
    &\quad+cR^{2s}\dashint_{I^s_{3R/4}(t_0)}\left\|\frac{\delta_h f}{|h|^\beta}\right\|_{L^\infty(B_{3R/4}(x_0))}\,dt
\end{aligned}
\end{equation}
for some constant $c=c(n,s_0,\Lambda)$, and
\begin{equation}\label{hol.ineq2}
\begin{aligned}
    R^{\gamma}\left[\frac{\delta_h w}{|h|^\beta}\right]_{C^{0,\gamma}(Q^s_{R/2}(z_0))}&\leq c\widetilde{E}^q\left(\frac{\delta_h w}{|h|^\beta};Q^s_{3R/4}(z_0)\right)\\
    &\quad+cR^{2s}\left(\dashint_{I^s_{3R/4}(t_0)}\left\|\frac{\delta_h f}{|h|^\beta}\right\|_{L^\infty(B_{3R/4}(x_0))}^q\,dt\right)^{\frac1q}
\end{aligned}
\end{equation}
for some constant $\gamma=\gamma(n,s_0,\Lambda,q)\in(0,1)$, where $c=c(n,s_0,\Lambda,q)$. Applying the inequality 
\begin{align*}
    \left\|\frac{\delta_h f(\cdot,t)}{|h|^\beta}\right\|_{L^\infty(B_{3R/4}(x_0))}\leq [f(\cdot,t)]_{C^{0,\beta}(B_{R}(x_0))}
\end{align*}
to the second terms in the right-hand side of \eqref{hol.ineq1} and \eqref{hol.ineq2}, respectively, we obtain the desired estimates.
\end{proof}

Using Lemma \ref{lem.loc} and Lemma \ref{lem.hol}, we now prove our first main theorem.
\begin{proof}[Proof of Theorem \ref{thm.hol}.]
Let us fix the constant $\gamma=\gamma(n,s_0,\Lambda,q)=\min\{\gamma_1,\gamma_2\}\in(0,1)$, where the constants $\gamma_1$ and $\gamma_2$ are determined in Lemma \ref{lem.lin.hol} and Lemma \ref{lem.lo.hol}, respectively. Then we have $u\in C^{0,\gamma}_{\mathrm{loc}}(\Omega_T)$ with the estimate 
\begin{align*}
    R^{\gamma}[u]_{C^{0,\gamma}(Q^s_{R/2}(z_1))}\leq cE^q(u/R;Q^s_R(z_1))
\end{align*}
for some constant $c=c(n,s_0,\Lambda,q)$, whenever $Q^s_R(z_1)\Subset \Omega_T$.
Let us assume that $Q^s_{5R}(z_1)\Subset \Omega_T$. Then $u_R(z)=u(x_1+Rx,t_1+R^{2s}t)/R^{s}$ is a weak solution to 
\begin{align*}
    \partial_t u_R+\mathcal{L}u_R=0\quad\text{in }Q^s_5.
\end{align*}
By the localization argument given in Lemma \ref{lem.loc}, we obtain that
\begin{align*}
    w=u_R\xi
\end{align*}
is a weak solution to 
\begin{equation*}
    \partial_t w+\mathcal{L}w=f\quad\text{in }Q^s_2,
\end{equation*}
where $f\in L^q(I^s_{5/2}(t_0);C^{0,\gamma}(B_{5/2}(x_0)))$ and the cutoff function $\xi$ is determined in Lemma \ref{lem.loc} with $R=1$ and $z_0=0$. 
By Lemma \ref{lem.hol}, we get 
\begin{align*}
    \left[\frac{\delta_h w}{|h|^\gamma}\right]_{C^{0,\gamma}(Q^s_{1/2})}&\leq c{E}^q\left(\frac{\delta_h w}{|h|^\gamma};Q^s_1\right)
    +c\left(\dashint_{I^s_1}[f(\cdot,t)]^q_{C^{0,\gamma}(B_{1})}\,dt\right)^{\frac1q}\eqqcolon J
\end{align*}
for some constant $c=c(n,s_0,\Lambda,q)$, where $h\in B_{1/100}\setminus\{0\}$. We first observe
\begin{align*}
    {E}^q\left(\frac{\delta_h w}{|h|^\gamma};Q^s_1\right)&\leq \left\|\frac{\delta_h u_R}{|h|^\gamma}\right\|_{L^q(I^s_1;L^\infty(B_4))}+\widetilde{E}^q(u_R;Q^s_4)\\
    &\leq c[u_R]_{L^q(I^s_1;C^{0,\gamma}(B_4))}+\widetilde{E}^q(u_R;Q^s_4).
\end{align*}
Using this and \eqref{esth.loc.tail} with $z_0=0$ and $R=1$, we further estimate $J$ as 
\begin{align*}
    \left[\frac{\delta_h w}{|h|^\gamma}\right]_{L^\infty(I^s_{1/2};C^{0,\gamma}(B_{1/2}))}+\left[\frac{\delta_h w}{|h|^\gamma}\right]_{C^{0,\gamma}(Q^s_{1/2})}\leq c\|u_R\|_{L^q(I^s_1;C^{0,\gamma}(B_4))}+c\widetilde{E}^q(u_R;Q^s_4)
\end{align*}
for some constant $c=c(n,s_0,\Lambda,q)$.
We now fix a positive integer $i_0=i_0(n,s_0,\Lambda,q)$ such that
\begin{align*}
    i_0 \gamma\leq 1<(i_0+1)\gamma.
\end{align*}
If $i_0=1$, then by \cite[Lemma 3.7]{DieKimLeeNow24}, we get 
\begin{align*}
    [w]_{L^\infty(I^s_{1/2};C^{0,1}(B_{1/2}))}\leq c\|u_R\|_{L^q(I^s_1;C^{0,\gamma}(B_4))}+c\widetilde{E}^q(u_R;Q^s_4)
\end{align*}
for some constant $c=c(n,s_0,\Lambda,q)$.
Using a scaling argument,
\begin{align*}
    R[\nabla u]_{L^\infty(Q^s_{R/2}(z_1))}\leq c\widetilde{E}^{q}(u;Q^s_{8R}(z_1))
\end{align*}
holds. A standard covering argument yields 
\begin{align}\label{ineq.gra.est}
    R[\nabla u]_{L^\infty(Q^s_{R/2}(z_0))}\leq c\widetilde{E}^{q}(u;Q^s_{R}(z_0)),
\end{align}
where $c=c(n,s_0,\Lambda,q)$ whenever $Q^s_R(z_0)\Subset \Omega_T$.

Suppose if $i_0>1$, then by \cite[Lemma 3.7]{DieKimLeeNow24}, we obtain 
\begin{align*}
    [w]_{L^\infty(I^s_{1/2};C^{0,2\gamma}(B_{1/2}))}\leq c\|u_R\|_{L^q(I^s_1;C^{0,\gamma}(B_4))}+c\widetilde{E}^q(u_R;Q^s_4)
\end{align*}
for some constant $c=c(n,s_0,\Lambda,q)$.
In light of standard scaling and covering arguments, we arrive at the estimate
\begin{align*}
    R^{2\gamma}[u]_{L^\infty(I^s_{R/2}(t_0);C^{0,2\gamma}(B_{R/2}(x_0)))}\leq c\widetilde{E}^{q}(u;Q^s_{R}(z_0)),
\end{align*}
where $c=c(n,s_0,\Lambda,q)$, whenever $Q^s_R(z_0)\Subset \Omega_T$. By iterating the above procedure $i_0-1$ times, we deduce \eqref{ineq.gra.est}.

We are now ready to prove the H\"older continuity of $\nabla u$. By Lemma \ref{lem.hol}, we first observe
\begin{align*}
\left[\frac{\delta_h w}{|h|}\right]_{C^{0,\gamma}(Q^s_{1/2})}\leq c{E}^q\left(\frac{\delta_h w}{|h|};Q^s_1\right)+c\left(\dashint_{I^s_1}[f(\cdot,t)]^q_{C^{0,1}(B_{1})}\,dt\right)^{\frac1q}
\end{align*}
for any $h\in B_{1/100}\setminus\{0\}$, where $c=c(n,s_0,\Lambda,q)$. As in the estimate of $J$, we further estimate the right-hand side of the above inequality as  
\begin{align*}
    \left[\frac{\delta_h w}{|h|}\right]_{C^{0,\gamma}(Q^s_{1/2})}&\leq c[u_R]_{L^q(I^s_1;C^{0,1}(B_4))}+c\widetilde{E}^q(u_R;Q^s_4)
\end{align*}
for some constant $c=c(n,s_0,\Lambda,q)$, where we have used \eqref{esth.loc.tail}. Let us now choose $h=\overline{h}e_i$, where $\overline{h}\in(0,1/100)$ and $e_i$ ($i=1,2,\dots,n$) are the canonical unit vectors in $\mathbb{R}^n$.
By taking $\overline{h}\to 0$, we deduce
\begin{align*}
    [\nabla w]_{C^{0,\gamma}(Q^s_{1/2})}\leq c[u_R]_{L^q(I^s_1;C^{0,1}(B_4))}+c \widetilde{E}^q(u_R;Q^s_4).
\end{align*}
By a scaling argument together with \eqref{ineq.gra.est}, we see that
\begin{align*}
    R^{\gamma}[\nabla u]_{C^{0,\gamma}(Q^s_{R/2}(z_1))}\leq c \widetilde{E}^q(u/R;Q^s_{8R}(z_1)).
\end{align*}
By standard covering arguments and the fact that $u-(u)_{Q^s_R(z_0)}$ is also a weak solution to \eqref{eq:nonlocaleq} with $\mu=0$, we get 
\begin{align*}
    R^{\gamma}[\nabla u]_{C^{0,\gamma}(Q^s_{R/2}(z_0))}\leq c {E}^q(u/R;Q^s_{R}(z_0))
\end{align*}
for some constant $c=c(n,s_0,\Lambda,q)$, whenever $Q^s_R(z_0)\Subset \Omega_T$. This completes the proof.
\end{proof}

We end this section with two lemmas. The first one yields oscillation estimates for the difference quotients $\frac{\delta_hv}{|h|}$, where $v$ is a weak solution to \eqref{eq:nonlocaleq} with $\mu=0$.
\begin{lemma}\label{lem.exc.hom}
    Let $v$ be a weak solution to 
    \begin{align}\label{eq.exc.hom}
        \partial_t v+\mathcal{L}v=0\quad\text{in }Q^s_{2R}(z_0).
    \end{align}
    For any $0<|h|<R/100$, we have 
     \begin{align*}
        &\left\|\frac{\delta_h v}{|h|}-\left(\frac{\delta_h v}{|h|}\right)_{Q^s_R(z_0)}\right\|_{L^\infty(Q^s_{R/2}(z_0))}\leq cE\left(\frac{\delta_h v}{|h|};Q^s_R(z_0)\right)
    \end{align*}
    for some constant $c=c(n,s_0,\Lambda)$.
    If $v\in L^q(I^s_{R}(t_0);L^1_{2s}(\bbR^n))$ with $q>1$,
    then there is a constant $\gamma=\gamma(n,s_0,\Lambda,q)$ such that
    \begin{align*}
       R^{\gamma}\left[\frac{\delta_h v}{|h|}\right]_{C^{\gamma}(Q^s_{R/2}(z_0))}\leq cE^q\left(\frac{\delta_h v}{|h|};Q^s_R(z_0)\right)
    \end{align*}
    for some constant $c=c(n,s_0,\Lambda,q)$.
\end{lemma}
\begin{proof}
    Since in view of Lemma \ref{lem.linsol}, $\frac{\delta_h v}{|h|}-\left(\frac{\delta_h v}{|h|}\right)_{Q^s_R(z_0)}$ is a weak solution to \eqref{eq.gra.dfq} in $Q^s_R(z_0)$ with $f=0$, the desired estimates follow from Lemma \ref{lem.lin.hol}.
\end{proof}
Next, we provide a lemma which yields H\"older estimates involving affine functions and will be an essential tool to obtain fractional Sobolev regularity of the gradient of weak solutions to \eqref{eq.exc.hom} with respect to the time variable.

\begin{lemma}\label{lem.linmsol} Let $s_0>1/2$, $s\in[s_0,1)$, and $v$ be a weak solution to \eqref{eq.exc.hom}. For any affine function $l=A\cdot x+b$ with $A\in \bbR^{n}$ and $ b\in \bbR$, we have
\begin{equation}\label{bdd.lin.ftn}
\begin{aligned}
    \|v-l\|_{L^\infty(Q^s_{R/2}(z_0))}\leq  cE(v-l;Q^s_{R}(z_0)),
\end{aligned}
\end{equation}
where $c=c(n,s_0,\Lambda)$. If $v\in L^q(I^s_{R}(t_0);L^1_{2s}(\bbR^n))$ with $q>1$, then there is constant $\gamma=\gamma(n,s_0,\Lambda,q)$ such that
\begin{equation}\label{hol.lin.ftn}
\begin{aligned}
    R^{\gamma}[v-l]_{C^{\gamma}(Q^s_{R/2}(z_0))}&\leq cE^q(v-l;Q^s_{R}(z_0))
\end{aligned}
\end{equation}
for some constant $c=c(n,s_0,\Lambda,q)$.
\end{lemma}
\begin{proof}
Let us fix a affine function $l(x)=A\cdot x+b$. Since $s \in (1/2,1)$, in view of \cite[Remark 3.4]{KuuSimYan22} we have
$l\in W^{s,2}_{\mathrm{loc}}(\bbR^n)\cap L^1_{2s}(\bbR^n)$. In addition, 
\begin{align*}
	\partial_tl+\mathcal{L}l=0\quad\text{in }Q^s_{2R}(z_0)
\end{align*}
holds weakly, as $\Phi$ in \eqref{eq:nonlocalop} is an odd function.
As in the proof of Lemma \ref{lem.lo.hol}, we obtain $v-l$ is a weak solution to
\begin{align}\label{eq.solmlin}
    \partial_t (v-l)+\widetilde{\mathcal{L}}_{K} (v-l)=0,
\end{align}
where 
\begin{align*}
    K(x,y,t)=\frac{\Phi\left(\frac{v(x,t)-v(y,t)}{|x-y|^s}\right)-\Phi\left(\frac{l(x,t)-l(y,t)}{|x-y|^s}\right)}{\frac{v(x,t)-v(y,t)}{|x-y|^s}-\frac{l(x,t)-l(y,t)}{|x-y|^s}}|x-y|^{-n-2s}
\end{align*}
satisfies \eqref{lin.kern.cond} and \eqref{lin.kern.cond2}. In light of Lemma \ref{lem.lin.hol}, we obtain \eqref{hol.lin.ftn}.
\end{proof}

\section{Comparison estimates and existence of SOLA} \label{sec3}
In this section, we establish zero-order comparison estimates that will be used crucially in the remainder of the paper. As a first application, we will use them to deduce the existence of SOLA to parabolic nonlinear nonlocal initial boundary value problems, that is, to prove Theorem \ref{thm: existence}.

\subsection{Zero-order comparison estimates}
We now prove robust comparison estimates which are stable as $s\to1$.
\begin{lemma}\label{lem.comp}
	Let $s_0 \in (0,1)$, $s \in [s_0,1)$, fix $p\in\left[1,\frac{n+2s_0}{n+s_0}\right)$ and let $$u\in L^{2}(I^s_{2R}(t_0);W^{s,2}(B_{2R}(x_0)))\cap C(I^s_{2R}(t_0);L^2(B_{2R}(x_0))) \cap L^p(I^s_{2R}(x_0);L^{1}_{2s}(\bbR^n))$$ be a weak solution to 
	\begin{equation*}
		\partial_t u+\mathcal{L}u=\mu\quad\text{in }Q^s_{2R}(z_0)
	\end{equation*}with $\mu\in L^1(I^s_{2R}(t_0);L^\infty(B_{2R}(x_0)))$. Then there exists a unique weak solution 
	\begin{equation*}
		v\in L^{2}(I^s_R(t_0);W^{s,2}(B_{R}(x_0)))\cap C(I^s_R(t_0);L^2(B_R(x_0))) \cap L^p(I^s_R(x_0);L^{1}_{2s}(\bbR^n))
	\end{equation*}to 
	\begin{equation}\label{eq.ivp}
		\left\{
		\begin{alignedat}{3}
			\partial_t v+\mathcal{L}v&= 0&&\qquad \mbox{in  $Q^s_{R}(z_0)$}, \\
			v&=u&&\qquad  \mbox{in $(\mathbb{R}^{n}\setminus B_{R}(x_0))\times I^s_R(t_0)$}, \\
			v(\cdot,t_0-R^{2s})&=u(\cdot,t_0-R^{2s}) &&\qquad \mbox{in  $B_R(x_0)$}
		\end{alignedat} \right.
	\end{equation}
	such that 
	\begin{align*}
		\sup_{t\in I^s_R(t_0)}\dashint_{B_R(x_0)}|(u-v)(x,t)|\,dx\leq cR^{-n}|\mu|(Q^s_R(z_0))
	\end{align*}
	holds for some constant $c=c(n,s_0,\Lambda)$.
	In addition, we have 
	\begin{align*}
		\left(\dashint_{Q^s_R(z_0)}|u-v|^p\,dz\right)^{\frac1p}\leq cR^{-n}|\mu|(Q^s_R(z_0))
	\end{align*}
	for some constant $c=c(n,s_0,\Lambda,p)$. 
\end{lemma}
\begin{remark}\label{rmk.exi} \normalfont
	We remark that the existence of the weak solution $v$ to \eqref{eq.ivp} is obtained in \cite[Lemma A.1]{ByuKimKim23} when $p=2$. Nevertheless, an inspection of the proof reveals that this existence result remains valid for $p\in[1,\infty]$. 
\end{remark}
\begin{proof}
	We may assume $R=1$ and $z_0=0$. Let us fix $p\in\left[1,\frac{n+2s_0}{n+s_0}\right)$. 
	By Remark \ref{rmk.exi}, there is a unique weak solution 
	\begin{equation*}
		v\in C(I^s_1;L^2(B_1))\cap  L^{2}(I^s_1;W^{s,2}(B_{2}))\cap L^p(I^s_1;L^{1}_{2s}(\bbR^n))
	\end{equation*}
	to \eqref{eq.ivp} with $R=1$ and $z_0=0$. Let us denote $w=u-v$ \and consider
	\begin{equation*}
		\varphi^{\pm}_{1,\epsilon}=\pm\min\{1,w_{\pm}/\epsilon\}.
	\end{equation*}
	Since 
	\begin{align*}
		|\min\{1,g_\pm\}(x)-\min\{1,g_\pm\}(y)|\leq |g_+(x)-g_+(y)|
	\end{align*}
	holds, we have $\varphi^{\pm}_{1,\epsilon}\in L^2(I^s_1;W^{s,2}(B_1))$.
	By following the same lines as in the proof of \cite[Lemma 3.1]{NguNowSirWei23} along with the standard mollification, we deduce 
	\begin{equation}\label{comp.ineq1}
		\begin{aligned}
			\sup_{t\in I^s_1}\int_{B_1}|w|\,dz       \leq c|\mu|(Q^s_1)
		\end{aligned}
	\end{equation}
	for some constant $c=c(n,\Lambda)$. In addition, we obtain
	\begin{equation}\label{comp.ineq2}
		\begin{aligned}
			(1-s)\int_{I^s_1}\int_{\bbR^n}\int_{\bbR^n}\left(\Phi(d_{x,y}^su)-\Phi(d_{x,y}^sv)\right)\frac{\phi^{\pm}_{1,\epsilon}(x,t)-\phi^{\pm}_{1,\epsilon}(y,t)}{|x-y|^{n+s}}\,dy\,dz\leq c|\mu|(Q^s_1)
		\end{aligned}
	\end{equation}
	for some constant $c=c(n,\Lambda)$, where we denote 
	\begin{align*}
		d^s_{x,y}u=\frac{u(x,t)-u(y,t)}{|x-y|^s}\quad\text{and}\quad d^s_{x,y}v=\frac{v(x,t)-v(y,t)}{|x-y|^s}.
	\end{align*}
	Let us fix the constant $\sigma\in(0,s)$ satisfying
	\begin{equation}\label{cond.sigma0}
		\frac{n+2\sigma}{n+\sigma}>p.
	\end{equation}    
	We are now going to prove
	\begin{equation}\label{ineq.zero.wp}
		\begin{aligned}
			&(1-\sigma)\int_{I^s_1}\int_{B_2}\int_{B_2}\frac{|w(x,t)-w(y,t)|^p}{|x-y|^{n+\sigma p}}\,dy\,dx\,dt\\
			&\leq c\left(\frac{(1-\sigma)}{(1-s)^{\frac{p}2}(s-\sigma)^{\frac{2-p}2}}\right)^{\frac{2(n+\sigma)}{n+2\sigma}}(|\mu|(Q^s_1))^p
		\end{aligned}
	\end{equation}
	with some constant $c=c(n,\Lambda)$. For $m>1$, let us write
	\begin{equation*}
		\varphi^{\pm}_{2,\epsilon}=\pm\left(d^{1-m}-(d+w_{\pm})^{1-m}\right)\varphi^{\pm}_{1,\epsilon}.
	\end{equation*}
	Since
	\begin{equation*}
		\pm\left(d^{1-m}-(d+w_{\pm})^{1-m}\right),\,\varphi^{\pm}_{1,\epsilon}\in L^\infty(Q^s_1)\cap L^2(I^s_1;W^{s,2}(B_1)),
	\end{equation*}
	we get $\varphi^{\pm}_{2,\epsilon}\in L^2(I^s_1;W^{s,2}(B_1))$. Using a standard mollification argument (see \cite[Lemma 4.1]{KuuMin14w}), we deduce 
	\begin{align*}
		-\int_{Q^s_1}w \partial_t\varphi^{\pm}_{2,\epsilon}\,dz=\int_{Q^s_1}\partial_t\left(\int_{0}^{w_{\pm}}\left(d^{1-m}-(d+\tau)^{1-m}\right){\min\{1,\tau/\epsilon\}}\,d\tau\right) \,dz.
	\end{align*}
	Using \eqref{comp.ineq1}, we obtain
	\begin{align*}
		\left|-\int_{Q}w \partial_t\varphi^{\pm}_{2,\epsilon}\,dz\right|\leq \sup_{t\in I^s_1}d^{1-m}\int_{B_1}|w|\,dx\leq cd^{1-m}|\mu|(Q^s_1).
	\end{align*}
	We next observe that
	\begin{align*}
		&\int_{I^s_1}\int_{\bbR^n}\int_{\bbR^n}\left(\Phi(d_{x,y}^su)-\Phi(d_{x,y}^sv)\right)\frac{\phi^{\pm}_{2,\epsilon}(x,t)-\phi^{\pm}_{2,\epsilon}(y,t)}{|x-y|^{n+s}}\,dy\,dz\\
		&=\int_{I^s_1}\int_{\bbR^n}\int_{\bbR^n}\left(\Phi(d_{x,y}^su)-\Phi(d_{x,y}^sv)\right)\frac{(\zeta(x,t)-\zeta(y,t))\phi^{\pm}_{1,\epsilon}(x,t)}{|x-y|^{n+s}}\,dy\,dz\\
		&\quad+ \int_{I^s_1}\int_{\bbR^n}\int_{\bbR^n}\left(\Phi(d_{x,y}^su)-\Phi(d_{x,y}^sv)\right)\frac{(\phi^{\pm}_{1,\epsilon}(x,t)-\phi^{\pm}_{1,\epsilon}(y,t))\zeta(y,t)}{|x-y|^{n+s}}\,dy\,dz\\
		&\eqqcolon J_1+J_2,
	\end{align*}
	where we denote 
	\begin{align*}
		\zeta=\pm\left(d^{1-m}-(d+w_{\pm})^{1-m}\right).
	\end{align*}
	By letting $\epsilon\to 0$, we obtain $\phi_{1,\epsilon}^{\pm}=1$. By following the same lines as in the proof of \cite[Lemma 3.1]{KuuMinSir15}, we then get
	\begin{align*}
		\lim_{\epsilon\to0}J_1&=\int_{I^s_1}\int_{\bbR^n}\int_{\bbR^n}\left(\Phi(d_{x,y}^su)-\Phi(d_{x,y}^sv)\right)\frac{(\zeta(x,t)-\zeta(y,t))}{|x-y|^{n+s}}\,dy\,dz\\
		&\geq \frac{1}{c}\int_{I^s_1}\int_{B_2}\int_{B_2}\frac{|w(x,t)-w(y,t)|^2}{(d+|w(x)|+|w(y)|)^m}\frac{\,dy\,dz}{|x-y|^{n+2s}}
	\end{align*}
	for some constant $c=c(n,\Lambda)$.
	We next observe from \eqref{comp.ineq2} that 
	\begin{align*}
		|J_2|\leq c(1-s)^{-1}d^{1-m}|\mu|(Q^s_1).
	\end{align*}

	Combining all the above estimates for $J_1$ and $J_2$ yields
	\begin{align*}
		(1-s)\int_{I^s_1}\int_{B_2}\int_{B_2}\frac{|w(x,t)-w(y,t)|^2}{(d+|w(x)|+|w(y)|)^m}\frac{\,dy\,dz}{|x-y|^{n+2s}}\leq cd^{1-m}|\mu|(Q^s_1)
	\end{align*}
	for some constant $c=c(n,\Lambda)$.
	Using H\"older's inequality, we obtain
	\begin{align*}
		L_1&\coloneqq(1-\sigma)\int_{I^s_1}\int_{B_2}\int_{B_2}\frac{|w(x,t)-w(y,t)|^p}{|x-y|^{n+\sigma p}}\,dy\,dz\\
		&\leq (1-\sigma)\left(\int_{I^s_1}\int_{B_2}\int_{B_2}\frac{|w(x,t)-w(y,t)|^2}{(d+|w(x)|+|w(y)|)^m}\frac{\,dy\,dz}{|x-y|^{n+2s}}\right)^{\frac{p}2}\\
		&\quad\times \left(\int_{I^s_1}\int_{B_2}\int_{B_2}(d+|w(x,t)|+|w(y,t)|)^{\frac{mp}{2-p}}\frac{\,dy\,dz}{|x-y|^{n-\frac{2p(s-\sigma)}{2-p}}}\right)^{\frac{2-p}2}\\
		&\leq c(1-\sigma)\left(\frac{d^{1-m}|\mu|(Q^s_1)}{1-s}\right)^{\frac{p}2}\left(\frac{1}{s-\sigma}\int_{I^s_1}\int_{B_2}(d+|w(z)|)^{\frac{mp}{2-p}}\,dz\right)^{\frac{2-p}2},
	\end{align*}
	where $c=c(n,\Lambda)$. We now choose 
	\begin{equation*}
		m=(2-p)\frac{n+\sigma}{n}\quad\text{and}\quad d=\left(\int_{Q^s_2}|w(z)|^{\frac{p(n+\sigma)}{n}}\,dz\right)^{\frac{n}{p(n+\sigma)}}
	\end{equation*}
	to see that 
	\begin{align*}
		d\leq c\left((1-\sigma)\int_{I^s_1}\int_{B_2}\int_{B_2}\frac{|w(x,t)-w(y,t)|^p}{|x-y|^{n+\sigma p}}\,dy\,dz\right)^{\frac{n}{p(n+\sigma)}}
		\left(\sup_{t\in I^s_1}\int_{B_2}|w|\,dx\right)^{\frac{\sigma}{n+\sigma}}
	\end{align*}
	for some constant $c=c(n,\sigma,p)$, where we have used Lemma \ref{lem.inter}.
	In addition, by \eqref{cond.sigma0}, we see that $m>1$.
	Combining the above two inequalities with \eqref{comp.ineq1}, we obtain
	\begin{align*}
		L_1\leq \frac{c(1-\sigma)\left(|\mu|(Q^s_1)\right)^{\frac{p(n+2\sigma)}{2(n+\sigma)}}}{(1-s)^{\frac{p}2}(s-\sigma)^{\frac{2-p}2}}
		\left((1-\sigma)\int_{I^s_1}\int_{B_2}\int_{B_2}\frac{|w(x,t)-w(y,t)|^p}{|x-y|^{n+\sigma p}}\,dy\,dz\right)^{\frac{n}{2(n+\sigma)}} .
	\end{align*}
	Thus, we have proved \eqref{ineq.zero.wp}
	and 
	\begin{align}\label{est.comp.last}
		\int_{Q^s_1}|w|^p\,dz\leq c\left(\frac{(1-\sigma)}{(1-s)^{\frac{p}2}(s-\sigma)^{\frac{2-p}2}}\right)^{\frac{2(n+\sigma)}{n+2\sigma}}|\mu|(Q^s_1)^p
	\end{align}
	for some constant $c=c(n,\sigma,p)$, where we have used the fractional Poincar\'e-type inequality from \cite[Corollary 4.9]{Coz17}. In addition, the constant $c$ depends only on $n,\sigma_0$ and $p$, whenever $\sigma\in[\sigma_0,1)$. 
	By choosing $\sigma$ in a suitable way, we now provide stable estimates with respect to the parameter $s$.
	We first consider the case when $p=1$. In this case we choose $\sigma=2s-1$ to see that 
	\begin{align}\label{sta.cons}
		\frac{(1-\sigma)}{(1-s)^{\frac{p}2}(s-\sigma)^{\frac{2-p}2}}\leq c
	\end{align}
	for some constant $c=c(n,s_0)$.
	
	Let us consider $p>1$ and choose
	\begin{equation*}
		s_p=\frac{n(p-1)}{2-p}
	\end{equation*}
	to see that
	\begin{align*}
		p=\frac{n+2s_p}{n+s_p}.
	\end{align*}
	As
	\begin{equation*}
		f(a)= \frac{n+2a}{n+a} \quad\text{for any }a\in[0,\infty)
	\end{equation*}
	is an increasing function, we conclude that
	\begin{equation}\label{cond.sigma}
		\begin{aligned}
			\sigma=\begin{cases}
				\frac{ s_p+s_0}{2}&\quad\text{if }s_0\leq s\leq \frac{s_p+s_0+1}{3}\\
				(3s-1)/2&\quad\text{if }s\geq \frac{s_p+s_0+1}{3}
			\end{cases}
		\end{aligned}
	\end{equation}
	satisfies \eqref{cond.sigma0}. Moreover, we conclude that \eqref{sta.cons} holds with respect to some constant $c=c(n,s_0,p)$.
	Since $\sigma\geq (3s_0-1)/2$, plugging \eqref{sta.cons} into \eqref{est.comp.last} now yields the desired stable estimates, completing the proof.
\end{proof}

\subsection{Existence of SOLA}

\begin{proof}[Proof of Theorem \ref{thm: existence}]
	Since $f(a)=\frac{n+2a}{n+a}$ for $a\geq0$ is an increasing function, we find constants $\sigma_0\in(\sigma,s)$ such that
	\begin{equation}\label{rel.p0sig0}
		p_0<\frac{n+2\sigma_0}{n+\sigma_0},
	\end{equation}
	where we denote 
	\begin{equation}\label{choi.p0}
		p_0=\frac12\left(p+\frac{n+2s}{n+s}\right)>1.
	\end{equation}
	We now divide the proof into the three steps. 
	
	\textbf{Step 1: Regularizing the given data}. By a standard approximation via mollifiers, there are sequences $\{g_{i}\}\subset C^{\infty}_{c}\left(\mathbb{R}^{n}\times[0,T]\right)$ and $\{\mu_{i}\}\subset C^{\infty}_{c}\left(\mathbb{R}^{n}\times[0,T]\right)$ such that 
	\begin{equation}
		\label{sola : regg}
		\left\{
		\begin{alignedat}{3}
			&g_{i}\to g&&\quad\text{in }L^{2}\left(0,T;W^{s,2}_{\mathrm{loc}}(\mathbb{R}^{n})\right)\cap L^{1}\left(0,T;L^{1}_{2s}(\mathbb{R}^{n})\right), \\
			&\partial_t g_{i}\to \partial_t g&&\quad\text{in }\big(L^{2}\left(0,T;W^{s,2}(\Omega)\right)\big)^{*},\\
			&g_{i}(\cdot,0)\to g_{0}&&\quad\text{in }L^{2}(\Omega)
		\end{alignedat} \right.
	\end{equation}
	and 
	\begin{equation}
		\label{sola : regmu}
		\left\{
		\begin{alignedat}{3}
			&\mu_{i}\rightharpoonup \mu \text{ in }\mathcal{M}(\mathbb{R}^{n}\times(0,T)), \\
			&\limsup_{i\to\infty}|\mu_i|(Q)\leq |\mu|(\bar{Q})
		\end{alignedat} \right.
	\end{equation}
	for all $Q\subset\mathbb{R}^{n}\times(0,T)$. By following the same lines as in the proof of \cite[Lemma A.1]{ByuKimKim23}, there exists a solution $u_i\in C([0,T];L^{2}(\Omega))\cap L^{2}\left(0,T;W^{s,2}(\mathbb{R}^{n})\right)$ to \eqref{thm eq : regularized1} and $v_i\in C([0,T];L^{2}(\Omega))\cap L^{2}\left(0,T;W^{s,2}(\mathbb{R}^{n})\right)$ to \eqref{thm eq : regularized1} with $\mu_i=0$, respectively.
	As in the estimate of $I$ given in \cite[Equation (4.9)]{ByuKimKum24} with $w_i=v_i$ ignoring the gradient term, we have 
	\begin{equation}\label{ineq.1.app}
		\begin{aligned}
			&\sup_{t\in (0,T)}\int_{\Omega}|v_i-g_i|^2\,dx+[v_i-g_i]_{L^2(0,T;W^{s,2}(\bbR^n))}\\
			&\leq c\|\partial_t g_i\|^2_{\left(L^2(0,T;W^{s,2}(\Omega))\right)^{*}}\\
			&\quad+c\int_{0}^T\int_{\bbR^n}\int_{\bbR^n}\frac{|g_i(x,t)-g_i(y,t)||(g_i-v_i)(x,t)-(g_i-v_i)(y,t)|}{|x-y|^{n+2s}}\,dx\,dy\,dt\eqqcolon L
		\end{aligned}
	\end{equation}
	for some constant $c=c(n,\Lambda)$. Moreover, by the estimate of $J$ determined in \cite[Equation (4.9)]{ByuKimKum24}, we estimate $L$ as
	\begin{align*}
		L&\leq 1/2\sup_{t\in (0,T)}\int_{\Omega}|v_i-g_i|^2\,dx+1/2[v_i-g_i]^2_{L^2(0,T;W^{s,2}(\bbR^n))}\\
		&\quad+ c\|g_i\|^2_{L^2(0,T;W^{s,2}(B_R))}+c\|g_i\|^2_{L^1(0,T;L^1_{2s}(\bbR^n))}
	\end{align*}
	for some constant $c=c(n,s,\Lambda,\Omega,R)$, where we choose a sufficiently large ball $B_R$ such that $\Omega\Subset B_{R/2}$. Therefore, plugging this into \eqref{ineq.1.app} together with a few calculations, we have 
	\begin{equation}\label{ineq.2.app}
		\begin{aligned}
			&\sup_{t\in (0,T)}\int_{\Omega}|v_i|^2\,dx+[v_i]^2_{L^2(0,T;W^{s,2}(B_R))}\\
			&\leq c\|\partial_t g_i\|^2_{\left(L^2(0,T;W^{s,2}(\Omega))\right)^{*}}+c\|g_i\|^2_{L^2(0,T;W^{s,2}(B_R))}+c\|g_i\|^2_{L^1(0,T;L^1_{2s}(\bbR^n))},
		\end{aligned}
	\end{equation}
	where $c=c(n,s,\Lambda,\Omega)$.
	
	\textbf{Step 2: A priori estimates for approximating solutions.}
	By Lemma \ref{lem.comp}, \eqref{ineq.zero.wp} together with \eqref{rel.p0sig0}, H\"older's inequality and \eqref{ineq.2.app}, we have 
	\begin{align*}
		&\sup_{t\in (0,T)}\int_{\Omega}|u_i|\,dx+[u_i]_{L^{p_0}(0,T;W^{\sigma_0,p_0}(B_R))}\\
		&\leq \sup_{t\in (0,T)}\int_{\Omega}|v_i|\,dx+[v_i]_{L^{p_0}(0,T;W^{\sigma_0,p_0}(B_R))}+c|\mu_i|(\Omega_T)\\
		&\leq \left(\sup_{t\in (0,T)}\int_{\Omega}|v_i|^2\,dx\right)^{\frac12}+[v_i]_{L^2(0,T;W^{s,2}(B_R))}+c|\mu_i|(\Omega_T)\\
		&\leq c\|\partial_t g_i\|_{\left(L^2(0,T;W^{s,2}(\Omega))\right)^{*}}+c\|g_i\|_{L^2(0,T;W^{s,2}(B_R))}+c\|g_i\|_{L^1(0,T;L^1_{2s}(\bbR^n))}\\
		&\quad+c|\mu_i|(\Omega_T)
	\end{align*}
	for some constant $c=c(n,s,\Lambda,p,\sigma,\Omega,R)$.
	We next note 
	\begin{align*}
		\int_{0}^T\int_{\bbR^n}\frac{|u_i(x,t)|}{(1+|x|)^{n+2s}}\,dz\leq \sup_{t\in (0,T)}\int_{\Omega}|u_i|\,dx+c\|g_i\|_{L^1(0,T;L^1_{2s}(\bbR^n))}.
	\end{align*}
	In addition, by the Sobolev inequality given in \cite[Corollary 4.9]{Coz17}, we have 
	\begin{align*}
		\|u_i-g_i\|_{L^{p_0}(0,T;W^{\sigma_0,p_0}(B_R))}\leq [u_i-g_i]_{L^{p_0}(0,T;W^{\sigma_0,p_0}(B_R))}.
	\end{align*}
	Thus, combining the above three estimates along with \eqref{sola : regg} and \eqref{sola : regmu} yields
	\begin{equation}\label{uiest}
		\|u_i\|_{L^\infty(0,T;L^1(\Omega))}+\|u_i\|_{L^{p_0}(0,T;W^{\sigma_0,p_0}(B_R))}+\|u_i\|_{ L^1(0,T;L^1_{2s}(\bbR^n))}\leq c,
	\end{equation}
	where $c$ depends only on the given data. Let us consider the space 
	\begin{equation}\label{defn.space.x0}
		X_0^{2s-\sigma_0,p_0'}(\Omega,B_R)=\{f \in W^{2s-\sigma_0,p_0'}(B_R) \mid f =0\text{ on }\bbR^n\setminus \Omega\},
	\end{equation}
	which was introduced in \cite[Appendix A]{BLSt}. We also define a norm on $X_0^{2s-\sigma_0,p_0'}(\Omega,B_R)$ as follows
	\begin{equation*}
		\|f\|_{X_0^{2s-\sigma_0,p_0'}(\Omega,B_R)}\coloneqq\|f\|_{W^{2s-\sigma_0,p_0'}(B_R)}.
	\end{equation*}
	We now observe that for any smooth function $\phi$ with $\|\phi\|_{L^{\infty}(0,T;X_0^{2s-\sigma_0,p_0'}(\Omega,B_R))}\leq 1$,
	\begin{align*}
		&\left|\int_{0}^T\int_{\bbR^n}\int_{\bbR^n}\Phi\left(\frac{u_i(x,t)-u_i(y,t)}{|x-y|^s}\right)\frac{\phi(x,t)-\phi(y,t)}{|x-y|^{n+s}}\,dx\,dy\,dt\right|\\
		&\leq \int_{0}^T\int_{\overline{\Omega}}\int_{B_R}\left|\Phi\left(\frac{u_i(x,t)-u_i(y,t)}{|x-y|^s}\right)\frac{\phi(x,t)-\phi(y,t)}{|x-y|^{n+s}}\right|\,dx\,dy\,dt\\
		&\quad+ c\int_{0}^T\int_{\bbR^n\setminus B_R}\int_{\Omega}\left|\Phi\left(\frac{u_i(x,t)-u_i(y,t)}{|x-y|^s}\right)\frac{\phi(x,t)}{|x-y|^{n+s}}\right|\,dx\,dy\,dt\\
		&\leq c[u_i]_{L^{p_0}(0,T;W^{\sigma_0,p_0}(B_R))}[\phi]_{L^{p_0'}(0,T;W^{2s-\sigma,p_0'}(B_R))}\\
		&\quad+ c\left[\|u_i\|_{L^{p_0}(0,T;L^{p_0}(B_R))}+\|u_i\|_{L^1(0,T;L^1_{2s}(\bbR^n))}\right]\|\phi\|_{L^{\infty}(0,T;L^{p_0'}(B_R))}\\
		&\leq c\left([u_i]_{L^{p_0}(0,T;W^{\sigma_0,p_0}(B_R))}+\|u_i\|_{L^{p_0}(0,T;L^{p_0}(B_R))}+\|u_i\|_{L^1(0,T;L^1_{2s}(\bbR^n))}\right).
	\end{align*}
	Since for any $\phi\in C^1_{0}(\Omega_T)$, we have 
	\begin{align*}
		&\int_{\Omega_T}u_i\partial_t\phi\,dz\\
		&=(1-s)\int_{0}^T\int_{\bbR^n}\int_{\bbR^n}\Phi\left(\frac{u_i(x,t)-u_i(y,t)}{|x-y|^s}\right)\frac{\phi(x,t)-\phi(y,t)}{|x-y|^{n+s}}\,dx\,dy\,dt
		+\int_{\Omega_T}\mu_i\phi\,dz,
	\end{align*}
	we get 
	\begin{align*}
		\partial_t u_i=f_i+\mu_i
	\end{align*}
	in the distributional sense, where $\{f_i\}\subset L^1\left(0,T;\left(X_0^{2s-\sigma_0,p_0'}(\Omega,B_R)\right)^*\right)$ and $\{\mu_i\}\subset L^1(0,T;L^1(\bbR^n))$ are uniformly bounded in each space, respectively. We now define
	\begin{equation*}
		p_1 \coloneqq \min \left \{\frac{n}{n-(2s-\sigma_0)},p_0 \right \}, \quad p_1^\prime \coloneqq \max \left \{\frac{n}{2s-\sigma_0},p_0^\prime \right \}
	\end{equation*}
	to see that in view of the fractional Sobolev embedding we have
	\begin{equation*}
		X_0^{2s-\sigma_0/2,p_1^\prime}(\Omega,B_R)\subset L^\infty(B_R),
	\end{equation*}
	which implies
	\begin{align*}
		\mu_i \in L^1\left(0,T;\left(X_0^{2s-\sigma_0/2,p_1^\prime}(\Omega,B_R)\right)^*\right).
	\end{align*}
	By \cite[Proposition 2.1]{Now21}, we also see $X_0^{2s-\sigma_0/2,p_1^\prime}(\Omega,B_R)\subset X_0^{2s-\sigma_0,p_0'}(\Omega,B_R)$. Therefore, the sequence $\{\partial_t u_i\}$ is uniformly bounded in $L^1\left(0,T;\left(X_0^{2s-\sigma_0/2,p_1^\prime}(\Omega,B_R)\right)^*\right)$, as $\Omega\subset B_R$. We now use a compactness result given in \cite[Section 8]{Sim87} with $q=p_0$, $X=W^{\sigma_0,p_0}(B_R)$, $B=L^1(B_R)$ and $Y=\left(X_0^{2s-\sigma_0/2,p_1^\prime}(\Omega,B_R)\right)^*$ to see that 
	\begin{align}\label{conv1}
		u_i\to u \quad\text{in }L^1(0,T;L^1(B_R))\quad\text{and}\quad \lim_{i\to\infty}u_i(z)=u(z)\quad\text{for a.e. } z\in B_R\times (0,T).
	\end{align}
	Applying Fatou's lemma to \eqref{uiest}, we have 
	\begin{align*}
		u\in L^\infty(0,T;L^1(\Omega))\cap L^{p_0}(0,T;W^{\sigma_0,p_0}(B_R)).
	\end{align*}
	In addition, since $u\in L^\infty(0,T;L^1(\Omega))$ and $g\in L^q(0,T;L^1_{2s}(\bbR^n))$, we obtain that $u\in L^q(0,T;L^1_{2s}(\bbR^n))$.
	By following the same lines as in Step 3 of the proof of \cite[Theorem 1.1]{KuuMinSir15} along with \eqref{conv1} and \eqref{sola : regmu}, we deduce that $u$ satisfies \eqref{defn.sola}.  
	In addition,  \eqref{conv1} along with \eqref{sola : regg} and \eqref{sola : regmu} yields \eqref{regular limit}.
	
	\textbf{Step 3: Initial data.} In essentially the same way as in Step 3 of the proof of \cite[Theorem 4.2]{ByuKimKum24}, we obtain \eqref{defn sola : initial data}. 
	Therefore, by \cite[Proposition 2.1]{Now21} together with the facts that $\sigma<\sigma_0$ and $p<p_0$, we have 
	\begin{align*}
		u\in L^\infty(0,T;L^1(\Omega))\cap L^{p}(0,T;W^{\sigma,p}(B_R))\cap L^q(0,T;L^1_{2s}(\bbR^n)),
	\end{align*}
	which completes the proof.
\end{proof} 

We end this section with the following remark.
\begin{remark}\label{rmk.relsolaweak}
	We now show that any weak solution 
	\begin{equation*}
		u\in L^2_{\mathrm{loc}}(0,T;W^{s,2}_{\mathrm{loc}}(\Omega))\cap C_{\mathrm{loc}}(0,T;L^2_{\mathrm{loc}}(\Omega))\cap L^1_{\mathrm{loc}}(0,T;L^1_{2s}(\bbR^n))
	\end{equation*}
	to \eqref{eq:nonlocaleq} is also a SOLA in every parabolic cylinder $Q^s_{R}(z_0)\Subset\Omega_T$ such that $Q^s_{2R}(z_0)\Subset\Omega_T$. We may assume $\mu=0$ on $\bbR^{n+1}\setminus \Omega_T$.
	By the standard approximation, there are sequences
	\begin{equation*}
		g_i\in C_c^{\infty}(\Omega_T)\quad\text{and}\quad \mu_i\in C_c^{\infty}(\Omega_T)
	\end{equation*}
	such that 
	\begin{equation*}
		g_i\to u \quad\text{in } L^2(I_{2R}(t_0);W^{s,2}(B_{2R}(x_0)))\cap L^1(I_{2R}(t_0);L^1_{2s}(\bbR^n))
	\end{equation*}
	and \eqref{sola : regmu}.
	Then by \cite[Lemma A.1]{ByuKimKim23}, there exists a solution 
	\begin{equation*}
		u_i\in C(I_{R}(t_0);L^2(B_{R}(x_0)))\cap L^2(I_{R}(t_0);W^{s,2}(\bbR^n))
	\end{equation*} to 
	\begin{equation*}
		\left\{
		\begin{alignedat}{3}
			\partial_tu_i+\mathcal{L}u_{i}&= \mu_{i}&&\qquad \mbox{in  $Q_{R}(z_0)$}, \\
			u_{i}&=g_{i}&&\qquad  \mbox{in $\big(\mathbb{R}^{n}\setminus B_{R}(x_0)\big)\times I_{R}(t_0)$},\\
			u_{i}(\cdot,0)&=g_{0,i}&&\qquad  \mbox{in $B_{R}(x_0)$}.
		\end{alignedat} \right.
	\end{equation*}
	By following the same lines as in the proof of Theorem \ref{thm: existence}, \eqref{regular limit} holds with $\Omega$ and $(0,T)$ replaced by $B_{R}(x_0)$ and $I_{R}(t_0)$, respectively. Therefore, we get that $u$ is SOLA to \eqref{eq: IVP} with $\Omega$, $(0,T)$, $g$ and $g_0$ replaced by $B_{R}(x_0)$, $I_{R}(t_0)$, $u$ and $u(\cdot,t_0-R^{2s})$, respectively.
\end{remark}

\section{Higher differentiability} \label{sec5}
In this section, we prove various higher differentiability estimates for parabolic nonlinear nonlocal measure data problems, which will serve as crucial tools in order to prove our pointwise gradient potential estimates later. For the remainder of this paper, let us fix a parameter $s_0\in (1/2,1)$ and some
\begin{equation}\label{choi.s}
     s\in[s_0,1).
\end{equation}

\subsection{Gradient differentiability}
Using the above comparison estimate together with a parabolic and nonlocal version of the nonlinear atomic decomposition methods pioneered in \cite{KrMin,MinCZMD}, we now prove various gradient estimates for solutions $u$ to \eqref{eq:nonlocaleq}.
\begin{lemma}\label{lem.graex}
    Let us fix $p\in \left(1,\frac{n+2s_0}{n+1}\right)$ and let $$u\in L^{2}(I^s_{2R}(t_0);W^{s,2}(B_{2R}(x_0)))\cap C(I^s_{2R}(t_0);L^2(B_{2R}(x_0))) \cap L^p(I^s_{2R}(x_0);L^{1}_{2s}(\bbR^n))$$ be a weak solution to \eqref{eq:nonlocaleq} in $Q^s_{2R}(z_0)$ with $\mu\in L^1(I^s_{2R}(t_0);L^\infty(B_{2R}(x_0))$.  Then we have
    \begin{align*}
        &\left(\dashint_{Q^s_{R}(z_0)}|\nabla u|^p\,dz\right)^{\frac1p}\leq cE^p(u/R;Q^s_{2R}(z_0))
        +cR^{-(n+1)}|\mu|(Q^s_{2R}(z_0))
    \end{align*}
    for some constant $c=c(n,s_0,\Lambda,p)$.
\end{lemma}
\begin{proof}
In view of a scaling argument based on Lemma \ref{lem:scaling}, we may assume $R=1/2$ and $z_0=0$.
Fix $h\in \bbR^n$ such that $0<|h|\leq s^{\frac1{2s}}/1000$. Let us fix $\beta\in(0,1)$ to be determined later. By Lemma \ref{lem.cov} with $R=1/2$, there are mutually disjoint coverings $\{B_{|h|^{\beta}}(x_{i})\}_{i\in\mathcal{I} }$ of $B_{5/8}$ and $\{I^s_{|h|^\beta}(t_j)\}_{j\in\mathcal{J}}$ of $I^s_{5/8}$ such that for any $k\in\setN$, we have $(x_{i},t_j)\in Q^s_{5/8}$,
\begin{equation}\label{max.cov.ineq00}
    |\mathcal{I}||h|^{n\beta}+|\mathcal{J}||h|^{2s\beta}\leq c,
\end{equation} 
\begin{align}
\label{max.cov.ineq0}
\sup_{x \in \mathbb{R}^{n}} \sum_{i\in \mathcal{I}}\bfchi_{B_{2^{k}|h|^{\beta}}(x_{i})}(x) \leq c2^{nk} 
\end{align}
and
\begin{align}
\label{max.cov.ineq1}
\sup_{t \in \mathbb{R}} \sum_{j\in \mathcal{J}}\bfchi_{I^s_{2^{k}|h|^{\beta}}(t_{j})}(t) \leq c2^{2sk} 
\end{align}
for some constant $c=c(n)$. Moreover, this implies 
\begin{align}\label{max.cov.ineq2}
    \sup_{z \in \mathbb{R}^{n+1}}\sum_{(i,j)\in \mathcal{I}\times\mathcal{J}}\bfchi_{B_{2^{k}|h|^{\beta}}(x_i)\times I^s_{4|h|^\beta}(t_j)}(z) \leq c2^{nk}, 
\end{align}
as for each $z\in \bbR^{n+1}$
\begin{align*}
    &|\{(i,j)\in\mathcal{I}\times\mathcal{J}\,:\, z\in B_{2^{k}|h|^{\beta}}(x_i)\times I^s_{4|h|^\beta}(t_j) \}|\\
    &\leq|\{i\in \mathcal{I}\,:\, x\in B_{2^{k}|h|^{\beta}}(x_i)\}|\times |\{j\in \mathcal{J}\,:\, t\in I^s_{4|h|^{\beta}}(t_j)\}|.
\end{align*}

We now fix a positive integer $m_{0}$ such that 
\begin{equation}
\label{el : choi.m0}
s^{\frac1{2s}}/8\leq2^{m_{0}+4}|h|^{\beta}<s^{\frac1{2s}}/4.
\end{equation}
By Lemma \ref{lem.comp}, there is the weak solution $v_{i,j}$ to \eqref{eq.ivp} with $z_0=(x_i,t_j)$ and $R=4|h|^\beta$
such that
\begin{equation}
\label{maxd.ineq0}
    \left(\dashint_{Q^s_{4|h|^{\beta}}(x_{i},t_j)}|u-v_{i,j}|^p\,dx\right)^{\frac1p}\leq c|h|^{-n\beta}|\mu|(Q^s_{4|h|^{\beta}}(x_{i},t_j))
\end{equation}
for some constant $c=c(n,s_0,\Lambda,p)$.
We first show $u\in L^p_{\mathrm{loc}}(I^s_2;W^{\varsigma,p}_{\mathrm{loc}}(B_2))$ for any $\varsigma\in(0,1)$ with the estimate
\begin{align}\label{maxd.goal1}
    r^{\varsigma}[u]_{L^p(I^s_r(t_0);W^{\varsigma,p}(B_r(x_0)))}\leq cE^p(u;Q^s_{2r}(z_0))+cr^{-n}|\mu|(Q^s_{2r}(z_0)),
\end{align}
where $c=c(n,s_0,\Lambda,p,\varsigma)$, whenever $Q^s_{2r}(z_0)\Subset Q^s_2$. 
 
To do this, let us fix $\varsigma\in(0,1)$. We first consider 
\begin{align*}
    \dashint_{Q^s_{|h|^{\beta}}(x_i,t_j)}|\delta_h u|^p\,dz&\leq c  \dashint_{Q^s_{|h|^{\beta}}(x_i,t_j)}|\delta_h (u-v_{i,j})|^p\,dz+c \dashint_{Q^s_{|h|^{\beta}}(x_i,t_j)}|\delta_h v_{i,j}|^p\,dz\\
    &\eqqcolon J_1+J_2.
\end{align*}
By \eqref{maxd.ineq0}, we have 
\begin{align}\label{ineq.highdf}
    J_1\leq c\left(|h|^{-n\beta}|\mu|(Q^s_{4|h|^{\beta}}(x_i,t_j))\right)^p,
\end{align}
where $c=c(n,s_0,\Lambda,p)$. Using Theorem \ref{thm.hol} with $q=p$ and \eqref{rel.exc}, we have
\begin{align*}
    J_2^{\frac1p}&\leq |h| \norm{\nabla v_{i,j}}_{L^{\infty}(Q^s_{2|h|^\beta}(x_i,t_j))}\\
    &\leq c|h|^{1-\beta}E^p( v_{i,j};Q^s_{3|h|^\beta}(x_i,t_j))\\
    &\leq c|h|^{1-\beta}E^p_{\mathrm{loc}}(v_{i,j};Q^s_{3|h|^\beta}(x_i,t_j))\\
    &\quad+c|h|^{1-\beta}\left(\dashint_{I^s_{3|h|^\beta}(t_j)}\mathrm{Tail}(v_{i,j}-(v_{i,j})_{B_{3|h|^\beta}(x_i)}(t);B_{3\abs{h}^{\beta}}(x_i))^p\,dt\right)^{\frac1p}
\end{align*}
for some constant $c=c(n,s_0,\Lambda,p)$. We now use 
\eqref{bdd.est.zero} along with the fact that $v_{i,j}-(v_{i,j})_{Q^s_{4|h|^\beta}(x_i,t_j)}$ is also a weak solution to \eqref{eq.ivp} with $z_0=(x_i,t_j)$ and $R=4|h|^\beta$, \eqref{rel.exc} and \cite[Lemma 2.2]{DieKimLeeNow24}, in order to see that 
\begin{align*}
    J_2^{\frac1p}&\leq c|h|^{1-\beta}E_{\mathrm{loc}}(v_{i,j};Q^s_{4|h|^\beta}(x_i,t_j))\\
    &\quad+c|h|^{1-\beta}\left(\dashint_{Q^s_{4|h|^\beta}(x_i,t_j)}|v_{i,j}-(v_{i,j})_{B_{4|h|^\beta}(x_i)}(t)|^p\,dz\right)^{\frac1p}\\
    &\quad+c|h|^{1-\beta}\left(\dashint_{I^s_{4|h|^\beta}(t_j)}\mathrm{Tail}(v_{i,j}-(v_{i,j})_{B_{4|h|^\beta}(x_i)}(t);B_{4|h|^\beta}(x_i))^p\,dt\right)^{\frac1p}\\
    &\eqqcolon |h|^{1-\beta}(J_{2,1}+J_{2,2}+J_{2,3})
\end{align*}
for some constant $c=c(n,s_0,\Lambda,p)$.

Let us recall the constant $\sigma\geq 2s-1$ determined in \eqref{cond.sigma} to see that 
\begin{align}\label{ineq.ssigma}
    (1-s)/{(1-\sigma)}\leq c
\end{align}
for some constant $c=c(n,s_0,p)$.
As in the estimate of $J_2$ given in \eqref{est.j.glu} with $q=1$, we obtain 
\begin{align*}
    {J}_{2,1}&\leq c\dashint_{Q^s_{4|h|^\beta}(x_i,t_j)}|v_{i,j}-(v_{i,j})_{B_{4|h|^\beta}(x_i)}(t)|\,dz+c{J}_{2,3}\\
    &\quad+c(1-s)|h|^{(2s-1)\beta}\dashint_{Q^s_{4|h|^\beta}(x_i,t_j)}\int_{B_{4|h|^\beta}(x_i)}\frac{|v_{i,j}(x,t)-v_{i,j}(y,t)|}{|x-y|^{n+2s-1}}\,dy\,dz\\
    &\leq c\dashint_{Q^s_{4|h|^\beta}(x_i,t_j)}|v_{i,j}-(v_{i,j})_{B_{4|h|^\beta}(x_i)}(t)|\,dz+c{J}_{2,3}\\
    &\quad+c|h|^{\sigma \beta}\left((1-\sigma)\dashint_{Q^s_{4|h|^\beta}(x_i,t_j)}\int_{B_{4|h|^\beta}(x_i)}\frac{|v_{i,j}(x,t)-v_{i,j}(y,t)|^p}{|x-y|^{n+\sigma p}}\,dy\,dz\right)^{\frac1p}
\end{align*}
for some constant $c=c(n,s_0,\Lambda,p)$, where we have used H\"older's inequality for the last inequality, \eqref{cond.sigma} and \eqref{ineq.ssigma}. We next use H\"older's inequality and a non-scaled version of \eqref{ineq.zero.wp} to estimate $J_{2,1}$ as 
\begin{equation}\label{est.j21.sd}
\begin{aligned}
    {J}_{2,1}&\leq cJ_{2,2}+c{J}_{2,3}+c|h|^{-n\beta}|\mu|(Q^s_{4|h|^\beta}(x_i,t_j))\\
    &\quad+c|h|^{\sigma\beta}\left((1-\sigma)\dashint_{Q^s_{4|h|^\beta}(x_i,t_j)}\int_{B_{4|h|^\beta}(x_i)}\frac{|u(x,t)-u(y,t)|^p}{|x-y|^{n+\sigma p}}\,dy\,dz\right)^{\frac1p}\\
    &\leq cJ_{2,2}+c{J}_{2,3}+c|h|^{-n\beta}|\mu|(Q^s_{4|h|^\beta}(x_i,t_j))\\
    &\quad+c|h|^{\beta{\gamma}}\left((1-\sigma)\dashint_{Q^s_{4|h|^\beta}(x_i,t_j)}\int_{B_{4|h|^\beta}(x_i)}\frac{|u(x,t)-u(y,t)|^p}{|x-y|^{n+{\gamma}p}}\,dy\,dz\right)^{\frac1p}
\end{aligned}
\end{equation}
for some constant $c=c(n,s_0,\Lambda,p)$, where we have also used \cite[Equation (2.6)]{DieKimLeeNow24d} for the last inequality.
We now estimate $J_{2,2}$ as 
\begin{equation}\label{est.j22.sd}
\begin{aligned}
    J_{2,2}&\leq c\left(\dashint_{Q^s_{4|h|^\beta}(x_i,t_j)}|v_{i,j}-u|^p\,dz\right)^{\frac1p}
   +c\left(\dashint_{Q^s_{4|h|^\beta}(x_i,t_j)}|u-(u)_{B_{4|h|^\beta}(x_i)}(t)|^p\,dz\right)^{\frac1p}\\
    &\leq c|h|^{-n\beta}|\mu|(Q^s_{4|h|^\beta}(x_i,t_j))+c\left(\dashint_{Q^s_{4|h|^\beta}(x_i,t_j)}|u-(u)_{B_{4|h|^\beta}(x_i)}(t)|^p\,dz\right)^{\frac1p}\\
    &\leq c|h|^{-n\beta}|\mu|(Q^s_{4|h|^\beta}(x_i,t_j))
   +c\left((1-\gamma)\dashint_{I^s_{4|h|^\beta}(t_j)}|h|^{\beta(\gamma p-n)}[u]^p_{W^{\gamma,p}(B_{4|h|^\beta}(x_i))}\,dt\right)^{\frac1p}
\end{aligned}
\end{equation}
for some constant $c=c(n,s_0,\Lambda,p)$ and $\gamma\in(0,1)$, where we have used Lemma \ref{lem.comp} and \cite[Lemma 2.2]{DieKimLeeNow24d}.
On the other hand, we note that in view of \cite[Lemma 2.2]{DieKimLeeNow24} we have
\begin{equation}\label{hd.j12}
\begin{aligned}
    &\left(\dashint_{I^s_{4|h|^\beta}(t_j)}\mathrm{Tail}\left(u-(u)_{B_{4|h|^\beta}(x_i)}(t);B_{4|h|^\beta}(x_i)\right)^p\,dt\right)^{\frac1p}\\
    &\leq c\left(\dashint_{I^s_{4|h|^\beta}(t_j)}\left(\sum_{k=0}^{m_0+2}2^{-2sk}\dashint_{B_{2^{k+2}|h|^{\beta}}(x_i)}|u-(u)_{B_{2^{k+2}|h|^{\beta}}(x_i)}(t)|\,dx\right)^{p}\,dt\right)^{\frac1p}\\
    &\quad +c\left(\dashint_{I^s_{4|h|^\beta}(t_j)}2^{-2sm_0p}{\mathrm{Tail}}\left(u-(u)_{B_{2^{m_0+4}|h|^{\beta}}(x_i)}(t);B_{2^{m_0+4}|h|^{\beta}}(x_i)\right)^p\,dt\right)^{\frac1p}\\
    &\eqqcolon J_{2,3,1}+J_{2,3,2}
\end{aligned}
\end{equation}
for some constant $c=c(n,s_0)$.
Using H\"older's inequality and the fractional Poincar\'e inequality from \cite[Lemma 2.3]{DieKimLeeNow24d}, we further estimate $J_{2,3,1}$ as 
\begin{align*}
    J_{2,3,1}&\leq c\left(\dashint_{I^s_{4|h|^\beta}(t_j)}\sum_{k=0}^{m_0+2}2^{-2sk}\dashint_{B_{2^{k+2}|h|^{\beta}}(x_i)}|u-(u)_{B_{2^{k+2}|h|^{\beta}}(x_i)}(t)|^p\,dz\right)^{\frac1p}\\
    &\leq c\left((1-\gamma)\dashint_{I^s_{4|h|^\beta}(t_j)}\sum_{k=0}^{m_0+2}2^{-2sk}(2^{k}|h|^\beta)^{\gamma p-n}[u]^p_{W^{\gamma,p}(B_{2^{k+2}|h|^\beta}(x_i))}\,dt\right)^{\frac1p},
\end{align*}
where $c=c(n,s_0,p)$.  
Since $B_{2^{m_0+4}|h|^{\beta}}(x_i)\subset B_{3/4}$ and $x_i\in B_{5/8}$ by \eqref{el : choi.m0}, we have 
\begin{align*}
    J_{2,3,2}&\leq c\left(\dashint_{I^s_{4|h|^\beta}(t_j)}2^{-2sm_0p}\left(\dashint_{B_{3/4}}|u-(u)_{B_{3/4}}(t)|\,dx\right)^p\,dt\right)^{\frac1p}\\
    &\quad+c\left(\dashint_{I^s_{4|h|^\beta}(t_j)}2^{-2sm_0p}{\mathrm{Tail}}\left(u-(u)_{B_{3/4}}(t);B_{3/4}\right)^p\,dt\right)^{\frac1p}
\end{align*}
for some constant $c=c(n,s_0,\Lambda,p)$. In light of the above estimates of $J_{2,3,1}$ and $J_{2,3,2}$, we now estimate $J_{2,3}$ as 
\begin{equation}\label{est.j23.sd}
\begin{aligned}
    J_{2,3}&\leq c\left(\dashint_{Q^s_{4|h|^\beta}(x_i,t_j)}|v_{i,j}-u|^p\,dz\right)^{\frac1p}\\
&\quad+c\left(\dashint_{I^s_{4|h|^\beta}(t_j)}\mathrm{Tail}\left(u-(u)_{B_{4|h|^\beta}(x_i)}(t);B_{4|h|^\beta}(x_i)\right)^p\,dt\right)^{\frac1p}\\
&\leq c|h|^{-n\beta}|\mu|(Q^s_{4|h|^\beta}(x_i,t_j))\\
&\quad+c\left((1-\gamma)\dashint_{I^s_{4|h|^\beta}(t_j)}\sum_{k=0}^{m_0+2}2^{-2sk}(2^{k}|h|^\beta)^{\gamma p-n}[u]^p_{W^{\gamma,p}(B_{2^{k+2}|h|^\beta}(x_i))}\,dt\right)^{\frac1p}\\
&\quad+c2^{-2sm_0}J^{\frac1p}
\end{aligned}
\end{equation}
for some constant $c=c(n,s_0,\Lambda,p)$, where we have used Lemma \ref{lem.comp} and \eqref{hd.j12} along with the estimates of $J_{2,3,1}$ and $J_{2,3,2}$. Moreover, we have denoted
\begin{align}\label{defn.J}
    J\coloneqq\dashint_{I^s_{4|h|^\beta}(t_j)}\left[\dashint_{B_{3/4}}|u-(u)_{B_{3/4}}(t)|^p\,dx+{\mathrm{Tail}}\left(u-(u)_{B_{3/4}}(t);B_{3/4}\right)^p\right]\,dt.
\end{align}

Combining all the estimates \eqref{ineq.highdf}, \eqref{est.j21.sd}, \eqref{est.j22.sd} and \eqref{est.j23.sd}, we obtain
\begin{equation*}
\begin{aligned}
    J_1+J_2
    &\leq  c|h|^{-np\beta }(|\mu|(Q^s_{4|h|^\beta}(x_i,t_j)))^p\\
    &\quad+c|h|^{p(1-\beta)}(1-\sigma)\dashint_{I^s_{4|h|^\beta}(t_j)}\sum_{k=0}^{m_0+2}2^{-2sk}(2^{k}|h|^\beta)^{\gamma p-n}[u]^p_{W^{\gamma,p}(B_{2^{k+2}|h|^\beta}(x_i))}\,dt\\
&\quad+c|h|^{p(1-\beta)}2^{-2sm_0p}J\eqqcolon \sum_{k=1}^3 L_{k}(i,j)
\end{aligned}
\end{equation*}
for some constant $c=c(n,s_0,\Lambda,p)$, where we have also used the fact that $\sigma\leq\gamma$. Using this, we get
\begin{align*}
    \int_{Q^s_{1/2}}|\delta_h u|^p\,dz&\leq |h|^{(n+2s)\beta}\sum_{j\in \mathcal{J}}\sum_{i\in \mathcal{I}}\dashint_{Q^s_{|h|^{\beta}}(x_i,t_j)}|\delta_h u|^p\,dz\\
    &\leq |h|^{(n+2s)\beta}\sum_{j\in \mathcal{J}}\sum_{i\in \mathcal{I}} \sum_{k=1}^{3}L_k(i,j)
\end{align*}
for some constant $c=c(n,s_0,\Lambda,p)$. 
We first use \eqref{max.cov.ineq2} with $k=2$ and the inequality
\begin{align*}
    \sum_{j\in \mathcal{J}}\sum_{i\in \mathcal{I}} |\mu|(Q^s_{4|h|^\beta}(x_i,t_j))^p\leq \left(\sum_{j\in \mathcal{J}}\sum_{i\in \mathcal{I}} |\mu|(Q^s_{4|h|^\beta}(x_i,t_j))\right)^p,
\end{align*}
which holds since $p\geq1$, in order to see that 
\begin{equation}\label{maxd.estl1}
\begin{aligned}
    &|h|^{(n+2s)\beta}\sum_{j\in \mathcal{J}}\sum_{i\in \mathcal{I}} L_1(i,j)\\
    &\leq c|h|^{(2s+n(1-p))\beta}\left(\sum_{(i,j)\in \mathcal{I}\times\mathcal{J}}\int_{Q^s_1}\bfchi_{Q^s_{4|h|^\beta}(x_i,t_j)}(z)|\mu|\,dz\right)^p\\
    &\leq c|h|^{(2s+n(1-p))\beta}(|\mu|(Q^s_{3/4}))^p.
    \end{aligned}
\end{equation}
We observe from \eqref{max.cov.ineq2} and a few algebraic inequalities that
\begin{equation}\label{maxd.estl2}
\begin{aligned}
     &\frac{|h|^{(n+2s)\beta}}{1-\sigma}\sum_{j\in \mathcal{J}}\sum_{i\in \mathcal{I}} L_2(i,j)\\
     &\leq c|h|^{\gamma p\beta+p(1-\beta)}\sum_{j\in \mathcal{J}}\sum_{i\in \mathcal{I}}\sum_{k=0}^{m_0+2}2^{(-2s+\gamma p-n)k}\int_{I^s_{4|h|^\beta}(t_j)}[u]^p_{W^{\gamma,p}(B_{2^{k+2}|h|^\beta}(x_i))}\,dt\\
     &\leq  c|h|^{\gamma p\beta+p(1-\beta)}\sum_{k=0}^{m_0+2}\sum_{(i,j)\in \mathcal{I}\times \mathcal{J}}2^{(-2s+\gamma p-n)k}\int_{{Q}_k}\int_{B_{3/4}}|U_\gamma|^p\,dy\,dz\\
    &\leq c|h|^{\gamma p\beta+p(1-\beta)}\sum_{k=0}^{m_0+2}\sum_{(i,j)\in \mathcal{I}\times \mathcal{J}}2^{(-2s+\gamma p-n)k}\int_{Q^s_{3/4}}\bfchi_{{Q}_k}(z)\int_{B_{3/4}}|U_\gamma|^p\,dy\,dz\\
    &\leq c|h|^{\gamma p\beta+p(1-\beta)}\sum_{k=0}^{m_0+2}2^{(-2s+\gamma p)k}\int_{Q^s_{3/4}}\int_{B_{3/4}}|U_\gamma|^p\,dy\,dz\\
    &\leq c|h|^{\gamma p\beta+p(1-\beta)}[u]^p_{L^p(I^s_{3/4};W^{\gamma,p}(B_{3/4}))}
\end{aligned}
\end{equation}
for some constant $c=c(n,s_0,\Lambda,p,\gamma)$,
where we denote ${Q}_k\coloneqq B_{2^{k+2}|h|^\beta}(x_i)\times I^s_{4|h|^\beta}(t_j)$ and
\begin{align*}
    |U_\gamma|^p(z)=\int_{{B_{3/4}}}\frac{|u(x,t)-u(y,t)|^p}{|x-y|^{n+\gamma p}}\,dy\in L^1(Q^s_{{3/4}}),
\end{align*}
and we have used the fact that $\gamma p<2s$ by the assumptions that $p<\frac{n+2s_0}{n+1}$ and $s>1/2$.
We next estimate 
\begin{equation}\label{maxd.estl3}
\begin{aligned}
    &|h|^{(n+2s)\beta}\sum_{j\in \mathcal{J}}\sum_{i\in \mathcal{I}} L_3(i,j)\\
    &\leq c|h|^{n\beta+p(1-\beta)}|\mathcal{I}||h|^{2s\beta}\sum_{j\in\mathcal{J}}\int_{I_{4|h|^\beta}(t_j)}{\mathrm{Tail}}\left(u-(u)_{B_{3/4}}(t);B_{{3/4}}\right)^p\,dt\\
    &\quad +c|h|^{n\beta+p(1-\beta)}|\mathcal{I}||h|^{2s\beta}\sum_{j\in\mathcal{J}}\int_{I_{4|h|^\beta}(t_j)}\int_{B_2}|u-(u)_{B_{3/4}}(t)|^p\,dx\,dt\\
    &\leq c|h|^{p((1-\beta)+2s\beta)}{E}^p(u;Q^s_{3/4})^p
\end{aligned}
\end{equation}
for some constant $c=c(n,s_0,\Lambda,p)$, where we have used the fact that $|\mathcal{I}||h|^{n\beta}\leq c$, \eqref{max.cov.ineq1}, \eqref{el : choi.m0} and 
\begin{align*}
    \widetilde{E}^p(u-(u)_{B_{3/4}}(t);Q^s_{3/4})\leq c{E}^p(u;Q^s_{3/4}).
\end{align*}

Combining all the estimates \eqref{maxd.estl1}, \eqref{maxd.estl2} and \eqref{maxd.estl3} yields
\begin{equation}\label{ineq2.hd}
\begin{aligned}
    \int_{Q^s_{5/8}}|\delta_h u|^p\,dz&\leq c|h|^{(2s+n(1-p))\beta}(|\mu|(Q^s_{3/4}))^p\\
    &\quad+c(1-\sigma)|h|^{{\gamma} p\beta+p(1-\beta)}[u]^p_{L^p(I^s_{3/4};W^{{\gamma},p}(B_{3/4}))}\\
    &\quad+ c|h|^{p((1-\beta)+2s\beta)}{E}^p(u;Q^s_{3/4})^p
\end{aligned}
\end{equation}
for some constant $c=c(n,s_0,\Lambda,p,\gamma)$. We are ready to prove \eqref{maxd.goal1}. We first choose the constant $\beta$ such that
\begin{align*}
    \beta=\frac{p(1+\varsigma)}{2(2s_0+n(1-p))}>0
\end{align*}
to see that $(2s_0+n(1-p))\beta= p(1+\varsigma)/2$. Since $p<(n+2s_0)/(n+1)$, $\beta\in(0,1)$ holds.
Let us define a sequence
\begin{equation}\label{seq.gammak}
    \gamma_0\coloneqq\sigma\quad\text{and}\quad \gamma_k\coloneqq\gamma_0\beta^{k-1}+\sum_{i=0}^{k-1}(1-\beta-\epsilon)\beta^{i}\quad\text{if }k\geq1,
\end{equation}
where the constant $\sigma$ is determined in \eqref{cond.sigma} and
\begin{equation*}
    \epsilon={(1-\beta)(1-\varsigma)}/{4}
\end{equation*}
to see that $\lim\limits_{k\to\infty}\gamma_k=(3+\varsigma)/4$. Thus, there is a positive number $i_\varsigma$ such that 
\begin{align}\label{choi.iv}
    \gamma_{i_\varsigma}>(1+\varsigma)/2.
\end{align}
We first apply Lemma \ref{lem.firemb} into \eqref{ineq2.hd} with $\gamma=\gamma_0$ to see that 
\begin{equation}\label{ind.frst.sd}
\begin{aligned}
    [u]_{L^p(I^s_{1/2};W^{\gamma_1,p}(B_{1/2}))}&\leq c (|\mu|(Q^s_{3/4}))^p+c(1-\gamma_0)[u]^p_{L^p(I^s_{3/4};W^{{\gamma_0},p}(B_{3/4}))}\\
    &\quad+ c{E}^p(u;Q^s_{3/4})^p
\end{aligned}
\end{equation}
for some constant $c=c(n,s_0,\Lambda,p,\varsigma)$.
To estimate the second term in the right-hand side of \eqref{ind.frst.sd}, we consider the weak solution $v$ to \eqref{eq.ivp} with $Q^s_R(z_0)$ replaced by $Q^s_1$ to see that 
\begin{align*}
    &(1-\gamma_0)\int_{Q^s_{3/4}}\int_{B_{3/4}}\frac{|u(x,t)-u(y,t)|^p}{|x-y|^{n+\gamma_0 p}}\,dy\,dz\\
    &\leq c(|\mu|(Q^s_{1}))^p+c(1-\gamma_0)\int_{Q^s_{3/4}}\int_{B_{3/4}}\frac{|v(x,t)-v(y,t)|^p}{|x-y|^{n+\gamma_0 p}}\,dy\,dz
\end{align*}
for some constant $c=c(n,s_0,\Lambda)$, where we have used \eqref{seq.gammak}, \eqref{ineq.zero.wp} and \eqref{sta.cons}. Using H\"older's inequality and the standard energy inequality given in the proof of \cite[Theorem 1.8]{KaWe23} along with \eqref{bdd.est.zero}, we get 
\begin{align*}
    &(1-\gamma_0)\int_{Q^s_{3/4}}\int_{B_{3/4}}\frac{|v(x,t)-v(y,t)|^p}{|x-y|^{n+\gamma_0 p}}\,dy\,dz\\
    &\leq c(1-s)(s-\gamma_0)^{-\frac12}\left(\int_{Q^s_{3/4}}\int_{B_{3/4}}\frac{|v(x,t)-v(y,t)|^2}{|x-y|^{n+2s}}\,dy\,dz\right)^{\frac12}\\
    &\leq cE^2_{\mathrm{loc}}(v;Q^s_{13/16})+c\int_{I^s_{13/16}}\mathrm{Tail}(v-(v)_{Q^s_{13/16}};B_{13/16})\,dt
    \leq cE(v;Q^s_{7/8}),
\end{align*}
where we have used \eqref{sta.cons} to see that $(1-s)(s-\gamma_0)^{-\frac12}\leq c(1-s)^{\frac12}$ for some constant $c=c(n,s_0,p)$.
By Lemma \ref{lem.comp} and H\"older's inequality, we deduce 
\begin{align*}
    E(v;Q^s_{7/8})\leq cE(u;Q^s_{7/8})+c|\mu|(Q^s_1)\leq cE^p(u;Q^s_{15/16})+c|\mu|(Q^s_1),
\end{align*}
where $c=c(n,s_0,\Lambda)$.

Combining the above three inequalities with \eqref{ind.frst.sd}, we get
\begin{align*}
    [u]_{L^p(I^s_{1/2};W^{\gamma_1,p}(B_{1/2}))}&\leq c{E}^p(u;Q^s_{15/16})^p+ c (|\mu|(Q^s_{1}))^p
\end{align*}
for some constant $c=c(n,s_0,\Lambda,p)$. By standard covering arguments, we deduce
\begin{align}\label{ind.r}
    r^{\gamma_1}[u]_{L^p(I^s_{r}(t_0);W^{\gamma_1,p}(B_{r}(x_0)))}&\leq  c{E}^p(u;Q^s_{2r}(z_0))^p+c r^{-n}(|\mu|(Q^s_{2r}(z_0)))^p
\end{align}
for some constant $c=c(n,s_0,\Lambda,p,\varsigma)$, whenever $Q^s_{2r}(z_0)\Subset Q^s_2$. We now use Lemma \ref{lem.firemb}, \eqref{ineq2.hd} with $\gamma=\gamma_1$ and \eqref{ind.r} to see that 
\begin{align*}
    [u]_{L^p(I^s_{1/2};W^{\gamma_2,p}(B_{1/2}))}&\leq  c{E}^p(u;Q^s_{1})^p+c (|\mu|(Q^s_{1}))^p
\end{align*}
for some constant $c=c(n,s_0,\Lambda,p,\varsigma)$. By standard covering arguments, we obtain
\begin{align*}
   r^{\gamma_2} [u]_{L^p(I^s_{r}(t_0);W^{\gamma_2,p}(B_{r}(x_0)))}&\leq  c{E}^p(u;Q^s_{2r}(z_0))^p+cr^{-n} (|\mu|(Q^s_{2r}(z_0)))^p
\end{align*}
for some constant $c=c(n,s_0,\Lambda,p,\varsigma)$, whenever $Q^s_{2r}(z_0)\Subset Q^s_2$. By iterating the above procedure $i_{\varsigma}-2$ times, we obtain
\begin{align*}
   r^{\gamma_{i_\varsigma}} [u]_{L^p(I^s_{r}(t_0);W^{\gamma_{i_\varsigma},p}(B_{r}(x_0)))}&\leq  c{E}^p(u;Q^s_{2r}(z_0))^p+c r^{-n}(|\mu|(Q^s_{2r}(z_0)))^p
\end{align*}
for some constant $c=c(n,s_0,\Lambda,p,\varsigma)$, whenever $Q^s_{2r}(z_0)\Subset Q^s_2$. By \cite[Equation (2.6)]{DieKimLeeNow24d}, we arrive at
\begin{align*}
   r^{{\varsigma}} [u]_{L^p(I^s_{r}(t_0);W^{{\varsigma},p}(B_{r}(x_0)))}&\leq  c{E}^p(u;Q^s_{2r}(z_0))^p+c r^{-n}(|\mu|(Q^s_{2r}(z_0)))^p,
\end{align*}
which completes the proof of \eqref{maxd.goal1}.

Using \eqref{maxd.goal1}, we now prove our desired gradient estimate, which requires the use of second-order difference quotients. Similar to the proof of \eqref{maxd.goal1}, we observe that
\begin{align*}
    \dashint_{Q^s_{|h|^{\beta}}(x_i,t_j)}|\delta^2_h u|^p\,dz&\leq c\left(|h|^{-n\beta}|\mu|(Q^s_{4|h|^{\beta}}(x_i,t_j))\right)^p  +c \dashint_{Q^s_{|h|^{\beta}}(x_i,t_j)}|\delta^2_h v_{i,j}|^p\,dz
\end{align*}
for some constant $c=c(n,s_0,\Lambda,p)$. We now use Theorem \ref{thm.hol} to see that 
there is a constant $\alpha=\alpha(n,s_0,\Lambda,p)$ such that
\begin{align*}
   \dashint_{Q^s_{|h|^{\beta}}(x_i,t_j)}|\delta^2_h v_{i,j}|^p\,dz &\leq |h|^{p(1+\alpha)}\left[v_{i,j}\right]_{C^{1,\alpha}(Q^s_{2|h|^\beta}(x_i,t_j))}\\
   &\leq c|h|^{p(1+\alpha)(1-\beta)}E^p(v_{i,j};Q^s_{4|h|^\beta}(x_i,t_j))^p.
\end{align*}
Therefore, as in the estimate of \eqref{ineq2.hd}, we have 
\begin{equation*}
\begin{aligned}
    \int_{Q^s_{5/8}}|\delta^2_h u|^p\,dz&\leq c|h|^{(2s+n(1-p))\beta}(|\mu|(Q^s_{1}))^p\\
    &\quad+c(1-\sigma)|h|^{{\gamma} p\beta+p(1+\alpha)(1-\beta)}[u]^p_{L^p(I^s_{3/4};W^{{\gamma},p}(B_{3/4}))}\\
    &\quad+ c|h|^{p((1+\alpha)(1-\beta)+2s\beta)}{E}^p(u;Q^s_{3/4})^p
\end{aligned}
\end{equation*}
for some constant $c=c(n,s_0,\Lambda,p)$. We now choose the constant $\beta$ as
\begin{equation}\label{choi2.beta}
    \beta=\frac{p}{2(2s_0+n(1-p))}+\frac12
\end{equation}
to see that $(2s_0+n(1-p))\beta>p$ and $\beta\in(1/2,1)$, as $p<(n+2s_0)/(n+1)$. We next select the constant $\gamma$ such that 
\begin{align*}
    \gamma=\frac{1-(1+\alpha)(1-\beta)}{2\beta}+\frac12<1
\end{align*}
to see that $\gamma \beta+(1+\alpha)(1-\beta)>1$. We now choose
\begin{align*}
    \widetilde{\gamma}=\min\{(2s+n(1-p))\beta/p,\gamma\beta+(1+\alpha)(1-\beta)\}-1>0
\end{align*}
to obtain
\begin{equation}\label{emb.forsd}
\begin{aligned}
    \int_{Q^s_{5/8}}|\delta^2_h u|^p\,dz&\leq c|h|^{p(1+\widetilde{\gamma})}(|\mu|(Q^s_{1}))^p\\
    &\quad+c|h|^{p(1+\widetilde{\gamma})}[u]^p_{L^p(I^s_{3/4};W^{{\gamma},p}(B_{3/4}))}+ c|h|^{p(1+\widetilde{\gamma})}{E}^p(u;Q^s_{3/4})^p
\end{aligned}
\end{equation}
for some constant $c=c(n,s_0,\Lambda,p)$. Using Lemma \ref{lem.secemb}, \eqref{emb.forsd}, \eqref{ineq2.hd} and \eqref{maxd.goal1} with $\gamma=\widetilde{\gamma}$, we arrive at
\begin{align*}
    \|\nabla u\|_{L^p(Q^s_{1/2})}\leq cE^p(u;Q^s_1)+c|\mu|(Q^s_1)
\end{align*}
for some constant $c=c(n,s_0,\Lambda,p)$, completing the proof.
\end{proof}

To obtain estimates of the gradient in fractional Sobolev spaces that are suitable to prove gradient potential estimates, we need the following estimate for solutions to homogeneous parabolic nonlinear nonlocal equations.
\begin{lemma}\label{lem.del.cac}
    Let us fix $h\in B_{s^{\frac1{2s}}/1000}\setminus\{0\}$ and $\beta\in(0,1)$. Let 
    \begin{equation*}
        v\in L^{2}(I^s_{4|h|^\beta}(t_{0});W^{s,2}(B_{4|h|^{\beta}}(x_{0})))\cap L^{p}(I^s_{4|h|^\beta}(t_{0});L^{1}_{2s}(\bbR^{n}))
    \end{equation*}
    be a weak solution to 
    \begin{align*}
        \partial_t v+\mathcal{L}v=0\quad\text{in }Q^s_{4|h|^\beta}(z_0).
    \end{align*}
    Then we have 
    \begin{align*}
        \dashint_{Q^s_{|h|^\beta}(z_0)}|\delta_h^2 v|^p\,dx\leq c|h|^{sp(1-\beta)}E(\delta_h v;Q^s_{4|h|^\beta}(z_0))^p
    \end{align*}
    for some constant $c=c(n,s_0,\Lambda,p)$.
\end{lemma}
\begin{proof}
In view of a scaling argument, we may assume $z_0=0$. We first observe from H\"older's inequality and \cite[Proposition 2.6]{BraLin17} that 
\begin{align*}
    \dashint_{Q^s_{|h|^\beta}}|\delta_h^2 v|^p\,dx&\leq |h|^{sp}\sup_{0<\overline{h}<|h|}\dashint_{Q^s_{|h|^\beta}}\left|\frac{\delta_{\overline{h}}}{|\overline{h}|^s}\left({\delta_h v}-(\delta_h v)_{Q^s_{2|h|^\beta}}\right)\right|^p\,dz\\
    &\leq |h|^{sp}\sup_{0<\overline{h}<|h|^\beta}\dashint_{Q^s_{|h|^\beta}}\left|\frac{\delta_{\overline{h}}}{|\overline{h}|^s}\left({\delta_h v}-(\delta_h v)_{Q^s_{2|h|^\beta}}\right)\right|^p\,dz\\
    &\leq |h|^{sp}\left(\dashint_{I^s_{|h|^\beta}}\sup_{0<\overline{h}<|h|^\beta}\dashint_{B_{|h|^\beta}}\left|\frac{\delta_{\overline{h}}}{|\overline{h}|^s}\left({\delta_h v}-(\delta_hv)_{Q^s_{2|h|^\beta}}\right)\right|^2\,dx\,dt\right)^{\frac{p}2}\\
    &\leq c|h|^{sp}\left((1-s)|h|^{-(n+2s)\beta}[\delta_h v]^2_{L^2(I^s_{2|h|^\beta};W^{s,2}(B_{2|h|^\beta}))}\right)^{\frac{p}2}\\
    &\quad+c|h|^{sp}\left(|h|^{-2s\beta}E^{2}_{\mathrm{loc}}(\delta_h v;Q^s_{2|h|^\beta})^2\right)^{\frac{p}2}\eqqcolon J
\end{align*}
for some constant $c=c(n,s_0,p)$. We note from Lemma \ref{lem.linsol} that $\delta_hv$ is a weak solution to 
\begin{align*}
    \partial_t \delta_hv +\mathcal{L}_{A_h}\delta_h v=0\quad\text{in }Q^s_{3|h|^\beta}
\end{align*}
for some coefficient $A_h$ satisfying \eqref{lin.kern.cond} and \eqref{lin.kern.cond2}.
Thus, by the standard energy inequalities given in the proof of \cite[Theorem 1.8]{KaWe23}, we get 
\begin{align*}
    &(1-s)|h|^{-(n+2s)\beta}[\delta_h v]^2_{L^2(I^s_{2|h|^\beta};W^{s,2}(B_{2|h|^\beta}))}\\
    &\leq c|h|^{-2s\beta}E^{2}_{\mathrm{loc}}(\delta_h v;Q^s_{3|h|^\beta})^2
    +c|h|^{-2s\beta}\left(\dashint_{I^s_{3|h|^\beta}}\mathrm{Tail}(\delta_h v-(\delta_h v)_{Q^s_{3|h|^\beta}};B_{3|h|^\beta})\,dt\right)^2.
\end{align*}
Plugging this inequality into the term $J$ along with Lemma \ref{lem.exc.hom} yields
\begin{align*}
    \dashint_{Q^s_{|h|^\beta}(z_0)}|\delta_h^2 v|^p\,dx&\leq c|h|^{sp}\left(|h|^{-2s\beta}E^{2}_{\mathrm{loc}}(\delta_h v;Q^s_{3|h|^\beta})^2\right)^{\frac{p}{2}}\\
    &\quad+c|h|^{sp(1-\beta)}\left(\dashint_{I^s_{3|h|^\beta}}\mathrm{Tail}(\delta_hv-(\delta_hv)_{Q^s_{3|h|^\beta}};B_{3|h|^\beta})\,dt\right)^p\\
    &\leq c|h|^{sp(1-\beta)}E(\delta_h v;Q^s_{4|h|^\beta})^p
\end{align*}
for some constant $c=c(n,s_0,\Lambda,p)$, which completes the proof.
\end{proof}

Next, we prove higher differentiability of the gradient with an explicit estimate.
\begin{lemma}\label{lem.highsu}
    Fix $p\in\left(1,\frac{n+2s_0}{n+1}\right)$, let $\mu\in L^1(I^s_{2R}(t_0);L^\infty(B_{2R}(x_0))$ and assume that
    $$
    u \in L^2(I^s_{2R}(t_0);W^{s,2}(B_{2R}(x_0)))\cap C(I^s_{2R}(t_0);L^2(B_{2R}(x_0))) \cap L^{p}(I^s_{2R}(t_0);L^{1}_{2s}(\bbR^n))
    $$
    is a weak solution to \eqref{eq:nonlocaleq} in $Q^s_{2R}(z_0)$ with $\nabla u\in L^p(\bbR^n\times I^s_{2R}(t_0))$ and $\mu\in L^1(I^s_{2R}(t_0);L^\infty(B_{2R}(x_0)))$. Then there is a constant $\sigma_0=\sigma_0(n,s_0,\Lambda,p) \in (0,1)$ such that $\nabla u\in L^{p}(I^s_{R}(t_0);W^{\sigma_0,p}(B_R(x_0)))$ with the estimate 
    \begin{align*}
       & R^{\sigma_0}\left(\dashint_{Q^s_{R}(z_0)}\int_{B_{R}(x_0)}\frac{|\nabla u(x,t)-\nabla u(y,t)|^p}{|x-y|^{n+\sigma_0 p }}\,dy\,dz\right)^{\frac1p}\\
        &\leq cE^p(\nabla u;Q^s_{2R}(z_0))+cR^{-(n+1)}|\mu|(Q^s_{2R}(z_0))
    \end{align*}
    for some constant $c=c(n,s_0,\Lambda,p)$.
\end{lemma}
\begin{remark} \normalfont
    The assumption that $\nabla u \in L^p(\bbR^n\times I^s_{2R}(t_0))$ is not restrictive for our purposes, since in view of Lemma \ref{lem.graex} and the localization argument given in Lemma \ref{lem.loc} this assumption can always be removed in the end.
\end{remark}
\begin{proof}
We may assume $R=1/2$ and $z_0=0$. First fix $\beta$ determined in \eqref{choi2.beta}.
We now choose $h\in \bbR$ such that
\begin{equation*}
   0<|h|\leq s^{\frac1{2s}}/{1000}.
\end{equation*}
Then by Lemma \ref{lem.cov}, there is a covering 
$\{Q^s_{|h|^\beta}(x_i,t_j)\}_{(i,j)\in \mathcal{I}\times \mathcal{J}}$ of $Q^s_{5/8}$ such that $(x_i,t_j)\in Q^s_{5/8}$, \eqref{max.cov.ineq00}, \eqref{max.cov.ineq0} and \eqref{max.cov.ineq1}. In addition, there is a positive integer $m_0$ such that \eqref{el : choi.m0} holds.
Let $v_{i,j}$ be a weak solution given in \eqref{eq.ivp} with $z_0=(x_i,t_j)$ and $R=4|h|^\beta$.
We first observe from Lemma \ref{lem.del.cac} that
\begin{align*}
   \dashint_{ Q^s_{|h|^\beta}(x_i,t_j)}|{\delta}^2_hv_{i,j}|^p\,dz
   &\leq c|h|^{sp(1-\beta)}E(\delta_h v_{i,j};Q^s_{4|h|^\beta}(x_i,t_j))^p
\end{align*}
for some constant $c=c(n,s_0,\Lambda,p)$.
Using this and Lemma \ref{lem.comp}, we next observe
\begin{equation}\label{first.ss}
\begin{aligned}
    \dashint_{Q^s_{|h|^\beta}(x_i,t_j)}|\delta^2_h u|^p\,dz&\leq  c|h|^{-n\beta p}(|\mu|(Q^s_{4|h|^\beta}(x_i,t_j)))^p\\
    &\quad+c|h|^{sp(1-\beta)}E(\delta_h v_{i,j};Q^s_{4|h|^\beta}(x_i,t_j))^p\\
    &\leq c|h|^{-n\beta p}(|\mu|(Q^s_{4|h|^\beta}(x_i,t_j)))^p\\
    &\quad+c|h|^{sp(1-\beta)}E(\delta_h u;Q^s_{4|h|^\beta}(x_i,t_j))^p,
\end{aligned}
\end{equation}
where $c=c(n,s_0,\Lambda,p)$. We note from H\"older's inequality and \eqref{rel.exc} that
\begin{align*}
    &{E}(\delta_h u;Q^s_{4|h|^\beta}(x_i,t_j))^p\\
    &\leq E^p_{\mathrm{loc}}(\delta_h u;Q^s_{4|h|^\beta}(x_i,t_j))^p\\
    &\quad+c\left(\dashint_{I^s_{4|h|^\beta}(t_j)}\mathrm{Tail}\left(\delta_hu-(\delta_hu)_{B_{4|h|^\beta}(x_i)}(t);B_{4|h|^\beta}(x_i)\right)^p\,dt\right)^{\frac1p}\eqqcolon J_1+J_2
\end{align*}
for some constant $c=c(n,s_0)$.  As in \eqref{hd.j12} with $u$ replaced by $\delta_h u$, we deduce that
\begin{equation*}
\begin{aligned}
    &J_2\leq c\left(\dashint_{I^s_{4|h|^\beta}(t_j)}\left(\sum_{k=0}^{m_0+2}2^{-2sk}\dashint_{B_{2^{k+2}|h|^{\beta}}(x_i)}|\delta_hu-(\delta_hu)_{B_{2^{k+2}|h|^{\beta}}(x_i)}(t)|\,dx\right)^{p}\,dt\right)^{\frac1p}\\
    &\quad +c2^{-2sm_0}\left(\dashint_{I^s_{4|h|^\beta}(t_j)}{\mathrm{Tail}}\left(\delta_hu-(\delta_hu)_{B_{2^{m_0+4}|h|^{\beta}}(x_i)}(t);B_{2^{m_0+4}|h|^{\beta}}(x_i)\right)^p\,dt\right)^{\frac1p}\\
    &\eqqcolon J_{2,1}+J_{2,2},
\end{aligned}
\end{equation*}
where $c=c(n,s_0)$. We now use H\"older's inequality and a few algebraic inequalities to see that 
\begin{align*}
    J_{2,1}&\leq c\left(\sum_{k=0}^{m_0+2}2^{-2sk}\dashint_{I^s_{4|h|^\beta}(t_j)}\dashint_{B_{2^{k+2}|h|^{\beta}}(x_i)}|\delta_hu-(\delta_hu)_{Q^s_{3/4}}|^p\,dz\right)^{\frac1p}
\end{align*}
and
\begin{align*}
    J_{2,2}&\leq c2^{-2sm_0}\left(\dashint_{I^s_{4|h|^\beta}(t_j)}{\mathrm{Tail}}\left(\delta_h u-(\delta_h u)_{Q^s_{3/4}};B_{3/4}\right)^p\,dt\right)^{\frac1p}\\
    &\quad+c2^{-2sm_0}\left(\dashint_{I^s_{4|h|^\beta}(t_j)}\dashint_{B_{3/4}}|\delta_hu-(\delta_hu)_{Q^s_{3/4}}|^p\,dz\right)^{\frac1p}
\end{align*}
for some constant $c=c(n,s_0)$.

Combining all the estimates $J_{2,1}$ and $J_{2,2}$, we obtain 
\begin{align*}
    E(\delta_h u;Q^s_{4|h|^\beta}(x_i,t_j))^p&\leq \sum_{k=0}^{m_0+2}2^{-2sk}\dashint_{I^s_{4|h|^\beta}(t_j)}\dashint_{B_{2^{k+2}|h|^{\beta}}(x_i)}|\delta_hu-(\delta_hu)_{Q^s_{3/4}}|^p\,dz\\
    &\quad+ c2^{-2sm_0p}\dashint_{I^s_{4|h|^\beta}(t_j)}\dashint_{B_{3/4}}|\delta_hu-(\delta_hu)_{Q^s_{3/4}}|^p\,dz\\
    &\quad+c2^{-2sm_0p}\dashint_{I^s_{4|h|^\beta}(t_j)}{\mathrm{Tail}}\left(\delta_h u-(\delta_h u)_{Q^s_{3/4}};B_{3/4}\right)^p\,dt.
\end{align*}
Using this, \eqref{first.ss}, \eqref{max.cov.ineq0}, \eqref{max.cov.ineq1} and \eqref{el : choi.m0} along with \eqref{maxd.estl1} and \eqref{maxd.estl3}, we deduce
\begin{align*}
    \dashint_{Q^s_{5/8}}|\delta_h^2 u|^p\,dz&\leq c|h|^{(n+2s)\beta}\sum_{i\in\mathcal{I}}\sum_{j\in\mathcal{J}}\dashint_{Q^s_{|h|^\beta}(x_i,t_j)}|\delta_h^2 u|^p\,dz\\
    &\leq c|h|^{sp(1-\beta)}E^p(\delta_hu;Q^s_{3/4})^p+c|h|^{2s\beta+(1-p)n\beta}(|\mu|(Q^s_1))^p
\end{align*}
for some constant $c=c(n,s_0,\Lambda,p)$. As in the proof of \cite[Lemma 2.6]{DieKimLeeNow24}, we further estimate 
\begin{align*}
    \dashint_{Q^s_{5/8}}|\delta_h^2 u|^p\,dz\leq c|h|^{p+sp(1-\beta)}E^p(\nabla u;Q^s_{7/8})^p+c|h|^{2s\beta+(1-p)n\beta}(|\mu|(Q^s_1))^p.
\end{align*}
Using Lemma \ref{lem.thiemb} along with the choice of the constant $\beta$ given in \eqref{choi2.beta}, we get 
\begin{align*}
    [\nabla u]_{L^p(I^s_{1/2};W^{\sigma_0,p}(B_{1/2}))}\leq cE^p(\nabla u;Q^s_{1})^p+c|\mu|(Q^s_1)
\end{align*}
for some constant $c=c(n,s_0,\Lambda,p)$ by taking 
\begin{align*}
    2\sigma_0=\min\{(2s\beta+(1-p)n\beta)/p,s(1-\beta)+1\}-1.
\end{align*}
This completes the proof.
\end{proof}

\subsection{Gradient oscillation decay for homogeneous initial boundary value problems}
The aim of this subsection is to establish decay estimates for nonlinear nonlocal homogeneous initial boundary value problems that are suitable for obtaining gradient potential estimates in the presence of general measure data. We accomplish this by first proving suitable higher differentiability estimates for the spatial gradient with respect to the spatial and temporal variables separately and then interpolating between these estimates.
We begin with a such higher differentiability estimate with respect to the spatial direction.
\begin{lemma}\label{lem.sobgrasp}
    Fix $p\in\left(1,\frac{n+2s_0}{n+1}\right)$, let $\mu\in L^1(I^s_{2R}(t_0);L^\infty(B_{2R}(x_0))$ and assume that
    $$
    	u \in L^2(I^s_{2R}(t_0);W^{s,2}(B_{2R}(x_0)))\cap C(I^s_{2R}(t_0);L^2(B_{2R}(x_0))) \cap L^{p}(I^s_{2R}(t_0);L^{1}_{2s}(\bbR^n))
    $$
    is a weak solution to \eqref{eq:nonlocaleq} in $Q^s_{2R}(z_0)$ with $\nabla u\in L^p(\bbR^n\times I^s_{2R}(t_0))$. Let 
    $$
        v\in L^2(I^s_{R}(t_0);W^{s,2}(B_{R}(x_0)))\cap C(I^s_R(t_0);L^2(B_R(x_0))) \cap L^{p}(I^s_{R}(t_0);L^{1}_{2s}(\bbR^n))
    $$
    be the unique weak solution to 
    \begin{equation}\label{eq.ivp.dec}
\left\{
\begin{alignedat}{3}
\partial_t v+\mathcal{L}v&= 0&&\qquad \mbox{in  $Q^s_{R}(z_0)$}, \\
v&=u&&\qquad  \mbox{in $(\mathbb{R}^{n}\setminus B_{R}(x_0))\times I^s_R(t_0)$}, \\
v(\cdot,t_0-R^{2s})&=u(\cdot,t_0-R^{2s}) &&\qquad \mbox{in  $B_R(x_0)$}.
\end{alignedat} \right.
\end{equation}
Then there is a constant $\kappa=\kappa(n,s_0,\Lambda,p)\in(0,1)$ such that 
\begin{equation}\label{gra.mainest}
\begin{aligned}
    &R^{\frac{-(n+2s)}{q}+\kappa}[\nabla v]_{L^q(I^s_{r}(t_0);W^{\kappa,q}(B_r(x_0)))}\\
    &\leq \frac{cR^{n+4s}}{(\rho-r)^{n+4s}}\left[E^q_{\mathrm{loc}}(\nabla v;Q^s_{\rho}(z_0))+ cE^p(\nabla u;Q^s_{\rho}(z_0))\right]\\
    &\quad +\frac{cR^{n+4s}}{(\rho-r)^{n+4s}}\frac{|\mu|(Q^s_{R}(z_0))}{R^{n+1}}
\end{aligned}
\end{equation}
holds for any $q\in [p,\infty)$ and all $R/2< r< \rho\leq3R/4$, where $c=c(n,s_0,\Lambda,p,q)$.
\end{lemma}
\begin{proof}
Let us assume $z_0=0$ and $R=1$. By the assumption that $s_0>1/2$, we observe 
\begin{equation}\label{gra.cond.beta}
    2s_0\beta-1>0,
\end{equation}
where $\beta\coloneqq\frac{2s_0+1}{4s_0}<1$.
We now choose $h\in \bbR^n\setminus \{0\}$ such that 
\begin{equation}\label{h.exc}
    |h|\leq s^{\frac1{2s}}(\rho-r)/1000
\end{equation}
to see that
\begin{equation*}
    Q^s_{2^5|h|^\beta}(z)\subset Q^s_{(3r+\rho)/{4}}\quad\text{if }z\in Q^s_r.
\end{equation*}
Then by Lemma \ref{lem.cov}, there is a covering 
$\{Q^s_{|h|^\beta}(x_i,t_j)\}_{(i,j)\in \mathcal{I}\times \mathcal{J}}$ of $Q^s_{(7r+\rho)/8}$ such that $(x_i,t_j)\in Q^s_{(7r+\rho)/8}$, \eqref{max.cov.ineq00}, \eqref{max.cov.ineq0} and \eqref{max.cov.ineq1} hold. In addition, there is a constant $m_0\geq1$ such that 
\begin{align}\label{gra.cond.m0}
    s^{\frac1{2s}}(\rho-r)/8\leq2^{m_0+4}|h|^\beta<s^{\frac1{2s}}(\rho-r)/4.
\end{align}
We first note that by Lemma \ref{lem.exc.hom}, we have
\begin{equation}\label{gra.start}
\begin{aligned}
   \dashint_{Q^s_{|h|^{\beta}}(x_i,t_j)}|\delta^2_hv|^q\,dz &\leq |h|^{q\alpha}\left[\delta_h v\right]_{C^{\alpha}(Q^s_{|h|^\beta}(x_i,t_j))}^q\\
   &\leq c|h|^{q\alpha(1-\beta)}E^p(\delta_hv;Q^s_{4|h|^\beta}(x_i,t_j))^q
\end{aligned}
\end{equation}
for some constant $c=c(n,s_0,\Lambda,p,q)$. As in \eqref{hd.j12} with $u$ replaced by $\delta_h v$ along with H\"older's inequality, we deduce 
\begin{equation*}
\begin{aligned}
    &E^p(\delta_hv;Q^s_{4|h|^\beta}(x_i,t_j))\\
    &\leq c\left(\sum_{k=0}^{m_0+2}2^{-2sk}\dashint_{I^s_{4|h|^\beta}(t_j)}\dashint_{B_{2^{k+2}|h|^{\beta}}(x_i)}|\delta_hv-(\delta_hv)_{Q^s_{\frac{3r+\rho}{4}}}|^q\,dz\right)^{\frac1q}\\
    &\quad +c2^{-2sm_0}\left(\dashint_{I^s_{4|h|^\beta}(t_j)}{\mathrm{Tail}}\left(\delta_hv-(\delta_hv)_{B_{2^{m_0+4}|h|^{\beta}}(x_i)}(t);B_{2^{m_0+4}|h|^{\beta}}(x_i)\right)^p\,dt\right)^{\frac1p}\\
    &\eqqcolon J_{2,1}+J_{2,2}
\end{aligned}
\end{equation*}
for some constant $c=c(n,s_0,\Lambda,p)$, where we have also used that
\begin{align*}
    &\left(\dashint_{B_{2^{k+2}|h|^{\beta}}(x_i)}|\delta_hv-(\delta_hv)_{B_{2^{k+2}|h|^{\beta}}(x_i)}(t)|^q\,dx\right)^{\frac1q}\\
    &\leq c\left(\dashint_{B_{2^{k+2}|h|^{\beta}}(x_i)}|\delta_hu-(\delta_hu)_{Q^s_{\frac{3r+\rho}{4}}}|^q\,dx\right)^{\frac1q}.
\end{align*}
As in the estimate of the term $I_2$ given in \cite[Equation (4.18)]{DieKimLeeNow24}, we estimate $J_{2,2}$ as 
\begin{align*}
    J_{2,2}&\leq \frac{c2^{-2sm_0}}{(\rho-r)^{n+2s}}\left(\dashint_{I^s_{4|h|^\beta}(t_j)}{\mathrm{Tail}}\left(\delta_hv-(\delta_hv)_{B_{\frac{3r+\rho}{4}}}(t);B_{\frac{3r+\rho}{4}}\right)^p\,dt\right)^{\frac1p}\\
    &\quad+\frac{c2^{-2sm_0}}{(\rho-r)^{n}}\left(\dashint_{I^s_{4|h|^\beta}(t_j)}\dashint_{B_{\frac{3r+\rho}4}}|\delta_h v-(\delta_h v)_{B_{\frac{3r+\rho}{4}}}(t)|^q\,dz\right)^{\frac1q}\eqqcolon J_{2,2,1}+J_{2,2,2}
\end{align*}
for some constant $c=c(n,s_0,\Lambda,p)$, where we have used H\"older's inequality. After a few calculations along with \eqref{gra.cond.m0} and Lemma \ref{lem.comp}, we estimate $J_{2,2,1}$ as 
\begin{align*}
    J_{2,2,1}&\leq \frac{c2^{-2sm_0}}{(\rho-r)^{n+2s}}\left(\dashint_{I^s_{4|h|^\beta}(t_j)}{\mathrm{Tail}}\left(\delta_hu-(\delta_hu)_{B_{\frac{3r+\rho}{4}}}(t);B_{\frac{3r+\rho}{4}}\right)^p\,dt\right)^{\frac1p}\\
    &\quad+ \frac{c|h|^{2s\beta}}{(\rho-r)^{n+4s}}\sup_{t\in I^s_{4|h|^\beta}(t_j)}\int_{B_1}|u-v|\,dx\\
    &\leq \frac{c2^{-2sm_0}}{(\rho-r)^{n+2s}}\left(\dashint_{I^s_{4|h|^\beta}(t_j)}{\mathrm{Tail}}\left(\delta_hu-(\delta_hu)_{B_{\frac{3r+\rho}{4}}}(t);B_{\frac{3r+\rho}{4}}\right)^p\,dt\right)^{\frac1p}\\
    &\quad+ \frac{c|h|^{2s\beta}}{(\rho-r)^{n+4s}}|\mu|(Q^s_1)
\end{align*}
for some constant $c=c(n,s_0,\Lambda,p)$. Combining all the estimates \eqref{gra.start}, $J_{2,1}$ and $J_{2,2}$, we get 
\begin{align*}
    &\int_{Q^s_{(7r+\rho)/8}}|\delta_h^2v|^q\,dz\\
    &\leq c|h|^{(n+2s)\beta}\sum_{i\in\mathcal{I}}\sum_{j\in\mathcal{J}}\dashint_{Q^s_{|h|^{\beta}}(x_i,t_j)}|\delta^2_hv|^q\,dz\\
    &\leq \frac{c|h|^{q\alpha(1-\beta)}}{(\rho-r)^{q(n+4s)}}\int_{Q^s_{\frac{3r+\rho}{4}}}\left|\delta_hv-(\delta_hv)_{Q^s_{\frac{3r+\rho}4}}\right|^q\,dz\\
    &\quad+ c|h|^{q\alpha(1-\beta)}\frac{c|h|^{2s\beta q}}{(\rho-r)^{q(n+4s)}}\left(\dashint_{I^s_{\rho}}{\mathrm{Tail}}\left(\delta_hu-(\delta_hu)_{B_{\frac{3r+\rho}{4}}}(t);B_{\frac{3r+\rho}{4}}\right)^p\,dt\right)^{\frac{q}p}\\
    &\quad+ \frac{c|h|^{2s\beta q+q\alpha(1-\beta)}}{(\rho-r)^{q(n+4s)}}(|\mu|(Q^s_1))^q
\end{align*}
for some constant $c=c(n,s_0,\Lambda,p,q)$, where we have also used \eqref{max.cov.ineq0}, \eqref{max.cov.ineq1} and \eqref{max.cov.ineq2}. As in the proof of \cite[Lemma 2.6]{DieKimLeeNow24} along with a few simple calculations, we deduce 
\begin{align*}
    \int_{Q^s_{(7r+\rho)/8}}|\delta_h^2v|^q\,dz
    &\leq \frac{c|h|^{q\alpha(1-\beta)+q}}{(\rho-r)^{n+4s}}\int_{Q^s_{\frac{r+\rho}{2}}}|\nabla v-(\nabla v)_{Q^s_{\frac{r+\rho}2}}|^q\,dz\\
    &\quad+ \frac{c|h|^{q\alpha(1-\beta)+q}}{(\rho-r)^{n+4s}}E^p(\nabla u;Q^s_{\frac{r+\rho}{2}})^q+\frac{c|h|^{2s\beta q+q\alpha(1-\beta)}}{(\rho-r)^{n+4s}}(|\mu|(Q^s_1))^q
\end{align*}
for some constant $c=c(n,s_0,\Lambda,p,q)$. With the choice of $\beta$ given in \eqref{gra.cond.beta}, we now apply Lemma \ref{lem.thiemb} into the above inequality to see that 
\begin{align*}
    [\nabla v]_{L^q(I^s_{r};W^{\kappa,q}(B_r))}&\leq \frac{c}{(\rho-r)^{n+4s}}\left[E^q_{\mathrm{loc}}(\nabla v;Q^s_{\rho})+ cE^p(\nabla u;Q^s_{\rho})+|\mu|(Q^s_{1})\right]
\end{align*}
for some constant $c=c(n,s_0,\Lambda,p,q)$, where the constant
\begin{align*}
    \kappa\coloneqq\alpha(1-\beta)/2
\end{align*}
depends only on $n,s_0,\Lambda$ and $p$. This completes the proof.
\end{proof}



Next, we establish  an estimate of the H\"older seminorm of solutions to homogeneous initial boundary value problems that involves affine functions, making it suitable for executing iteration prodedures in our first-order setting. 
\begin{lemma}\label{lem.higsotimesol}
    Let us fix $p\in\left(1,\frac{n+2s_0}{n+1}\right)$. Let $v$ be a weak solution to \eqref{eq.ivp.dec}. For any affine function $l=A\cdot x+b$ with $A\in\bbR^{n}$ and $b\in\bbR$, we have 
    \begin{align*}
         &R^{\sigma_1}\left[{v-l}\right]_{C^{\sigma_1}(Q^s_{R/2}(z_0))}\\
         &\leq cE^p_{\mathrm{loc}}(v-l;Q^s_{3R/4}(z_0))\\
        &\quad+c\left(\dashint_{I^s_{3R/4}(t_0)}\mathrm{Tail}(u-l-(u-l)_{B_{3R/4}(x_0)}(t);B_{3R/4}(x_0))^p\,dt\right)^{\frac{1}p}
        +cR^{-n}|\mu|(Q^s_{R}(z_0))
    \end{align*}
    for some $\sigma_1=\sigma_1(n,s_0,\Lambda,p)\in(0,1)$ and some $c=c(n,s_0,\Lambda,p)$.
\end{lemma}
\begin{proof}
We may assume $R=1$ and $z_0=0$. We are first going to prove that for any $q\in[p,\infty)$,
\begin{equation}\label{osc.vminusl}
\begin{aligned}
        &[v-l]_{W^{\frac{\sigma_1}{s},q}(I^s_{9/16};L^q(B_{9/16}))}+[v-l]_{L^{q}(I^s_{9/16};W^{2\sigma_1,q}(B_{9/16}))}\\
        &\leq cE^p_{\mathrm{loc}}(v-l;Q^s_{3/4})\\
        &\quad+c\left(\dashint_{I^s_{3/4}}\mathrm{Tail}(u-l-(u-l)_{B_{3/4}}(t);B_{3/4})^p\,dt\right)^{\frac{1}p}+c|\mu|(Q^s_{1}),
    \end{aligned}
    \end{equation}
    holds for some constants $\sigma_1=\sigma_1(n,s_0,\Lambda,p)\in(0,1)$ and $c=c(n,s_0,\Lambda,p,q)$.
     Let us fix $h<s^{\frac1{2s}}/1000$ and $\beta=1/2$. As in Lemma \ref{lem.graex}, there are mutually disjoint coverings $\{B_{|h|^{\beta}}(x_{i})\}_{i\in\mathcal{I} }$ of $B_{5/8}$ and $\{I^s_{|h|^\beta}(t_j)\}_{j\in\mathcal{J}}$ of $I^s_{5/8}$ such that $(x_{i},t_j)\in Q^s_{5/8}$, \eqref{max.cov.ineq0}, \eqref{max.cov.ineq1} and \eqref{max.cov.ineq2} hold. In addition, there is a positive integer $m_0$ satisfying \eqref{el : choi.m0}. We observe from \eqref{defn.diff.t} and \eqref{hol.lin.ftn} that 
    \begin{equation}\label{ineq.fir.hss}
    \begin{aligned}
        \int_{Q^s_{5/8}}|\delta^t_h(v-l)|^q\,dz &\leq c\sum_{i,j}|h|^{(n+2s)\beta}\dashint_{Q^s_{|h|^\beta}(x_i,t_j)}|\delta^t_h(v-l)|^q\,dz \\
        &\leq c|h|^{(n+2s)\beta+q\gamma(1-\beta)}\sum_{i,j}E^p(v-l;Q^s_{4|h|^\beta}(x_i,t_j))^q,
    \end{aligned}
    \end{equation}
    where $c=c(n,s_0,\Lambda,p)$. We first note from \eqref{rel.exc} that
    \begin{align*}
        &E^p(v-l;Q^s_{4|h|^\beta}(x_i,t_j))^p\\
        &\leq c\dashint_{Q^s_{4|h|^\beta}(x_i,t_j)}|v-l-(v-l)_{Q^s_{11/16}}|^p\,dz\\
        &\quad+\dashint_{I^s_{4|h|^\beta}(t_j)}\mathrm{Tail}(v-l-(v-l)_{B_{4|h|^{\beta}}(x_i)}(t);B_{4|h|^\beta}(x_i))^p\,dt\eqqcolon J_1+J_2
    \end{align*}
    for some constant $c=c(n,s_0,p)$.
    In light of \cite[Lemma 2.2]{DieKimLeeNow24} along with H\"older's inequality, we next estimate $J_2$ as 
    \begin{align*}
        J_2&\leq c\sum_{k=2}^{m_0}2^{-2sk}\dashint_{I^s_{|h|^\beta}(t_j)}\dashint_{B_{2^{k}|h|^{\beta}}(x_i)}|v-l-(v-l)_{B_{2^{k}|h|^\beta}(x_{j})}(t)|^p\,dz\\
        &\quad +c2^{-2sm_0p}\dashint_{I^s_{4|h|^\beta}(t_j)}\mathrm{Tail}(v-l-(v-l)_{B_{11/16}}(t);B_{11/16})^p\,dt\\
        &\leq c\sum_{k=0}^{m_0}2^{-2sk}\dashint_{I^s_{|h|^\beta}(t_j)}\dashint_{B_{2^k|h|^{\beta}}(x_i)}|v-l-(v-l)_{Q^s_{11/16}}|^p\,dz\\
        &\quad+c\int_{I^s_{4|h|^\beta}(t_j)}\mathrm{Tail}(u-l-(u-l)_{B_{11/16}}(t);B_{11/16})^p\,dt+c|\mu|(Q^s_1)
        \end{align*}
        for some constant $c=c(n,s_0,\Lambda,p)$, where we have used \eqref{el : choi.m0} and Lemma \ref{lem.comp} for the last inequality.
        Combining the above two inequalities with \eqref{max.cov.ineq0}, \eqref{max.cov.ineq1}, \eqref{max.cov.ineq2},  H\"older's inequality and the fact that 
        \begin{align*}
            \sum_{i,j}|a_{i,j}|^q\leq \left(\sum_{i,j}|a_{i,j}|\right)^q \quad(q\geq 1),
        \end{align*} we get 
        \begin{align*}
            &\sum_{i\in\mathcal{I}}\sum_{j\in\mathcal{J}}E^p(v-l;Q^s_{4|h|^\beta}(x_i,t_j))^q\\
            &\leq c\sum_{i\in\mathcal{I}}\sum_{j\in\mathcal{J}}\sum_{k=2}^{m_0}2^{-2sk}\dashint_{I^s_{|h|^\beta}(t_j)}\dashint_{B_{2^k|h|^{\beta}}(x_i)}|v-l-(v-l)_{Q^s_{11/16}}|^q\,dz\\
            &\,\, +c\sum_{i\in\mathcal{I}}\left(\sum_{j\in\mathcal{J}}\int_{I^s_{4|h|^\beta}(t_j)}\mathrm{Tail}(u-l-(u-l)_{B_{11/16}}(t);B_{11/16})^p\,dt\right)^{\frac{q}p}
            +c\sum_{i\in\mathcal{I}}\sum_{j\in\mathcal{J}}(|\mu|(Q^s_1))^q\\
            &\leq c|h|^{-(n+2s)\beta}\int_{Q^s_{11/16}}|v-l-(v-l)_{Q^s_{11/16}}|^q\,dz\\
            &\,\,+ c|h|^{-n\beta}\left(\int_{{I^s_{11/16}}}\mathrm{Tail}(u-l-(u-l)_{B_{11/16}}(t);B_{11/16})^p\,dt\right)^{\frac{q}p}
            +|h|^{-(n+2s)\beta}(|\mu|(Q^s_1))^q
        \end{align*}
        for some constant $c=c(n,s_0,\Lambda,p)$.

        On the other hand, we note from Lemma \ref{lem.linmsol}, \eqref{rel.exc} and Lemma \ref{lem.comp} that
        \begin{align*}
            &\left(\int_{Q^s_{11/16}}|v-l-(v-l)_{Q^s_{11/16}}|^q\,dz\right)^{\frac1q}\\
            &\leq cE^p(v-l;Q^s_{3/4})\\
            &\leq cE^p_{\mathrm{loc}}(v-l;Q^s_{3/4})+c|\mu|(Q^s_1)+c\left(\int_{I^s_{3/4}}\mathrm{Tail}(u-l-(u-l)_{B_{3/4}}(t);B_{3/4})^p\,dt\right)^{\frac{1}p}
        \end{align*}
        for some constant $c=c(n,s_0,\Lambda,p,q)$.
        Using the above two inequalities together with \eqref{ineq.fir.hss}, 
        we obtain
        \begin{align*}
            \int_{Q^s_{5/8}}|\delta^t_h(v-l)|^q\,dz&\leq c|h|^{q\gamma(1-\beta)}\left(\int_{Q^s_{3/4}}|v-l-(v-l)_{Q^s_{3/4}}|^p\,dz\right)^{\frac{q}p}\\
            &\quad+c|h|^{q\gamma(1-\beta)}\left(\int_{I^s_{3/4}}\mathrm{Tail}(u-l-(u-l)_{B_{3/4}}(t);B_{3/4})^p\,dt\right)^{\frac{q}p}\\
            &\quad +c|h|^{q\gamma(1-\beta)}(|\mu|(Q^s_1))^q.
        \end{align*}
        We now use the embedding result given by Lemma \ref{lem.fouemb} to see 
        \begin{equation}\label{ineq.tq}
        \begin{aligned}
            [v-l]_{W^{\sigma_1/s,q}(I^s_{9/16};L^q(B_{9/16}))}&\leq cE^p_{\mathrm{loc}}(v-l;Q^s_{3/4})\\
            &\quad+c\left(\dashint_{{I^s_{3/4}}}\mathrm{Tail}(u-l-(u-l)_{B_{3/4}}(t);B_{3/4})^p\,dt\right)^{\frac1p}\\
            &\quad+c|\mu|(Q^s_1)
        \end{aligned}
        \end{equation}
        for some constant $c=c(n,s_0,\Lambda,p,q)$, where $\sigma_1=\gamma(1-\beta)/4$ depends only on $n,s_0,\Lambda$ and $p$. Similarly, by replacing $\delta_h^t(v-l)$ with $\delta_h (v-l)$ and following the same lines as in the proof of \eqref{ineq.tq}, we deduce 
        \begin{equation*}
        \begin{aligned}
            [v-l]_{L^q(I^s_{9/16};W^{2\sigma_1,q}(B_{9/16}))}&\leq cE^p_{\mathrm{loc}}(v-l;Q^s_{3/4})\\
            &\quad+c\left(\dashint_{{I^s_{3/4}}}\mathrm{Tail}(u-l-(u-l)_{B_{3/4}}(t);B_{3/4})^p\,dt\right)^{\frac1p}\\
            &\quad+c|\mu|(Q^s_1).
        \end{aligned}
        \end{equation*}
        Therefore, the proof of \eqref{osc.vminusl} is complete. We now fix $q=\frac{n+2s}{\sigma_1}$ to see that if $Q^s_r(z_1)\subset Q^s_{9/16}$, then
        \begin{align*}
            \dashint_{Q^s_r(z_1)}|(v-l)-(v-l)_{Q^s_r(z_1)}|\,dz&\leq cr^{\sigma_1}[v-l]_{W^{\sigma_1/s,q}(I^s_{1/2};L^q(B_{1/2}))}\\
            &\quad+cr^{\sigma_1}[v-l]_{L^{q}(I^s_{1/2};W^{2\sigma_1,q}(B_{1/2}))}
        \end{align*}
        holds for some constant $c=c(n)$, where we have used H\"older's inequality and Lemma \ref{lem.spwa}. By combining Lemma \ref{lem.camp} with the above inequality and \eqref{osc.vminusl}, we obtain the desired estimate. 
\end{proof}

Using the previous lemma along with fractional Gagliardo-Nirenberg-type inequalities given in \cite{BrezisM} yields the following reverse H\"older-type inequality.
\begin{lemma}\label{lem.rev}
    Let $v$ be a weak solution to \eqref{eq.ivp.dec}. Then for any $q \in [p,\infty)$, we have 
    \begin{equation}\label{est.redorder}
    \begin{aligned}
        &E^q_{\mathrm{loc}}(\nabla v;Q^s_{R/2}(z_0))\\
        &\leq cE^p_{\mathrm{loc}}(\nabla v;Q^s_{3R/4}(z_0))+cE^p(\nabla u;Q^s_{3R/4}(z_0))\\
        &\quad+\frac{c}{R}\left(\dashint_{I^s_{3R/4}(t_0)}\mathrm{Tail}(u-(\nabla u)_{Q^s_{3R/4}(z_0)}\cdot (y-x_0)-(u)_{B_{3R/4}(x_0)}(t);B_{3R/4}(x_0))^p\,dt\right)^{\frac1p}\\
        &\quad+c\frac{|\mu|(Q^s_R(z_0))}{R^{n+1}}
    \end{aligned}
    \end{equation}
    for some constant $c=c(n,s_0,\Lambda,p,q)$.
\end{lemma}
\begin{proof}
    We may assume $z_0=0$ and $R=1$. Recall the constant $\kappa=\kappa(n,s_0,\Lambda,p)\in(0,1)$ determined in Lemma \ref{lem.sobgrasp}.
    We fix the constant
    \begin{equation}\label{defn.theta}
        \theta=\frac{1+\kappa/2}{1+\kappa}.
    \end{equation}
    We next choose $\epsilon \coloneqq \frac{1}{(1-\theta)/(2q)+\theta/q}-q>0$    to see that
    \begin{align}\label{eq.y}
        \frac{1-\theta}{2q}+\frac{\theta}{q}=\frac{1}{q+\epsilon}.
    \end{align}
    Let us also choose 
    \begin{equation}\label{choi.ep11}
        \epsilon_1 \coloneqq \min\{\epsilon,\kappa/(2n)\}>\kappa/(16n).
    \end{equation}
    Then there is a positive integer $i_0=i_0(n,s_0,\Lambda,p,q)$ such that 
    \begin{align}\label{i00}
          p+i_0\epsilon_1\geq q.
    \end{align}
    We define a sequence by
    \begin{equation}\label{ris}
        r_i={1}/{2}+i/(32(i_0+2)),\quad 0\leq i\leq i_0+2
    \end{equation}
    and observe that $r_{i_0+2}=17/32$.
    We note from \cite[Theorem 1]{BrezisM} that
    \begin{align}\label{ineq.garni}
        \|g(\cdot,t)\|_{W^{1+\kappa/2,q}(B_{r_i})}\leq c\|g(\cdot,t)\|_{L^q(B_{r_i})}^{1-\theta}\|g(\cdot,t)\|_{W^{1+\kappa,q}(B_{r_i})}^{\theta}
    \end{align}
    for some $c=c(n,s_0,\Lambda,p)$,
    where $g\coloneqq v-l$ with $l\coloneqq (\nabla v)_{Q^s_{r_i}}\cdot x+(v)_{Q^s_{r_i}}$.
    Applying H\"older's inequality together with \eqref{eq.y} into \eqref{ineq.garni} yields
    \begin{equation}\label{ineq2.sst}
    \begin{aligned}
        &\left(\int_{I^s_{r_i}}\|g(\cdot,t)\|_{W^{1+\kappa/2,q}(B_{r_i})}^{q+\epsilon}\,dt\right)^{\frac1{q+\epsilon}}\\
        &\leq c\left(\int_{I^s_{r_i}}\|g(\cdot,t)\|_{L^q(B_{r_i})}^{2q}\,dt\right)^{\frac{1-\theta }{2q}}\left(\int_{I^s_{r_i}}\|g(\cdot,t)\|_{W^{1+\kappa,q}(B_{r_i})}^{q}\,dt\right)^{\frac{\theta}{q}}\\
        &\leq c\left(\int_{I^s_{r_i}}\|g(\cdot,t)\|_{L^{2{q}}(B_{r_i})}^{2{q}}\,dt\right)^{\frac{1-\theta }{2{q}}}\left(\int_{I^s_{r_i}}\|g(\cdot,t)\|_{W^{1+\kappa,q}(B_{r_i})}^{q}\,dt\right)^{\frac{\theta}{q}},
    \end{aligned}
    \end{equation}
    where $c=c(n,s_0,\Lambda,p,q)$. 
    
    Using the fact that $v-l$ is a solution to \eqref{eq.solmlin} and Lemma \ref{lem.glu}, we deduce that
    \begin{equation}\label{glu.gra}
    \begin{aligned}
        E^{2q}_{\mathrm{loc}}(v-l;Q^s_{r_i})&\leq c\dashint_{Q^s_{r_{i+1}}}|\nabla(v-l)|\,dz+c\dashint_{I^s_{r_{i+1}}}\mathrm{Tail}(v-l-(v-l)_{B_{r_{i+1}}}(t))\,dt\\
        &\leq cE_{\mathrm{loc}}(\nabla v;Q^s_{3/4})\\
        &\quad+c\dashint_{I^s_{r_{i+1}}}\mathrm{Tail}(v-(\nabla v)_{Q^s_{r_{i+1}}}\cdot y-(v)_{B_{r_{i+1}}}(t))\,dt
    \end{aligned}
    \end{equation}
    holds, where $c=c(n,s_0,\Lambda,p,q)$, as $r_{i+1}-r_i$ depends only on $n,s_0,\Lambda,p,q$. In addition, using Lemma \ref{lem.tail1} together with a few simple calculations, we obtain
    \begin{equation}\label{glu.tailest}
    \begin{aligned}
        &\dashint_{I^s_{r_{i+1}}}\mathrm{Tail}(v-(\nabla v)_{Q^s_{r_{i+1}}}\cdot y-(v)_{B_{r_{i+1}}}(t))\,dt\\
        &\leq c\dashint_{Q^s_{5/8}}|\nabla v-\nabla u|\,dz+c\dashint_{Q^s_{5/8}}|v-u|\,dz\\
        &\quad+c\dashint_{I^s_{5/8}}\mathrm{Tail}(u-(\nabla u)_{Q^s_{5/8}}\cdot y-(u)_{B_{5/8}}(t);B_{5/8})\,dt,
    \end{aligned}
    \end{equation}
   where $c=c(n,s_0,\Lambda,p,q)$. Combining \eqref{glu.gra} and \eqref{glu.tailest} with Lemma \ref{lem.comp}, Lemma \ref{lem.comp.gra}, Lemma \ref{lem.tail1} and H\"older's inequality, we obtain 
   \begin{equation}\label{ineq.solgrai}
   \begin{aligned}
       &E^{2q}_{\mathrm{loc}}(v-l;Q^s_{r_i})\\
       &\leq cE^p_{\mathrm{loc}}(\nabla v;Q^s_{3/4})+c|\mu|(Q^s_1)+cE^p(\nabla u;Q^s_{3/4})\\
       &\quad+c\left(\dashint_{I^s_{3/4}}\mathrm{Tail}(u-(\nabla u)_{Q^s_{3/4}}\cdot y-(u)_{B_{{3/4}}}(t);B_{3/4})^p\,dt\right)^{\frac1p},
   \end{aligned}
   \end{equation}
   where $c=c(n,s_0,\Lambda,p,q,i)$. Using \eqref{ineq2.sst} and \eqref{ineq.solgrai} along with H\"older's inequality, we obtain 
   \begin{equation}\label{indd1}
   \begin{aligned}
       &\left(\dashint_{I^s_{r_i}}\|(v-l)(\cdot,t)\|_{W^{1+\kappa/2,q}(B_{r_i})}^{q+\epsilon}\,dt\right)^{\frac1{q+\epsilon}}\\
       &\leq cE^q_{\mathrm{loc}}(\nabla v;Q^s_{r_{i+1}})+c|\mu|(Q^s_1)+cE^p(\nabla v;Q^s_{3/4})+cE^p(\nabla u;Q^s_{3/4})\\
       &\quad+c\left(\dashint_{I^s_{3/4}}\mathrm{Tail}(u-(\nabla u)_{Q^s_{3/4}}\cdot y-(u)_{B_{3/4}}(t);B_{3/4})^p\,dt\right)^{\frac1p}
   \end{aligned}
   \end{equation}
   for some constant $c=c(n,s_0,\Lambda,p,q)$, where we have also used \eqref{gra.mainest} with $R=1$, $r=r_i$, $\rho=r_{i+1}$ and $z_0=0$. By the fractional Sobolev embedding as in e.g.\ \cite[Theorem 6.7]{DinPalVal12} together with \eqref{choi.ep11} and \eqref{indd1}, we get 
   \begin{align*}
       E_{\mathrm{loc}}^{q+\epsilon_1}(\nabla v;Q^s_{r_i})&\leq \left(\dashint_{I^s_{r_{i}}}\|\nabla(v-l)(\cdot,t)\|_{W^{\kappa/2,q}(B_{r_i})}^{q+\epsilon_1}\,dt\right)^{\frac1{q+\epsilon_1}}\\
       &\leq cE^q_{\mathrm{loc}}(\nabla v;Q^s_{r_{i+1}})+c|\mu|(Q^s_1)+cE_{\mathrm{loc}}^p(\nabla v ;Q^s_{3/4})+cE^p(\nabla u;Q^s_{3/4})\\
       &\quad+c\left(\dashint_{I^s_{3/4}}\mathrm{Tail}(u-(\nabla u)_{Q^s_{3/4}}\cdot y-(u)_{B_{3/4}}(t);B_{3/4})^p\,dt\right)^{\frac1p}
   \end{align*}
   for some constant $c=c(n,s_0,\Lambda,p,q,i)$. Therefore, we have
   \begin{equation}\label{ind.ineq.q}
   \begin{aligned}
       E^{q_{k+1}}_{\mathrm{loc}}(\nabla v;Q^s_{r_i})&\leq cE^{q_k}_{\mathrm{loc}}(\nabla v;Q^s_{r_{i+1}})+c|\mu|(Q^s_1)+cE_{\mathrm{loc}}^p(\nabla v ;Q^s_{3/4})+cE^p(\nabla u;Q^s_{3/4})\\
       &\quad+c\left(\dashint_{I^s_{3/4}}\mathrm{Tail}(u-(\nabla u)_{Q^s_{3/4}}\cdot y-(u)_{B_{3/4}}(t);B_{3/4})^p\,dt\right)^{\frac1p}
   \end{aligned}
   \end{equation}
   for some constant $c=c(n,s_0,\Lambda,p,q,i)$, where we denote 
   \begin{align*}
       q_k=p+k\epsilon_1\quad\text{for any }k\geq0.
   \end{align*}
   By H\"older's inequality, \eqref{ind.ineq.q}, \eqref{i00} and \eqref{ris}, we arrive at
   \begin{align*}
       E^{q}_{\mathrm{loc}}(\nabla v;Q^s_{1/2})&\leq E^{q_{i_0+2}}_{\mathrm{loc}}(\nabla v;Q^s_{1/2})\\
       &\leq cE^{q_{i_0+1}}_{\mathrm{loc}}(\nabla v;Q^s_{r_{1}})+c|\mu|(Q^s_1)+cE_{\mathrm{loc}}^p(\nabla v ;Q^s_{3/4})+cE^p(\nabla u;Q^s_{3/4})\\
       &\quad+c\left(\dashint_{I^s_{3/4}}\mathrm{Tail}(u-(\nabla u)_{Q^s_{3/4}}\cdot y-(u)_{B_{3/4}}(t);B_{3/4})^p\,dt\right)^{\frac1p}\\
       &\leq cE^{p}_{\mathrm{loc}}(\nabla v;Q^s_{3/4})+c|\mu|(Q^s_1)+cE_{\mathrm{loc}}^p(\nabla v ;Q^s_{3/4})+cE^p(\nabla u;Q^s_{3/4})\\
       &\quad+c\left(\dashint_{I^s_{3/4}}\mathrm{Tail}(u-(\nabla u)_{Q^s_{3/4}}\cdot y-(u)_{B_{3/4}}(t);B_{3/4})^p\,dt\right)^{\frac1p}
   \end{align*}
   for some constant $c=c(n,s_0,\Lambda,p,q)$, as we iterate $i_0+1$ times and the positive integer $i_0$ depends only on $n,s_0,\Lambda,p$ and $q$. This completes the proof.  
\end{proof}

Using Lemma \ref{lem.higsotimesol} and Lemma \ref{lem.rev} in conjunction with Lemma \ref{lem.holspI} and Lemma \ref{lem.interpovec}, we now prove higher differentiability with respect to the time variable of the spatial gradient of solutions to homogeneous problems. 
\begin{lemma}\label{lem.higsobtimgra}
    Fix $p\in\left(1,\frac{n+2s_0}{n+1}\right)$, let $v$ be a weak solution to \eqref{eq.ivp.dec} and let $q \in [p,\infty)$. Then there is a constant $\varkappa=\varkappa(n,s_0,\Lambda,p)\in(0,1)$ independent of $q$ such that
    \begin{equation}\label{ineq.higsobtimgra}
    \begin{aligned}
        &R^{\varkappa}[\nabla v]_{W^{\frac{\varkappa}{2s},q}(I^s_{R/2}(t_0);L^q(B_{R/2}(x_0))}\\
        &\leq cE^p_{\mathrm{loc}}(\nabla v;Q^s_{3R/4}(z_0))+cE^p(\nabla u;Q^s_{3R/4}(z_0))+c\frac{|\mu|(Q^s_R(z_0))}{R^{n+1}}\\
        &\quad+\frac{c}{R}\left(\dashint_{I^s_{3R/4}(t_0)}\mathrm{Tail}(u-(\nabla u)_{Q^s_{3R/4}(z_0)}\cdot (y-x_0)-(u)_{B_{3R/4}(x_0)}(t);B_{3R/4}(x_0))^p\,dt\right)^{\frac1p}
    \end{aligned}
    \end{equation}
    for some constant $c=c(n,s_0,\Lambda,p,q)$.
\end{lemma}
\begin{proof}
    We may assume $z_0=0$ and $R=1$. 
    We are going to prove that there is a constant $\varkappa=\varkappa(n,s_0,\Lambda,p)\in(0,1)$ such that
    \begin{equation}\label{est.timesob}
    \begin{aligned}
        &\|\nabla v\|_{W^{\varkappa/s,q}(I^s_{1/2};L^q(B_{1/2}))}\\
        &\leq cE^p_{\mathrm{loc}}(\nabla v;Q^s_{3/4})+cE^p(\nabla u;Q^s_{3/4})+c|\mu|(Q^s_1)\\
        &\quad+c\left(\dashint_{I^s_{3/4}}\mathrm{Tail}(u-(\nabla u)_{Q^s_{3/4}}\cdot y-(u)_{B_{3/4}}(t);B_{3/4})^p\,dt\right)^{\frac1p}
    \end{aligned}
    \end{equation}
    holds for some constant $c=c(n,s_0,\Lambda,p,q)$. We divide the proof into two parts depending on the range of $q$.
    
    \textbf{Step 1: In case of $q\leq 2$.} 
    Let us fix some $\phi\in C_c^{\infty}(Q^s_{5/8})$ with $\phi\equiv 1$ on $Q^s_{1/2}$.
    In view of Lemma \ref{lem.sobgrasp} and Lemma \ref{lem.higsotimesol}, we have
    \begin{align}\label{aa}
        w\coloneqq (v-l)\phi&\in W^{\widetilde{\kappa}/2s,2}(\bbR;L^{{2}}(\bbR^n))\cap L^{2}(\bbR;W^{1+\widetilde{\kappa},2}(\bbR^n)),
    \end{align}
     where 
     \begin{equation}\label{cond.tkappa}
         \widetilde{\kappa}\coloneqq\min\{\kappa,\sigma_1/2\}\quad\text{and}\quad l\coloneqq(\nabla v)_{Q^s_{11/16}}\cdot y+(v)_{Q^s_{11/16}}.
     \end{equation} We point out that the constants $\kappa$ and $\sigma_1$ are determined in Lemma \ref{lem.sobgrasp} and Lemma \ref{lem.higsotimesol}, respectively. In addition, using Lemma \ref{lem.higsotimesol}, \eqref{ineq.solgrai}, Lemma \ref{lem.sobgrasp} and Lemma \ref{lem.rev} along with a few simple calculations, we obtain
    \begin{equation}\label{ineq.12g}
    \begin{aligned}
        \|w\|_{W^{\widetilde{\kappa}/2s,2}(\bbR;L^2(\bbR^n))}
        &\leq c\|v-l\|_{L^{2}(Q^s_{5/8})}+c[v-l]_{C^{\sigma_1}(Q^s_{5/8})}\\
        &\leq cE_{\mathrm{loc}}^p(v-l;Q^s_{11/16})+c|\mu|(Q^s_1)\\
        &\quad+c\left(\dashint_{I^s_{11/16}}\mathrm{Tail}(u-l-(u-l)_{B_{11/16}}(t);B_{11/16})^p\,dt\right)^{\frac1p}\\
        &\leq cM
    \end{aligned}
    \end{equation}
    and
    \begin{align*}
        \|w\|_{L^2(\bbR;W^{1+\widetilde{\kappa},2}(\bbR^n))}&\leq c\|v-l\|
        _{L^2(Q^s_{5/8})}+c[v-l]_{C^{\sigma_1}(Q^s_{5/8})}+c\|\nabla(v-l)\|_{L^2(I^s_{5/8};W^{\widetilde{\kappa},2}(B_{5/8}))}\\
        &\leq cM,
    \end{align*}
    where $c=c(n,s_0,\Lambda,p)$ and
    \begin{align*}
        M&\coloneqq E^p_{\mathrm{loc}}(\nabla v;Q^s_{3/4})+|\mu|(Q^s_1)+cE^p(\nabla u;Q^s_{3/4})\\
       &\quad+\left(\dashint_{I^s_{3/4}}\mathrm{Tail}(u-(\nabla u)_{Q^s_{3/4}}\cdot y-(u)_{B_{3/4}}(t);B_{3/4})^p\,dt\right)^{\frac1p}.
    \end{align*}
	By using the Fourier transform, we observe that
	\begin{align*}
		\|f\|_{W^{\kappa/2s,2}(\bbR;L^2(\bbR^n))}\eqsim \|(|\tau|+1)^{\kappa/2s}\mathcal{F}\,\widehat{f}(\xi,\tau)\|_{L^2(\bbR^{n+1})}
	\end{align*}
	and 
	\begin{align*}
		\|f\|_{L^{2}(\bbR;W^{1+\kappa,2}(\bbR^n))}\eqsim \|(|\xi|+1)^{1+\kappa}\mathcal{F}\,\widehat{f}(\xi,\tau)\|_{L^2(\bbR^{n+1})},
	\end{align*}
	where $\widehat{f}$ is the Fourier transform with respect to the spatial variables and $\mathcal{F}f$ is the Fourier transform with respect to the time variable.
    We next observe from \eqref{aa} and H\"older's inequality that
    \begin{align*}
        &\|\partial_i w\|_{W^{\varkappa_0/s,2}(\bbR^n;L^{2}(\bbR^n))}\\
        &\eqsim \|(|\tau|+1)^{{\varkappa_0}/s}\xi_i(\mathcal{F}\, \widehat{w})(\xi,\tau)\|_{L^{2}(\bbR^{n+1})} \\
        &\leq \|(|\tau|+1)^{\widetilde{\kappa}/2s}(\mathcal{F}\, \widehat{w})(\xi,\tau)\|_{L^{2}(\bbR^{n+1})}^{1-\vartheta}\|(|\xi_i|+1)^{1+\widetilde{\kappa}}(\mathcal{F}\, \widehat{w})(\xi,\tau)\|_{L^{2}(\bbR^{n+1})}^{\vartheta}\\
        &\leq c\|w\|_{W^{\widetilde{\kappa}/2s,2}(\bbR;L^2(\bbR^n))}^{1-\vartheta}\|w\|_{L^2(\bbR;W^{1+\widetilde{\kappa},2}(\bbR^n))}^{\vartheta},
    \end{align*}
    where 
    \begin{align}\label{cond.varkappa0}
        \varkappa_0\coloneqq\frac{\widetilde{\kappa}^2}{2(1+\widetilde{\kappa})}\quad\text{and}\quad\vartheta\coloneqq\frac{1}{1+\widetilde{\kappa}}.
    \end{align}
    Combining the above three estimates along with the fact that $w=(v-l)\varphi$ yields 
    \begin{equation}\label{est.p2}
    \begin{aligned}
        \|\nabla (v-\overline{l})\|_{W^{{\varkappa_0}/{s},2}(I^s_{1/2};L^2(B_{1/2}))}
        &\leq \|\nabla w\|_{W^{\varkappa_0/s,2}(\bbR;L^{2}(\bbR^n))}+E^p_{\mathrm{loc}}(\nabla v;Q^s_{3/4})\\
        &\leq cM
    \end{aligned}
    \end{equation}
    for some constant $c=c(n,s_0,\Lambda,p)$, where we write 
     \begin{align*}
     \overline{l}(y,t)\coloneqq(\nabla v)_{Q^s_{3/4}}\cdot y+(v)_{Q^s_{3/4}}.
    \end{align*}
    If $q\leq 2$, then by H\"older's inequality, we deduce
    \begin{equation}\label{ineq.fir.t}
    \begin{aligned}
        \|\nabla (v-\overline{l})\|_{W^{\varkappa_0/2s,q}(I^s_{1/2};L^q(B_{1/2}))}
        \leq [\nabla v]_{W^{\varkappa_0/s,2}(I^s_{1/2};L^2(B_{1/2}))}\leq cM
    \end{aligned}
    \end{equation}
    for some constant $c=c(n,s_0,\Lambda,p,q)$. Thus, by taking $\varkappa\leq \varkappa_0$, we obtain \eqref{est.timesob}.

    \textbf{Step 2: In case of $q>2$.} 
    We now use the fractional Gagliardo-Nirenberg-type inequality given by \cite[Theorem 1]{BrezisM} and H\"older's inequality to see that 
    \begin{equation}\label{mis}
    \begin{aligned}
        &\|\nabla (v-\overline{l})\|_{W^{\varkappa_0/(sq),q}(I^s_{1/2};L^q(B_{1/2}))}\\
        &=\left(\int_{B_{1/2}}\|\nabla (v-\overline{l})(x,\cdot)\|_{W^{\varkappa_0/(qs),q}(I^s_{1/2})}^{q}\,dx\right)^{\frac1q}\\
        &\leq c\left(\int_{B_{1/2}}\|\nabla (v-\overline{l})(x,\cdot)\|_{W^{\varkappa_0/s,2}(I^s_{1/2})}\|\nabla (v-\overline{l})(x,\cdot)\|_{L^{2(q-1)}(I^s_{1/2})}^{q-1}\,dx\right)^{\frac1q}\\
        &\leq c\|\nabla (v-\overline{l})\|_{W^{\varkappa_0/s,2}(I^s_{1/2};L^2(B_{1/2}))} ^{1/q} \|\nabla (v-\overline{l})\|_{L^{2(q-1)}(Q^s_{1/2})}^{(q-1)/q}
    \end{aligned}
    \end{equation}
    for some constant $c=c(n,\kappa,q)$.
    By \eqref{est.p2} and Lemma \ref{lem.rev}, we further estimate the right-hand side of \eqref{mis} as 
    \begin{align}\label{ineq.sec.t}
      &\|\nabla (v-\overline{l})\|_{W^{\varkappa_0/(sq),q}(I^s_{1/2};L^q(B_{1/2}))}\leq cM,
    \end{align}
   where $c=c(n,s_0,\Lambda,p,q)$. We now choose 
   \begin{align}\label{q0.cond.hight}
       q_0\coloneqq \frac{10n(2+\varkappa_0)}{\varkappa_0^2(1-\varkappa_0)}.
   \end{align}
   Then by taking 
   \begin{equation}\label{choi.varkappa}
       \varkappa\leq \varkappa_0/q_0,
   \end{equation}
   we obtain \eqref{ineq.higsobtimgra} in the case when $q\leq q_0$. 
   Therefore, we may assume that $q>q_0$. First, we are going to prove that
   \begin{align}\label{fin.goalfort}
       \|v-\overline{l}\|_{W^{\varkappa_0/4s,q}(I^s_{1/2};W^{\varkappa_0/2,q}(B_{1/2}))}+\|v-\overline{l}\|_{W^{\varkappa_0/2sq,q}(I^s_{1/2};W^{1+\varkappa_0^2/4,q}(B_{1/2}))}\leq cM
   \end{align}
   for some constant $c=c(n,s_0,\Lambda,p,q)$.
   Using the same reasoning as in \eqref{ineq.12g} along with \eqref{cond.tkappa}, \eqref{cond.varkappa0}, Lemma \ref{lem.holspI} and Lemma \ref{lem.higsotimesol}, we obtain
   \begin{equation} \label{eq:vest}
       \|v-\overline{l}\|_{W^{\varkappa_0/4s,q}(I^s_{1/2};W^{\varkappa_0/2,q}(B_{1/2}))}\leq [v-l]_{C^{2\varkappa_0}(Q^s_{1/2})}+\|v-\overline{l}\|_{L^q(Q^s_{1/2})}
       \leq cM.
   \end{equation}
   Combining \cite[Equation (2.6)]{DieKimLeeNow24d}, Lemma \ref{lem.sobgrasp}, Lemma \ref{lem.rev},  and \eqref{ineq.sec.t} yields
   \begin{equation}\label{ineq.bddvv}
   \begin{aligned}
       &\|v-\overline{l}\|_{W^{\varkappa_0/2sq,q}(I^s_{1/2};W^{\varkappa_0^2/4,q}(B_{1/2}))}+\|\nabla (v-\overline{l})\|_{L^q(Q^s_{1/2})}\\
       &\quad+[\nabla(v-\overline{l})]_{L^q(I^s_{1/2};W^{\varkappa_0^2/4,q}(B_{1/2}))}+[\nabla(v-\overline{l})]_{W^{\varkappa_0/(2sq),q}(I^s_{1/2};L^{q}(B_{1/2}))}\\
       &\leq 
       c\|v-\overline{l}\|_{W^{\varkappa_0/4s,q}(I^s_{1/2};W^{\varkappa_0/2,q}(B_{1/2}))}+cE^q_{\mathrm{loc}}(\nabla v;Q^s_{3/4})\\
       &\quad+ c[\nabla(v-\overline{l})]_{L^q(I^s_{1/2};W^{\varkappa_0,q}(B_{1/2}))}+c[\nabla(v-\overline{l})]_{W^{\varkappa_0/(sq),q}(I^s_{1/2};L^{q}(B_{1/2}))} \\
       &\leq cM.
   \end{aligned}
   \end{equation}
   Next, we estimate the follwing quantity
   \begin{align*}
    J\coloneqq\left(\int_{I^s_{1/2}}\int_{I^s_{1/2}}\int_{B_{1/2}}\int_{B_{1/2}}\frac{|(G(x,t)-G(x,\tau))-(G(y,t)-G(y,\tau))|^q}{|x-y|^{n+q\varkappa_0^2/4}|t-\tau|^{1+\varkappa_0/(2s)}}\,dZ\right)^{\frac1q},
\end{align*}
where we write 
\begin{align*}
    G\coloneqq \nabla(v-\overline{l})\quad\text{and}\quad \,dZ=\,dx\,dy\,dt\,d\tau.
\end{align*}
We now choose
\begin{align*}
    \vartheta_0=\frac{\varkappa_0/s}{2(1+\varkappa_0/s)},
\end{align*}
so that
\begin{align*}
    (1+\varkappa_0/(2s))/(1-\vartheta_0)=1+\varkappa_0/s.
\end{align*}
By \eqref{q0.cond.hight}, we see
\begin{align*}
    (n+q\varkappa_0^2/4)/\vartheta_0<n+q\varkappa_0.
\end{align*}
Thus, H\"older's inequality along with a few simple calculations and \eqref{ineq.bddvv} yields
\begin{align*}
    J&\leq c\left(\int\frac{|G(x,t)-G(y,t)|^q}{|x-y|^{(n+q\varkappa_0^2/4)/\vartheta_0}}\,dZ\right)^{\vartheta_0}\left(\int\frac{|G(x,t)-G(x,\tau)|^q}{|t-\tau|^{(1+\varkappa_0/(2s)/(1-\vartheta_0)}}\,dZ\right)^{1-\vartheta_0}\\
     &\leq c\left(\int\frac{|G(x,t)-G(y,t)|^q}{|x-y|^{(n+q\varkappa_0)}}\,dZ\right)^{\vartheta_0}\left(\int\frac{|G(x,t)-G(x,\tau)|^q}{|t-\tau|^{1+\varkappa_0/s}}\,dZ\right)^{1-\vartheta_0}\\
     &\leq c[\nabla(v-\overline{l})]_{L^q(I^s_{1/2};W^{\varkappa_0,q}(B_{1/2}))}^{\vartheta_0}[\nabla(v-\overline{l})]_{W^{\varkappa_0/(sq),q}(I^s_{1/2};L^{q}(B_{1/2}))}^{1-\vartheta_0}\leq cM
\end{align*}
for some constant $c=c(n,s_0,\Lambda,p,q)$. Therefore, combining the previous display with \eqref{eq:vest} and \eqref{ineq.bddvv} yields \eqref{fin.goalfort}.

  Using \cite[Equation (2.6)]{DieKimLeeNow24d} and combining the interpolation inequality given by Lemma \ref{lem.interpovec} with \eqref{fin.goalfort}, we obtain
   \begin{align*}
       & \|v-\overline{l}\|_{W^{\gamma,q}(I^s_{1/2};W^{1+\varkappa_1,q}(B_{1/2}))}\\
       &\leq  c\|v-\overline{l}\|_{W^{\varkappa_0/(4s),q}(I^s_{1/2};W^{\varkappa_0/2,q}(B_{1/2}))}^{\Theta}\|v-\overline{l}\|_{W^{\varkappa_0/2sq,q}(I^s_{1/2};W^{1+\varkappa_0^2/4,q}(B_{1/2}))}^{1-\Theta}\leq cM
   \end{align*}
   for some constant $c=c(n,s_0,\Lambda,p,q)$, where we write 
   \begin{align*}
     \varkappa_1\coloneqq  \frac{\varkappa_0^2}{8}+\frac{\varkappa_0^3/16}{1+\varkappa_0^2/4},\quad\Theta\coloneqq\frac{\varkappa_0^2/8}{1+\varkappa_0^2/4}\quad\text{and}\quad \gamma\coloneqq\frac{\varkappa_0^3}{8s(4+\varkappa_0^2)}+\frac{\varkappa_0(8+\varkappa_0^2)}{4sq(4+\varkappa_0^2)}.
   \end{align*}

   We now use \cite[Equation (2.6)]{DieKimLeeNow24d} to see that if 
    \begin{equation}\label{varkappa.cond}
       \varkappa\leq\varkappa_0^3/64,
   \end{equation}
   then we have
   \begin{equation}\label{ineq.thi.t}
   \begin{aligned}
       [\nabla v]_{W^{\varkappa,q}(I^s_{1/2};L^q(B_{1/2}))} \leq c\|v-\overline{l}\|_{W^{\gamma,q}(I^s_{1/2};W^{1+\varkappa_1,q}(B_{1/2}))}\leq cM,
   \end{aligned}
   \end{equation}
   where $c=c(n,s_0,\Lambda,p,q)$. This completes the proof of \eqref{est.timesob} when $q>q_0$ and \eqref{varkappa.cond} holds.
   Therefore, if we take $\varkappa=\min\{\varkappa_0/q_0,\varkappa_0^3/32\}$, where $\varkappa_0=\varkappa_0(n,s_0,\Lambda,p)$ and $q_0=q_0(n,s_0,\Lambda,p)$ are determined in \eqref{cond.varkappa0} and \eqref{q0.cond.hight}, respectively, then the estimate \eqref{est.timesob} follows by combining \eqref{ineq.fir.t}, \eqref{ineq.sec.t}, \eqref{choi.varkappa}, \eqref{varkappa.cond} and \eqref{ineq.thi.t}.
\end{proof}

Combining Lemma \ref{lem.sobgrasp}, Lemma \ref{lem.higsotimesol} and Lemma \ref{lem.rev}, we are finally able to deduce suitable decay estimates for homogeneous problems at the gradient level.
\begin{lemma}\label{lem.exc.dec.gra}
    Fix $p\in\left(1,\frac{n+2s_0}{n+1}\right)$ and let
    \begin{equation*}
        u\in L^2(I^s_{2R}(t_0);W^{s,2}(B_{2R}(x_0)))\cap L^p(I^s_{2R}(t_0);W^{1,p}(\bbR^n)).
    \end{equation*} Let $v\in L^2(I^s_{R}(t_0);W^{s,2}(B_{2R}(x_0)))\cap L^{p}(I^s_{R}(t_0);L^{1}_{2s}(\bbR^n))$ be a unique solution to \eqref{eq.ivp.dec}. There is a constant $\alpha_0=\alpha_0(n,s_0,\Lambda,p)\in(0,1)$ such that if $\rho\in(0,1/4]$, then
\begin{equation}\label{gra.mainestf}
\begin{aligned}
    \osc_{Q^s_{\rho R}(z_0)}\nabla v
    &\leq c\rho^{\alpha_0}\left[E^p_{\mathrm{loc}}(\nabla v;Q^s_{3R/4}(z_0))+ cE^p(\nabla u;Q^s_{3R/4}(z_0))\right]\\
    &\quad +c\rho^{\alpha_0}\left(\dashint_{I^s_{3R/4}(t_0)}\mathrm{Tail}\left(\frac{u-l}{R};B_{3R/4}(x_0)\right)^p\,dt\right)^{\frac1p}\\
    &\quad +c\rho^{\alpha_0}\frac{|\mu|(Q^s_R(z_0))}{R^{n+1}}
\end{aligned}
\end{equation}
holds for any $q\in [1,\infty)$, where $c=c(n,s_0,\Lambda,p,q)$ and 
\begin{equation*}
    l(y,t)=(\nabla u)_{Q^s_{3R/4}(z_0)}\cdot (y-x_0)+(u)_{B_{3R/4}(x_0)}(t).
\end{equation*}
\end{lemma}
\begin{proof}
    The proof relies on Campanato's characterization of H\"older spaces. We are first going to prove that there are constants $c=c(n,s_0,\Lambda,p)$ and $\alpha_0=\alpha_0(n,s_0,\Lambda,p)\in(0,1)$ such that
    \begin{equation}\label{ineq.camp}
    \begin{aligned}
        \dashint_{Q^s_r(z_1)}|\nabla v-(\nabla v)_{Q^s_r(z_1)}|\,dz&\leq c\left(\dfrac{r}{R}\right)^{\alpha_0}\left[E^p_{\mathrm{loc}}(\nabla v;Q^s_{3R/4}(z_0))+ cE^p(\nabla u;Q^s_{3R/4}(z_0))\right]\\
    &\quad +c\left(\dfrac{r}{R}\right)^{\alpha_0}\left(\dashint_{I^s_{3R/4}(t_0)}\mathrm{Tail}\left({\dfrac{u-l}{R};B_{3R/4}(x_0)}\right)^p\,dt\right)^{\frac1p}\\
    &\quad +c\left(\dfrac{r}{R}\right)^{\alpha_0}\frac{|\mu|(Q^s_R(z_0))}{R^{n+1}}\eqqcolon c\left(\dfrac{r}{R}\right)^{\alpha_0}M
    \end{aligned}
    \end{equation}
    holds whenever $Q^s_r(z_1)\subset Q^s_{R/2}(z_0)$. We now choose 
    \begin{align*}
        \overline{\kappa}=\min\{\kappa,\varkappa\},
    \end{align*}
    where the constants $\kappa$ and $\varkappa$ are determined in Lemma \ref{lem.higsobtimgra} and Lemma \ref{lem.higsotimesol}, respectively. Let us fix the constant $q$ such that 
    \begin{equation*}
        q=(n+2s)/2\overline{\kappa}
    \end{equation*}
    to see that 
    \begin{equation}\label{chhhh}
        \overline{\kappa}-\frac{n+2s}{q}=\overline{\kappa}/2.
    \end{equation}
    By H\"older's inequality and Lemma \ref{lem.spwa} along with the fact that $Q^s_r(z_1)\subset Q^s_{R/2}(z_0)$, we have 
    \begin{align*}
        \dashint_{Q^s_r(z_1)}|\nabla v-(\nabla v)_{Q^s_r(z_1)}|\,dz&\leq E^q_{\mathrm{loc}}(\nabla v;Q^s_{r}(z_1)) \\
         &\leq cr^{\overline{\kappa}-\frac{n+2s}{q}}[\nabla v]_{W^{\frac{\overline{\kappa}}{2s},q}(I^s_{R/2}(t_0);L^q(B_{R/2}(x_0)))}\\
        &\quad+cr^{\overline{\kappa}-\frac{n+2s}{q}}[\nabla v]_{L^q(I^s_{R/2}(t_0);W^{\overline{\kappa},q}(B_{R/2}(x_0)))}
    \end{align*}
    for some constant $c=c(n)$.
In addition, by \eqref{chhhh} and applying the estimates given in Lemma \ref{lem.higsobtimgra} and Lemma \ref{lem.higsotimesol} into the right-hand side of the above inequality, we obtain \eqref{ineq.camp} with $\alpha_0=\overline{\kappa}/2$. Finally, we deduce from Lemma \ref{lem.camp} and \eqref{ineq.camp} that 
\begin{equation*}
    R^{\alpha_0}[\nabla v]_{C^{\alpha_0}(Q^s_{R/4}(z_0))}\leq c\left(M+E_{\mathrm{loc}}(\nabla v;Q^s_{3R/4}(z_0))\right)
\end{equation*}
for some constant $c=c(n,s_0,\Lambda,p)$. Together with an application of H\"older's inequality, the desired estimate follows, finishing the proof.
\end{proof}

\subsection{First-order comparison estimates}
By combining Lemma \ref{lem.comp}, Lemma \ref{lem.highsu} and Lemma \ref{lem.sobgrasp} with an interpolation argument, we are able to upgrade our zero-order comparison estimates from Lemma \ref{lem.comp} to the gradient level.
\begin{lemma}\label{lem.comp.gra}
	Fix $p\in\left(1,\frac{n+2s_0}{n+1}\right)$, let $\mu\in L^1(I^s_{2R}(t_0);L^\infty(B_{2R}(x_0))$ and assume that
	$$
	u \in L^2(I^s_{2R}(t_0);W^{s,2}(B_{2R}(x_0)))\cap C(I^s_{2R}(t_0);L^2(B_{2R}(x_0))) \cap L^{p}(I^s_{2R}(t_0);L^{1}_{2s}(\bbR^n))
	$$
	is a weak solution to \eqref{eq:nonlocaleq} in $Q^s_{2R}(z_0)$ with $\nabla u\in L^p(\bbR^n\times I^s_{2R}(t_0))$.
	Let \begin{equation*}
		v\in L^2(I^s_{R}(t_0);W^{s,2}(B_{R}(x_0)))\cap C(I^s_R(t_0);L^2(B_R(x_0))) \cap L^{p}(I^s_{R}(t_0);L^{1}_{2s}(\bbR^n))
	\end{equation*}
	be the unique solution to \eqref{eq.ivp.dec}. Then we have the comparison estimate 
	\begin{align*}
		\left(\dashint_{Q^s_{R/2}(z_0)}|\nabla u-\nabla v|^p\,dz\right)^{\frac1p}&\leq c\frac{|\mu|(Q^s_R(z_0))}{R^{n+1}}\\
		&\quad+c\left(\frac{|\mu|(Q^s_R(z_0))}{R^{n+1}}\right)^{1-\theta}E^p(\nabla u ;Q^s_R(z_0))^\theta,
	\end{align*}
	where $\theta=\theta(n,s_0,\Lambda,p) \in (0,1)$ and $c=c(n,s_0,\Lambda,p)$.
\end{lemma}
\begin{proof}
	Let us fix $1/2\leq r<\rho\leq 3/4$ and choose
	\begin{equation*}
		\widetilde{\kappa}=\min\{\sigma_0,\kappa\},
	\end{equation*}
	where the constants $\sigma_0$ and $\kappa$ are determined in Lemma \ref{lem.highsu} and Lemma \ref{lem.sobgrasp}, respectively. We next select $\theta=1/(1+\widetilde{\kappa})$.  By \cite[Theorem 1]{BrezisM}, we obtain
	\begin{align*}
		&\|\nabla( u -v)(\cdot,t)\|_{L^p(B_r)}\\
		&\leq \|(u-v)(\cdot,t)\|_{L^p(B_r)}^{1-\theta} \left(\|(u-v)(\cdot,t)\|_{L^p(B_r)}+\|\nabla (u-v)(\cdot,t)\|_{L^p(B_r)}\right)^{\theta}\\
		&\quad +c\|(u-v)(\cdot,t)\|_{L^p(B_r)}^{1-\theta }[\nabla (u-v)(\cdot,t)]_{W^{\widetilde{\kappa},p}(B_r)}^{\theta}
	\end{align*}
	a.e. $t\in I^s_r$, where $c=c(n,s_0,\Lambda,p)$, as $\widetilde{\kappa}$ depends only on $n,s_0,\Lambda$ and $p$.
	After a few simple calculations together with integrating both sides of the above inequality with respect to the time variables and H\"older's inequality, we obtain
	\begin{align*}
		\|\nabla u-\nabla v\|_{L^p(Q^s_r)}&\leq \|u-v\|_{L^p(Q^s_r)}^{1-\theta} \left(\|u-v\|_{L^p(Q^s_r)}+\|\nabla u-\nabla v\|_{L^p(Q^s_r)}\right)^{\theta}\\
		&\quad +c\|u-v\|_{L^p(Q^s_r)}^{1-\theta }[\nabla u-\nabla v]_{L^p(I^s_r;W^{\widetilde{\kappa},p}(B_r))}^{\theta}.
	\end{align*}
	Using Young's inequality, we have 
	\begin{align*}
		\|\nabla u-\nabla v\|_{L^p(Q^s_r)}&\leq c\|u-v\|_{L^p(Q^s_r)}+c\|u-v\|_{L^p(Q^s_r)}^{1-\theta }[\nabla u-\nabla v]_{L^p(I^s_r;W^{\widetilde{\kappa},p}(B_r))}^{\theta}.
	\end{align*}
	We further estimate the above right-hand side as 
	\begin{align*}
		& \|\nabla u-\nabla v\|_{L^p(Q^s_r)} \\ &\leq c|\mu|(Q^s_1)+c(|\mu|(Q^s_1))^{1-\theta}\left[[\nabla u]_{L^p(I^s_r;W^{\widetilde{\kappa},p}(B_r))}+[\nabla v]_{L^p(I^s_r;W^{\widetilde{\kappa},p}(B_r))}\right]^{\theta}\\
		&\leq c|\mu|(Q^s_1)+\frac{c(|\mu|(Q^s_1))^{1-\theta}}{(\rho-r)^{n+4s}}\left[E^p_{\mathrm{loc}}(\nabla v; Q^s_\rho)+E^p(\nabla u ;Q^s_1)+ |\mu|(Q^s_1)\right]^{\theta}\\
		&\leq c|\mu|(Q^s_1)+\frac{c(|\mu|(Q^s_1))^{1-\theta}}{(\rho-r)^{n+4s}}\left[ \|\nabla u-\nabla v\|_{L^p(Q^s_\rho)}+E^p(\nabla u ;Q^s_1)+ |\mu|(Q^s_1)\right]^{\theta}\\
		&\leq \frac{ \|\nabla u-\nabla v\|_{L^p(Q^s_\rho)}}{2}+\frac{c\left[|\mu|(Q^s_1)+(|\mu|(Q^s_1))^{1-\theta}E^p(\nabla u ;Q^s_1)^{\theta}\right]}{(\rho-r)^{(n+4s)/(1-\theta)}},
	\end{align*}
	where we have used Lemma \ref{lem.sobgrasp} and Young's inequality.
	By employing a classical iteration lemma given in \cite[Lemma 6.1]{Giu03}, we arrive at the desired estimate.
\end{proof}

\section{Gradient potential estimates} \label{sec6}
In this section, we finally establish our pointwise gradient estimates in terms of caloric Riesz potentials. Let us recall \eqref{choi.s}. First, we fix some
\begin{equation}\label{Defn.pp}
    p\in\left(1,\frac{n+2s_0}{n+1}\right).
\end{equation}
Next, we select
\begin{equation}\label{defn.qq}
    q\coloneqq\frac{4s_0}{2s_0-1}
\end{equation}to see that $q$ depends only on $n,s_0$  and
\begin{align}\label{defn.q}
     2s-1-\frac{2s}q\geq\frac{2s_0-1}{2}.
\end{align}

\subsection{Excess decay}
We now introduce the functional 
\begin{align*}
    \overline{E}(u;Q^s_R(z_0))&\coloneqq\left(\dashint_{I^s_R(t_0)}\left(\dashint_{B_R(x_0)}\frac{|u-(\nabla u)_{Q^s_R(z_0)}\cdot (x-x_0)- (u)_{B_R(x_0)}(t)|}{R}\,dx\right)^q\,dt\right)^{\frac1q}\\
    &\,+\left(\dashint_{I^s_R(t_0)}\mathrm{Tail}\left(\frac{u-(\nabla u)_{Q^s_R(z_0)}\cdot (y-x_0)- (u)_{B_R(x_0)}(t)}R;B_R(x_0)\right)^q\,dt\right)^{\frac1q},
\end{align*}
which will be used to handle the parabolic tail that appears in the second line of the right-hand side in \eqref{gra.mainestf}.
\begin{remark} \normalfont
    We point out that in light of the localization argument given in Lemma \ref{lem.loc}, it is not restrictive for our purposes to always assume that $u\in L^q(I^s_R(t_0);L^1_{2s}(\bbR^n))$ for any $q\in[1,\infty]$.
\end{remark}
We also define the following modified parabolic nonlocal excess functional
\begin{equation}\label{Defn.exc}
    {E}(u,\nabla; Q^s_R(z_0))\coloneqq E^p(\nabla u;Q^s_R(z_0))+\overline{E}(u;Q^s_R(z_0)).
\end{equation}
\begin{remark} \normalfont
    Let us show that the following basic inequality
    \begin{align}\label{bas.ee}
    {E}(u,\nabla; Q^s_{\rho R}(z_0))\leq c\rho^{-(n+2s+1)} {E}(u,\nabla; Q^s_R(z_0))
\end{align}
holds for some constant $c=c(n,s_0,p)$ and any $\rho\in(0,1]$. We first note that in view of simple calculations, we have
    \begin{align}\label{111}
    E^p(\nabla u;Q^s_{\rho R}(z_0))\leq c\rho^{-(n+2s)}E^p(\nabla u;Q^s_{ R}(z_0))
\end{align}
for some constant $c=c(n,s_0,p)$. On the other hand, we observe that
\begin{align*}
    \overline{E}(u;Q^s_{\rho R}(z_0))&\leq \left(\dashint_{I^s_{\rho R}(t_0)}\left(\dashint_{B_{\rho R}(x_0)}\frac{|u-l-(u-l)_{B_{\rho R}(x_0)}(t)|}{\rho R}\,dx\right)^q\,dt\right)^{\frac1q}\\
    &\quad+\left(\dashint_{I^s_{\rho R}(t_0)}\mathrm{Tail}\left(\frac{u-l-(u-l)_{B_{\rho R}(x_0)}(t)}{\rho R};B_{\rho R}(x_0)\right)^q\,dt\right)^{\frac1q}\\
    &\quad+c|(\nabla u)_{Q^s_{\rho R}(z_0)}-(\nabla u)_{Q^s_{R}(z_0)}|
\end{align*}
for some constant $c=c(n,s_0)$, where we denote 
\begin{equation*}
    l(y,t)=(\nabla u)_{Q^s_{R}(z_0)}\cdot (y-x_0)+ (u)_{B_{ R}(x_0)}(t).
\end{equation*}
After a few calculations, we arrive at
\begin{align}\label{222}
    \overline{E}(u;Q^s_{\rho R}(z_0))\leq c\rho^{-(n+2s+1)}\left[\overline{E}(u;Q^s_{R}(z_0))+E(\nabla u;Q^s_R(z_0))\right]
\end{align}
for some constant $c=c(n,s_0,p)$. Thus, by \eqref{111} and \eqref{222}, we have \eqref{bas.ee}.
\end{remark}

We now prove decay estimates of $E(u,\nabla;\cdot)$ defined in \eqref{Defn.exc}.
\begin{lemma}\label{lem.exc.dec.fin}
    Let 
    \begin{equation*}
        u\in L^2(I^s_{R}(t_0);W^{s,2}(B_{R}(x_0)))\cap L^p(I^s_{R}(t_0);W^{1,p}(\bbR^n))\cap L^{q}(I^s_{R}(t_0);L^1_{2s}(\bbR^n))
    \end{equation*}
    be a weak solution to \eqref{eq:nonlocaleq} in $Q^s_{R}(z_0)$, where the constants $p$ and $q$ are determined in \eqref{Defn.pp} and \eqref{defn.qq}, respectively. For any $\rho\in(0,1]$, we have 
\begin{equation}\label{exc.ineq}
\begin{aligned}
    {E}(u,\nabla;Q^s_{\rho R}(z_0))
    &\leq c\rho^{\alpha_1}{E}(u,\nabla;Q^s_{ R}(z_0))\\
    &\quad +c\rho^{-(n+2s+1)}\left(\frac{|\mu|(Q^s_R(z_0))}{R^{n+1}}\right)^{1-\theta}{E}(u,\nabla;Q^s_{ R}(z_0))^\theta\\
    &\quad+c\rho^{-(n+2s+1)}\frac{|\mu|(Q^s_R(z_0))}{R^{n+1}},
\end{aligned}
\end{equation}
where $\alpha_1=\alpha_1(n,s_0,\Lambda,p)\in(0,1)$ and $c=c(n,s_0,\Lambda,p)$. In particular, the constant $\theta=\theta(n,s_0,\Lambda,p) \in (0,1)$ is determined in Lemma \ref{lem.comp.gra}.
\end{lemma}
\begin{proof}
    We may assume $R=1$ and $z_0=0$. 
    Let us fix 
    \begin{equation}\label{choi.alpha1}
        \alpha_1 \coloneqq \min\{2s_0(p-1)/p,\alpha_0,(2s_0-1)/4,2s_0-1-2s_0/q\}>0,
    \end{equation}
    where the constant $\alpha_0$ is determined in Lemma \ref{lem.exc.dec.gra}. 
    If $\rho\in[2^{-6},1]$, \eqref{exc.ineq} directly follows by \eqref{bas.ee}.
    Let $v$ be a weak solution to \eqref{eq.ivp} with $R$ and $z_0$ replaced by $1/2$ and $0$, respectively.
    We now assume $\rho< 2^{-6}$. 
    Note that there exists a natural number $N_\rho$ such that 
    \begin{equation}\label{choi.nrho}
        2^{-5}< 2^{N_\rho}\rho\leq 2^{-4}.
    \end{equation}
    In view of \eqref{rel.exc} and \cite[Lemma 2.2]{DieKimLeeNow24} together with a simple calculation, we have
    \begin{align*}
        {E}^p(\nabla u ;Q^s_\rho)&\leq c\left(\dashint_{I^s_\rho}\left(\sum_{i=0}^{N_\rho}2^{-2si}\dashint_{B_{2^i\rho}}|\nabla u-(\nabla u)_{Q^s_{2^i\rho}}|\,dx\right)^p\,dt\right)^{\frac1p}\\
        &\quad+c\left(\dashint_{I^s_\rho}\left(2^{-2sN_\rho}\mathrm{Tail} \left (\nabla u-(\nabla u)_{Q^s_{2^{N_\rho}\rho}} ; B_{2^{N_\rho} \rho} \right )\right)^p\,dt\right)^{\frac1p}\eqqcolon J_1+J_2
    \end{align*}
    for some constant $c=c(n,s_0)$. Using \eqref{choi.nrho}, \eqref{111} and the fact that by \eqref{choi.alpha1} we have $\alpha_1\leq 2s_0(p-1)/p$, we observe
    \begin{align*}
        J_2\leq c\rho^{{2s(p-1)}/{p}}E^p(\nabla u ;Q^s_1)\leq c\rho^{\alpha_1}E^p(\nabla u;Q^s_1),
    \end{align*}
    where $c=c(n,s_0)$. We next observe that
    \begin{align*}
        J_1&\leq c\left(\dashint_{I^s_\rho}\left(\sum_{i=0}^{N_\rho}2^{-2si}\dashint_{B_{2^i\rho}}|\nabla u-\nabla v|\,dx\right)^p\,dt\right)^{\frac1p}\\
        &\quad+c\left(\dashint_{I^s_\rho}\left(\sum_{i=0}^{N_\rho}2^{-2si}\dashint_{B_{2^i\rho}}|(\nabla u)_{Q^s_{2^i\rho}}-(\nabla v)_{Q^s_{2^i\rho}}|\,dx\right)^p\,dt\right)^{\frac1p}\\
        &\quad+c\left(\dashint_{I^s_\rho}\left(\sum_{i=0}^{N_\rho}2^{-2si}\dashint_{B_{2^i\rho}}|\nabla v-(\nabla v)_{Q^s_{2^i\rho}}|\,dx\right)^p\,dt\right)^{\frac1p}\\
        &\leq c\left(\sum_{i=0}^{N_\rho}2^{-\frac{(2si)(p+1)}{2}}\dashint_{I^s_\rho}\dashint_{B_{2^i\rho}}|\nabla u-\nabla v|^p\,dx\,dt\right)^{\frac1p}\\
        &\quad+c\sum_{i=0}^{N_\rho}2^{-2si}\osc_{Q^s_{2^i\rho}}\nabla v\eqqcolon J_{1,1}+J_{1,2}
    \end{align*}
    for some constant $c=c(n,s_0,p)$, where we have used the fact that 
    \begin{align*}
        \left(\sum_{i=0}^{N_\rho}2^{-2si}a_i\right)^p\leq \left(\sum_{i=0}^{N_\rho}2^{-\frac{(2si)(p+1)}{2}}a_i^p\right)\left(\sum_{i=0}^{N_\rho}2^{-si}\right)^{{p-1}}.
    \end{align*}
    By Lemma \ref{lem.comp.gra} and \eqref{Defn.exc}, we further estimate $J_{1,1}$ as
    \begin{equation}\label{defn.m1}
    \begin{aligned}
        J_{1,1}&\leq c\sum_{i=0}^{N_\rho}2^{-\frac{2si(p-1)}{2p}}(2^i\rho)^{-(n+2s)/p}\|\nabla(u-v)\|_{L^{p}(Q^s_{2^i\rho})}\\
        &\leq c\sum_{i=0}^{N_\rho}2^{-\frac{2si(p-1)}{2p}}(2^i\rho)^{-(n+2s)/p}\|\nabla(u-v)\|_{L^{p}(Q^s_{1/4})}\\
        &\leq c\rho^{-(n+2s)/p}\sum_{i=0}^{N_\rho}2^{-\frac{2si(p-1)}{2p}}\left[\left({|\mu|(Q^s_{1/2})}\right)^{1-\theta}E^p(\nabla u ;Q^s_{1/2})^\theta+{|\mu|(Q^s_{1/2})}\right]\\
        &\leq c\rho^{-(n+2s+1)}\left[\left({|\mu|(Q^s_{1})}\right)^{1-\theta}{E}(u,\nabla ;Q^s_{1})^\theta+c{|\mu|(Q^s_{1})}\right]\eqqcolon M_1,
    \end{aligned}
    \end{equation}
   where $c=c(n,s_0,p)$. In light of Lemma \ref{lem.exc.dec.gra}, \eqref{defn.q} and the fact that $\alpha_1\leq \alpha_0$ by \eqref{choi.alpha1} and Lemma \ref{lem.comp.gra}, we obtain \begin{equation}\label{choic.m2}
    \begin{aligned}
        J_{1,2}&\leq c\sum_{i=0}^{N_\rho}2^{-2si}(2^i\rho)^{\alpha_1}\left(E^p_{\mathrm{loc}}(\nabla v;Q^s_{3/8})+E^p(\nabla u;Q^s_{3/8})+|\mu|(Q^s_{1/2})\right)\\
        &\quad+c\sum_{i=0}^{N_\rho}2^{-2si}(2^i\rho)^{\alpha_1}\left(\dashint_{I^s_{3/8}}\mathrm{Tail}(u-(\nabla u)_{Q^s_{3/8}}\cdot y-(u)_{B_{3/8}}(t);B_{3/8})^p\,dt\right)^{\frac1p}\\
         &\leq c\rho^{\alpha_1}\left[E^p(\nabla u;Q^s_1)+\overline{E}(u;Q^s_1)\right]+cM_1
    \end{aligned}
    \end{equation}
    for some $c=c(n,s_0,p)$. Combining all the estimates $J_1$ and $J_2$ along with \eqref{Defn.exc}, we arrive at the estimate 
    \begin{align}\label{ineq.ep}
        {E}^p(\nabla u ;Q^s_\rho)\leq c\rho^{\alpha_1}{E}(u,\nabla ; Q^s_1)+cM_1
    \end{align}
    for some constant $c=c(n,s_0,\Lambda,p)$, where the $M_1$ is defined in \eqref{defn.m1}.
    
    Using a non-scaled version of \eqref{tail.est}, we next observe that
    \begin{align*}
    \overline{E}(u;Q^s_\rho)&\leq c\rho^{-1}\left(\dashint_{I^s_\rho}\left(\sum_{i=0}^{N_\rho}2^{-2si}\dashint_{B_{2^i\rho}}|u-l_{i}|\,dy\right)^q\,dt\right)^{\frac1q}\\
    &\quad+c\sum_{i=0}^{N_\rho}2^{i(1-2s)}E_{\mathrm{loc}}(\nabla u;Q^s_{2^i\rho})\\
    &\quad+c\rho^{-1}2^{-2sN_\rho}\left(\dashint_{I^s_\rho}\mathrm{Tail}(u-l_{N_\rho};B_{2^{N_\rho}\rho})^q\,dt\right)^{\frac1q}\eqqcolon L_1+L_2+L_3
\end{align*}
for some constant $c=c(n,s_0)$, where we write 
\begin{align*}
    l_i(y,t)=(\nabla u)_{Q^s_{2^i\rho}}\cdot y+(u)_{B_{2^i\rho}}(t).
\end{align*}
We first estimate $L_1$ as 
\begin{align*}
    L_1&\leq c\rho^{-1}\sum_{i=0}^{N_\rho}2^{-2si}\sup_{t\in I^s_{2^i\rho}}\dashint_{B_{2^i\rho}}|u-l_i-(u-l_i)_{B_{2^i\rho}}(t)|\,dx\\
    &\leq c\rho^{-1}\sum_{i=0}^{N_\rho}2^{-2si}\sup_{t\in I^s_{2^i\rho}}\dashint_{B_{2^i\rho}}|v-\overline{l}_i-(v-\overline{l}_i)_{B_{2^i\rho}}(t)|\,dx\\
    &\quad+c\rho^{-1}\sum_{i=0}^{N_\rho}2^{-2si}\left[(2^i\rho)\dashint_{Q^s_{2^i\rho}}|\nabla(u-v)|\,dz+\sup_{t\in I^s_{2^i\rho}}\dashint_{B_{2^i\rho}}|u-v|\,dx\right]\\
    &\leq c\sum_{i=0}^{N_\rho}2^{(1-2s)i}\left[\sup_{t\in I^s_{2^i\rho}}\dashint_{B_{2^i\rho}}|\nabla v-(\nabla v)_{Q^s_{2^i\rho}}|\,dx+\dashint_{Q^s_{2^i\rho}}|\nabla (u-v)|\,dz\right]\\
    &\quad+c\rho^{-1}\sum_{i=0}^{N_\rho}2^{-2si}\sup_{t\in I^s_{2^i\rho}}\dashint_{B_{2^i\rho}}|u-v|\,dx
\end{align*}
for some constant $c=c(n,s_0,\Lambda)$, where we denote 
\begin{align*}
    \overline{l}_i(y,t)=(\nabla v)_{Q^s_{2^i\rho}}\cdot y+(v)_{B_{2^i\rho}}(t).
\end{align*}
Observe that we have also used the Poincar\'e inequality for the last inequality.
As in the estimate of $J_{1,1}$ and $J_{1,2}$ given in \eqref{defn.m1} and using the third condition in \eqref{choi.alpha1} as well as \eqref{choic.m2}, we further estimate $L_1$ as 
\begin{align*}
    L_1&\leq c\sum_{i=0}^{N_\rho}2^{(1-2s)i}\osc_{Q^s_{2^i\rho}}\nabla v+c\sum_{i=0}^{N_\rho}2^{(1-2s)i}(2^i\rho)^{-(n+2s)}\int_{Q^s_{1/4}}|\nabla(u-v)|\,dz\\
    &\quad+ c\rho^{-1}\sum_{i=0}^{N_\rho}2^{-2si}(2^i\rho)^{-n}|\mu|(Q^s_{1/2})\\
    &\leq c\rho^{\alpha_1}\left[E^p(\nabla u;Q^s_1)+\overline{E}(u;Q^s_1)\right]+cM_1,
\end{align*}
where $c=c(n,s_0,\Lambda,p)$ and the constant $M_1$ is determined in \eqref{defn.m1}. Since 
\begin{equation*}
    L_2\leq  c\sum_{i=0}^{N_\rho}2^{(1-2s)i}\left[\sup_{t\in I^s_{2^i\rho}}\dashint_{B_{2^i\rho}}|\nabla v-(\nabla v)_{Q^s_{2^i\rho}}|\,dx+\dashint_{Q^s_{2^i\rho}}|\nabla (u-v)|\,dz\right]
\end{equation*}
for some constant $c=c(n,s_0,\Lambda,p)$, as in the estimate of $L_1$, we have
\begin{align*}
    L_2\leq c\rho^{\alpha_1}\left[E^p(\nabla u;Q^s_1)+\overline{E}(u;Q^s_1)\right]+cM_1.
\end{align*}
By \eqref{choi.nrho}, the fourth condition given in \eqref{choi.alpha1} and \eqref{222}, we now estimate $L_{3}$ as
\begin{align*}
    L_{3}&\leq c\rho^{2s-1-2s/q}\left(\dashint_{I^s_{3/8}}\mathrm{Tail}\left(u-(\nabla u)_{Q^s_{2^{N_\rho}\rho}}\cdot y-(u)_{B_{2^{N_\rho}\rho}}(t);B_{2^{N_\rho}\rho}\right)^q\,dt\right)^{\frac1q}\\
    &\leq c\rho^{\alpha_1}\left[E^p(\nabla u;Q^s_1)+\overline{E}(u;Q^s_1)\right]
\end{align*}
for some constant $c=c(n,s_0,\Lambda,p)$.
Therefore, combining all the estimates $L_1$, $L_2$  and $L_3$ with \eqref{Defn.exc}, we arrive at the estimate
\begin{align}\label{ine.epo}
    \overline{E}(u;Q^s_\rho)\leq c\rho^{\alpha_1}{E}(u,\nabla ; Q^s_1)+cM_1
\end{align}
for some constant $c=c(n,s_0,\Lambda,p)$.
Finally, the estimate \eqref{exc.ineq} follows by combining \eqref{defn.m1}, \eqref{ineq.ep} and \eqref{ine.epo}.
\end{proof}

\subsection{Pointwise gradient estimates}
We are now able to prove that the averages of $\nabla u$ on any small cylinder can be uniformly controlled by the Riesz potential of $\mu$.
\begin{lemma}\label{lem.exc.sum}
Under the same assumptions as in Lemma \ref{lem.exc.dec.fin}, there is a positive integer $m=m(n,s_0,\Lambda,p)$ such that for any positive integer $j$, we have
\begin{align}\label{ava.ff}
    \left|(\nabla u)_{Q^s_{2^{-mj}R}(z_0)}-(\nabla u)_{Q^s_{R}(z_0)}\right|\leq c{E}(u,\nabla;Q^s_R(z_0))+cI^{|\mu|}_{2s-1,s}(z_0;R)
\end{align}
for some constant $c=c(n,s_0,\Lambda,p)$ which is independent of $j$. 
\end{lemma}
\begin{proof}
We may assume $R=1$ and $z_0=0$.
Let us fix a positive integer $m$ which will be determined later. 
We observe from Lemma \ref{lem.exc.dec.fin} that for any non-negative integer $k$, we have
\begin{align*}
    {E}(u,\nabla;Q^s_{2^{-(k+1)m}})&\leq c2^{-\alpha_1m}{E}(u,\nabla;Q^s_{2^{-km}})\\
    &\quad+c2^{-\alpha_1m(n+2s+1)}\left(\frac{|\mu|(Q^s_{2^{-km}})}{2^{-km(n+1)}}\right)^{1-\theta}{E}(u,\nabla;Q^s_{2^{-km}})^{\theta}\\
    &\quad+c2^{-\alpha_1 m(n+2s+1)}\frac{|\mu|(Q^s_{2^{-km}})}{2^{-km(n+1)}}
\end{align*}
for some constant $c=c(n,s_0,\Lambda,p)$, where the constant $\alpha_1=\alpha_1(n,s_0,\Lambda,p)$ is determined in Lemma \ref{lem.exc.dec.fin}.
We now choose $m=m(n,s_0,\Lambda,p)$ sufficiently large so that $c2^{-\alpha_1m}\leq 1/4$.
Then summing over $k$ yields 
\begin{align*}
    \sum_{k=0}^{j}{E}(u,\nabla;Q^s_{2^{-(k+1)m}})&\leq \frac{1}{4}\sum_{k=0}^{j}{E}(u,\nabla;Q^s_{2^{-km}})\\
    &\quad+c\sum_{k=0}^{j}\left(\frac{|\mu|(Q^s_{2^{-km}})}{2^{-km(n+1)}}\right)^{1-\theta}{E}(u,\nabla;Q^s_{2^{-km}})^{\theta}+c\sum_{k=0}^{j}\frac{|\mu|(Q^s_{2^{-km}})}{2^{-km(n+1)}}
\end{align*}
for some constant $c=c(n,s_0,\Lambda,p)$.
Applying Young's inequality on the second term on the right-hand side of the above inequality gives
\begin{align*}
    \sum_{k=0}^{j}{E}(u,\nabla;Q^s_{2^{-(k+1)m}})&\leq \frac{1}{2}\sum_{k=0}^{j}{E}(u,\nabla;Q^s_{2^{-km}})+c\sum_{k=0}^{j}\frac{|\mu|(Q^s_{2^{-km}})}{2^{-km(n+1)}},
\end{align*}
which implies 
\begin{align}\label{lastsss}
    \sum_{k=0}^{j}{E}(u,\nabla;Q^s_{2^{-(k+1)m}})&\leq c{E}(u,\nabla;Q^s_{1})+cI^{|\mu|}_{2s-1,s}(0;1)
\end{align}
for some constant $c=c(n,s_0,\Lambda,p)$.
Since 
\begin{align*}
    |(\nabla u)_{Q^s_{2^{-mj}}}-(\nabla u)_{Q^s_1}|&\leq \sum_{k=1}^{j}E^p(\nabla u;Q^s_{2^{-km}})+E^p(\nabla u;Q^s_{1})\\
    &\leq \sum_{k=0}^{j-1}{E}(u,\nabla;Q^s_{2^{-(k+1)m}})+{E}(u,\nabla;Q^s_{1}),
\end{align*}
\eqref{ava.ff} follows by combining the above inequality and \eqref{lastsss}. This completes the proof.
\end{proof}

\begin{lemma}\label{lem.weak.last}
    Let $\Omega \subset \mathbb{R}^n$ be a domain and let $T>0$. Moreover, assume that $\mu\in L^1(0,T;L^\infty(\Omega))$ and let 
    \begin{equation*}
        u\in L^2(0,T;W^{s,2}(\Omega))\cap C(0,T;L^2(\Omega))\cap L^1_{2s}(0,T;L^1_{2s}(\bbR^n))
    \end{equation*}
    be a weak solution to 
    \begin{equation*}
        \partial_t u+\mathcal{L}u=\mu\quad\text{in }\Omega_T.
    \end{equation*}
    Then for any positive integer $j$, any $z_0 \in \Omega_T$ and any $R>0$ such that $Q^s_{R}(z_0)\Subset \Omega_T$, we have  
    \begin{equation}\label{ineq.pot.weak}
    \begin{aligned}
        |(\nabla u)_{Q^s_{2^{-jm}R}(z_0)}|&\leq cE(u/R;Q^s_R(z_0))+cI_{2s-1,s}^{|\mu|}(z_0,R)\\
        &\quad+c\int_{0}^{R}\left(\int_{I^s_r(t_0)}R^{-2s}\mathrm{Tail}(u-(u)_{Q^s_R(z_0)};B_R(x_0))\,dt\right)\frac{dr}{r^2}
    \end{aligned}
    \end{equation}
    for some $c=c(n,s_0,\Lambda)$ and some positive integer $m=m(n,s_0,\Lambda)$.
\end{lemma}

\begin{proof}
Fix some $z_0 \in \Omega_T$ and some $R>0$ such that $Q^s_{R}(z_0)\Subset \Omega_T$. 
Then 
\begin{align*}
    u_R(x,t)=u((R/5)x+x_0,(R/5)^{2s}t+t_0)/(R/5)^s
\end{align*}
is a weak solution to 
\begin{align*}
    \partial_t u_R+\mathcal{L}u_R=\mu_R\quad\text{in }Q^s_5,
\end{align*}
where 
\begin{align*}
    \mu_R(x,t)=(R/5)^{s}\mu((R/5)x+x_0,(R/5)^{2s}t+t_0).
\end{align*}
We note that in view of Lemma \ref{lem.comp}, there exists a unique weak solution $v_R$ to \eqref{eq.ivp} with $R$ and $z_0$ replaced by $5$ and $0$, respectively, such that 
\begin{align*}
    \sup_{t\in I^s_5}\int_{B_5}|(u_R-v_R)(x,t)|\,dx\leq c|\mu_R|(Q^s_5)
\end{align*}
for some constant $c=c(n,s_0)$. In addition, using the previous display along with the standard energy inequality for $v_R$ given in the proof of \cite[Theorem 1.8]{KaWe23}, we obtain
\begin{equation}\label{ineq.sss}
\begin{aligned}
    \sup_{t\in I^s_4}\int_{B_{4}}|u_R(x,t)|\,dx&\leq c\sup_{t\in I^s_4}\int_{B_{4}}|v_R(x,t)|\,dx+c|\mu_R|(Q^s_5)\\
    &\leq c\widetilde{E}(v_R;Q^s_5)+c|\mu_R|(Q^s_5)\\
    &\leq c\widetilde{E}(u_R;Q^s_5)+c|\mu_R|(Q^s_5)
\end{aligned}
\end{equation}
for some constant $c=c(n,s_0,\Lambda)$, where the functional $\widetilde{E}(\cdot)$ is defined in \eqref{defn.exctilde}. Let us fix a cutoff function $\xi_1 \in C_c^{\infty}(B_{4})$ with $\xi_1 \equiv 1$ on $B_{3}$ and $|\nabla \xi_1|\leq c$ for some constant $c$. By Lemma \ref{lem.loc} with $R=1$ and $z_0=0$, we obtain that
\begin{align}\label{defn.u1}
    u_1=u_R\xi_1\in L^\infty(I^s_{4};L^1_{2s}(\bbR^n))
\end{align}
is a weak solution to 
\begin{align*}
    \partial_t u_1+\mathcal{L}u_1=\mu_R+f_1\quad\text{in }Q^s_2,
\end{align*}
where $f_1\in L^1(I^s_{5/2};L^\infty(B_{5/2}))$ with 
\begin{align}\label{f1.ineq}
    \|f_1\|_{L^1(I^s_{r};L^\infty(B_{5/2}))}\leq c\int_{I^s_{r}}\mathrm{Tail}(u_R;B_{5/2})\,dt
\end{align}
for some constant $c=c(n,\Lambda)$ and any $r\in(0,5/2]$. In addition, from \eqref{ineq.sss} we deduce 
\begin{align}\label{ineq.supu1}
    \sup_{t\in I^s_{5/2}}\|u_1(\cdot,t)\|_{L^1(\bbR^n)}\leq \sup_{t\in I^s_{4}}\|u_R(\cdot,t)\|_{L^1(B_4)} \leq c\widetilde{E}(u_R;Q^s_{5})+c|\mu_R|(Q^s_5).
\end{align}

We now fix 
\begin{equation}\label{choi.p.last}
    p=\frac12\left(1+\frac{n+2s_0}{n+1}\right)\in \left(1,\frac{n+2s_0}{n+1}\right).
\end{equation}
By Lemma \ref{lem.graex}, we have $\nabla u_1\in L^p(Q^s_2)$ with the estimate
\begin{align}\label{last.u1g}
    \left(\dashint_{Q^s_{r}(z_1)}|\nabla u_1|^p\,dz\right)^{\frac1p}&\leq cE^p(u_1/r;Q^s_{2r}(z_1))+cr^{-(n+1)}(|f_1|+|\mu|)(Q^s_{2r}(z_1))
\end{align}
for some $c=c(n,s_0,\Lambda)$ whenever $Q^s_{2r}(z_1)\subset Q^s_2$. Next, let us fix another cutoff function $\xi_2 \in C_c^{\infty}(B_{8/5})$ with $\xi_2 \equiv 1$ on $B_{6/5}$ and $|\nabla \xi_2|\leq c$ for some constant $c$.
We now employ our localization lemma given by Lemma \ref{lem.loc} once more, this time with $R=2/5$ and $z_0=0$ to obtain that
\begin{equation}\label{defn.u2}
    u_2=u_1\xi_2\in L^p(I^s_{8/5};W^{1,p}(\bbR^n))\cap L^\infty(I^s_{8/5};L^1_{2s}(\bbR^n))
\end{equation} is a weak solution to 
\begin{align}\label{sol.u2}
    \partial_t u_2+\mathcal{L}u_2=\mu_R+f_1+f_2\quad\text{in }Q^s_{4/5},
\end{align}
where $f_2\in L^1(I^s_{1};L^\infty(B_{1}))$ with 
\begin{align}\label{f2.ineq}
    \|f_2\|_{L^1(I^s_{r};L^\infty(B_{1}))}\leq c\int_{I^s_{r}}\mathrm{Tail}(u_1;B_{6/5})\,dt
\end{align}
for some constant $c=c(n,\Lambda)$ and any $r\in(0,1]$. In addition, we observe from \eqref{last.u1g}, \eqref{ineq.sss} and \eqref{f1.ineq} that
\begin{align}\label{1ineq}
    \|u_2\|_{L^p(I^s_{4/5};W^{1,p}(\bbR^n))}\leq c\|u_1\|_{L^p(I^s_{8/5};W^{1,p}(Q^s_{8/5}))}\leq c\widetilde{E}(u_R;Q^s_5)+c|\mu_R|(Q^s_5)
\end{align}
and 
\begin{align}\label{2ineq}
    \sup_{t\in I^s_{4/5}}\|u_2(\cdot,t)\|_{L^1(\bbR^n)}\leq c\sup_{t\in I^s_{4/5}}\|u_1(\cdot,t)\|_{L^1(\bbR^n)}\leq c\widetilde{E}(u_R;Q^s_5)+c|\mu_R|(Q^s_5)
\end{align}
for some constant $c=c(n,s_0,\Lambda)$.
We now apply Lemma \ref{lem.exc.sum} to see that for any $j\geq1$
\begin{align*}
    |(\nabla u_2)_{Q^s_{2^{-jm}4/5}}-(\nabla u_2)_{Q^s_{4/5}}|&\leq c{E}( u_2,\nabla ;Q^s_{4/5})+
cI^{|\nu|}_{2s-1,s}(0;4/5),
\end{align*}
where $m=m(n,s_0,\Lambda)$, the constant $q$ is determined in \eqref{defn.qq} and we write 
\begin{align*}
    |\nu|=|f_1|+|f_2|+|\mu_R|.
\end{align*}
After a few simple calculations together with \eqref{Defn.exc}, \eqref{1ineq} and \eqref{2ineq}, we obtain
\begin{align*}
     |(\nabla u_2)_{Q^s_{4/5}}|+{E}(u_2,\nabla;Q^s_{4/ 5})&\leq c\|u_2\|_{L^p(I^s_{4/5};W^{1,p}(\bbR^n))}+c\sup_{t\in I^s_{4/5}}\|u_2(\cdot,t)\|_{L^1(\bbR^n)}\\
     &\leq c\widetilde{E}(u_R;Q^s_5)+c|\mu_R|(Q^s_5)
\end{align*}
for some $c=c(n,s_0,\Lambda)$. In addition, by \eqref{ineq.sss}, we get
\begin{align*}
    I^{|f_1|+|f_2|}_{2s-1}(0;4/5)&\leq c\int_{0}^{4/5}\frac{1}{r^2}\int_{I^s_r}\left[\mathrm{Tail}(u_R;B_{3})+\mathrm{Tail}(u_1;B_{6/5})\right]\,dt\,dr\\
    &\leq c\int_{0}^{4/5}\frac{1}{r^2}\int_{I^s_r}\mathrm{Tail}(u_R;B_{3})\,dt\,dr\\
    &\quad+c\int_{0}^{4/5}r^{2s-2}\left[\widetilde{E}(u_R;Q^s_5)+|\mu_R|(Q^s_5)\right]\,dr\\
    &\leq c\int_{0}^{4/5}\frac{1}{r^2}\int_{I^s_r}\mathrm{Tail}(u_R;B_{3})\,dt\,dr +c\widetilde{E}(u_R;Q^s_5)+c|\mu_R|(Q^s_5)
\end{align*}
for some constant $c=c(n,s_0,\Lambda)$. Therefore, combining the above two inequalities, we arrive at
\begin{align*}
    |(\nabla u_R)_{Q^s_{2^{-jm}}}|\leq c\widetilde{E}(u_R;Q^s_5)+c\int_{0}^{5}\frac{1}{r^2}\int_{I^s_r}\mathrm{Tail}(u_R;B_{5})\,dt\,dr +cI^{|\mu_R|}_{2s-1,s}(0;5)
\end{align*}
for some constant $c=c(n,s_0,\Lambda)$. Using a scaling argument along with the fact that $u-(u)_{Q^s_R(z_0)}$ is also a weak solution to \eqref{eq:nonlocaleq}, we obtain the desired result.
\end{proof}

We are now ready to prove our gradient potential estimates for SOLA to initial boundary value problems.
\begin{proof}[Proof of Theorem \ref{thm.pt}]
We first fix the constant $p$ determined in \eqref{choi.p.last}.
    By the definition of SOLA, there is a weak solution $u_i$ to \eqref{thm eq : regularized1} with \eqref{regular limit}. Let us fix a parabolic cylinder $Q^s_{2R}(z_0)\Subset \Omega_T$. Then we first prove 
    \begin{equation}\label{ineq0.main}
        \|u_i\|_{L^p(I^s_{R/5}(t_0);W^{1+\sigma_0,p}(B_{R/5}(x_0)))}\leq c\widetilde{E}(u_i;Q^s_{2R}(z_0))+c|\mu_i|(Q^s_{2R}(z_0))
    \end{equation}
    for some constants $\sigma_0=\sigma_0(n,s_0,\Lambda) \in (0,1)$ and  $c=c(n,s_0,\Lambda,R)$ which are independent of $i$. As in the proof of Lemma \ref{lem.weak.last} with $u$ replaced by $u_i$, we get 
    \begin{equation}\label{ineq1.main}
    \begin{aligned}
    R\left(\dashint_{Q^s_{2R/5}(z_0)}|\nabla u_i|^p\,dz\right)^{\frac1p}&\leq c\left(\dashint_{Q^s_{R}(z_0)}|u_i|^p\,dz\right)^{\frac1p}+c\mathrm{Tail}(u_i;Q^s_R(z_0))\\
    &\quad+cR^{-n}|\mu_i|(Q^s_{R}(z_0)),
\end{aligned}
\end{equation}
where $c=c(n,s_0,\Lambda)$ and the constant $p$ is given in \eqref{choi.p.last}.
Indeed, we have also used \eqref{f1.ineq}, \eqref{ineq.supu1} and the fact that $\xi_1=0$ on $\bbR^n\setminus B_4$, where $\xi_1$ is given as in \eqref{defn.u1}. We now use Lemma \ref{lem.comp} with $Q^s_R(z_0)$ replaced by $Q^s_{2R}(z_0)$ and \eqref{bdd.est.zero} to see that 
\begin{equation}\label{ineq3.main}
\begin{aligned}
    \left(\dashint_{Q^s_{R}(z_0)}|u_i|^p\,dz\right)^{\frac1p}&\leq cR^{-n}|\mu_i|(Q^s_{2R}(z_0))+\|v_i\|_{L^\infty(Q^s_R(z_0))}\\
    &\leq cR^{-n}|\mu_i|(Q^s_{2R}(z_0))+c\widetilde{E}(v_i;Q^s_{2R}(z_0))\\
    &\leq c\widetilde{E}(u_i;Q^s_{2R}(z_0))+cR^{-n}|\mu_i|(Q^s_{2R}(z_0))
\end{aligned}
\end{equation}
for some constant $c=c(n,s_0,\Lambda)$. Combining \eqref{ineq3.main} with \eqref{ineq1.main} yields 
\begin{equation*}
    R^{1-n/p}\|\nabla u_i\|_{L^p(Q^s_R(z_0))}\leq c\widetilde{E}(u_i;Q^s_{2R}(z_0))+cR^{-n}|\mu_i|(Q^s_{2R}(z_0)).
\end{equation*}
Similarly, using \eqref{1ineq} with $u$ replaced by $u_i$ and Lemma \ref{lem.highsu}, we obtain 
\begin{equation}\label{ineq2.main}
    \begin{aligned}
    R^{1+\sigma_0-n/p}[\nabla u_i]_{L^p(I^s_{R/5}(t_0);W^{1+\sigma_0,p}(B_{R/5}(x_0)))}&\leq c\widetilde{E}(u_i;Q^s_{2R}(z_0))\\
    &\quad+cR^{-n}|\mu_i|(Q^s_{2R}(z_0)),
\end{aligned}
\end{equation}
where $c=c(n,s_0,\Lambda)$ and $\sigma_0=\sigma_0(n,s_,\Lambda) \in (0,1)$. Using \eqref{ineq1.main}, \eqref{ineq3.main} and \eqref{ineq2.main} now yields \eqref{ineq0.main}. Together with \eqref{regular limit}, this implies that the sequence $\{u_i\}$ is uniformly bounded in the space $L^p(I^s_{R/5}(t_0);W^{1+\sigma_0}(B_{R/5}(x_0)))$.

 Next, we observe from Step 2 in the proof of Theorem \ref{thm: existence} that the sequence $\{\partial_t u_i\}$ is uniformly bounded in $L^1\left(I^s_{R/5}(t_0);\left(X_0^{2s-\sigma_0/2,p_1^\prime}(B_{R/5}(x_0),B_{2R/5}(x_0))\right)^*\right)$, where the associated space $X_0^{2s-\sigma_0/2,p_1^\prime}(B_{R/5}(x_0),B_{2R/5}(x_0))$ is defined in \eqref{defn.space.x0}. Thus, employing a compactness result from \cite[Section 8]{Sim87} with $q=p$, $X=W^{1+\sigma_0,1}(B_{R/5}(x_0))$, $B=W^{1,1}(B_{R/5}(x_0))$ and $Y=\left(X_0^{2s-\sigma_0/2,p_1^\prime}(B_{R/5}(x_0),B_{2R/5}(x_0))\right)^*$ yields the convergence 
\begin{equation*}
    u_i\to u\quad\text{in }L^1(I^s_{R/5}(t_0);W^{1,1}(B_{R/5}(x_0)))
\end{equation*}
up to passing to a subsequence if necessary. Therefore, by standard covering arguments, we conclude that
\begin{equation*}
    \nabla u_i\to\nabla u\quad\text{in }L^1(Q^s_{R/2}(z_0)),
\end{equation*}
whenever $Q^s_{R}(z_0)\Subset \Omega_T$. Using this convergence along with \eqref{regular limit}, \eqref{ineq.pot.weak} and the Lebesgue differentiation theorem as $j \to \infty$ now yields the desired estimate, finishing the proof.
\end{proof}

We are now able to also deduce our gradient potential estimates on the whole space.

\begin{proof}[Proof of Theorem \ref{thm:gradpotwhole}]
We first note from Remark \ref{rmk.relsolaweak} and Theorem \ref{thm.pt} that
\begin{align*}
    |\nabla u(z_0)|&\leq cE(u/R;Q^s_R(z_0))+cI^{|\mu|}_{2s-1,s}(z_0,R)\\
    &\quad+cR^{-1}\sup_{t\in I^s_R(t_0)}\dashint_{\bbR^n\setminus B_R(x_0)}\frac{|u(y,t)-(u)_{Q^s_R(z_0)}|}{|y|^{n+2s}}\,dy
\end{align*}
holds for some constant $c=c(n,s_0,\Lambda)$. Indeed, we have used that $u\in C(0,T;L^2(\bbR^n))$ for the last term in the right-hand side of the above inequality. We next observe that
\begin{align*}
    &E(u/R;Q^s_R(z_0))+R^{-1+2s}\sup_{t\in I^s_R(t_0)}\int_{\bbR^n\setminus B_R(x_0)}\frac{|u(y,t)-(u)_{Q^s_R(z_0)}|}{|y|^{n+2s}}\,dy\\
    &\leq cR^{-(n+1)}\sup_{t\in (0,\infty)}\|u(\cdot,t)\|_{L^2(\bbR^n)},
\end{align*}
    where $c=c(n,s_0,\Lambda)$. The desired estimate follows by combining the above two inequalities and letting $R\to\infty$.\end{proof}

\subsection{Borderline gradient regularity via potentials}

We conclude by proving the various fine regularity results that follow from our gradient potential estimates.

\begin{proof}[Proof of Theorem \ref{thm.conti} and Theorem \ref{thm.vmo}] Since the general case of SOLA can always be obtained by essentially the same approximation procedure as in the proof of Theorem \ref{thm.pt}, we only prove Theorem \ref{thm.conti} and Theorem \ref{thm.vmo} in the case when $u$ is a weak solution to \eqref{eq:nonlocaleq} under the additional assumption that $\mu\in L^1(0,T;L^\infty(\Omega))$.
We first prove Theorem \ref{thm.vmo}. Assume that \eqref{ass1.vmo} and \eqref{ass2.vmo} hold. By Theorem \ref{thm.pt}, this implies that $\nabla u\in L^\infty(Q^s_R(z_0))$.
Let us consider the function $u_2$ which is specified in \eqref{defn.u2} and is a weak solution to \eqref{sol.u2}. Taking into account the fact that 
\begin{equation*}
   u_2(x,t)= u((R/5)x+x_0,(R/5)^{2s}t+t_0)/(R/5)^s \quad\text{on }Q^s_{4/5},
\end{equation*} 
we observe that
\begin{align}\label{ass.bdd}
   \nabla u_2\in L^\infty(Q^s_{4/5}).
\end{align}
We are now going to prove 
\begin{equation}\label{fir.vmo}
    \lim_{r\to0}\sup_{z_2\in Q^s_{1/25}}E( u_2,\nabla ;Q^s_r(z_2))=0.
\end{equation}
For convenience of notation, for the remainder of the proof all constants $c$ will depend only on $n,s_0,\Lambda$ and $R$.
We may assume that $r<1/1000$. We first prove that
\begin{align}\label{ineq0.vmo}
    E( u_2,\nabla ;Q^s_r(z_2))\leq cM,
\end{align}
where we denote 
\begin{align*}
    M\coloneqq\|\nabla u\|_{L^\infty(Q^s_{R}(z_0))}+c\widetilde{E}(u;Q^s_R(z_0))+|\mu|(Q^s_R(z_0)).
\end{align*} To do this, we recall
\begin{equation*}
    E( u_2,\nabla ;Q^s_r(z_2))=E^p(\nabla u_2;Q^s_r(z_2))+\overline{E}(u_2;Q^s_r(z_2)),
\end{equation*}
where $p$ is determined in \eqref{choi.p.last}. 
By \eqref{ass.bdd}, we have 
\begin{align*}
    E^p(\nabla u_2;Q^s_r(z_2))\leq c\|\nabla u_2\|_{L^\infty(Q^s_{4/5})}+\left(\dashint_{I^s_r(t_2)}\mathrm{Tail}(\nabla u_2;B_r(x_2))^p\,dt\right)^{\frac1p}
\end{align*}
for some constant $c=c(n,s_0)$. We further estimate 
\begin{align*}
    &\left(\dashint_{I^s_r(t_2)}\mathrm{Tail}(\nabla u_2;B_r(x_2))^p\,dt\right)^{\frac1p}\\
    &\leq c\|\nabla u_2\|_{L^\infty(Q^s_{4/5})}+\left(\dashint_{I^s_r(t_2)}\left(r^{2s}\int_{\bbR^n\setminus B_{4/5}}\frac{|\nabla u_2(y,t)|}{|y|^{n+2s}}\,dy\right)^p\,dt\right)^{\frac1p}\\
    &\leq c\|\nabla u_2\|_{L^\infty(Q^s_{4/5})}+c\|\nabla u_2\|_{L^p(\bbR^n\times I^s_{4/5})} \leq cM,
\end{align*}
where we have used \eqref{1ineq} to obtain the last inequality. Combining the previous two displays yields
\begin{align}\label{ineq1.vmo}
    E^p(\nabla u_2;Q^s_r(z_2))\leq cM.
\end{align}
In light of \eqref{tail.est}, we have
\begin{align*}
    &r\overline{E}(u_2;Q^s_r(z_2))\\
    &\leq
    c\left[\dashint_{I^s_{r}(t_2)}\left[\sum_{j=0}^{i}2^{-2sj}\dashint_{B_{2^{j}r}(x_2)}|u_2-(\nabla u_2)_{Q^s_{2^jr}(z_2)}\cdot (y-x_2)-(u)_{B_{2^jr}(x_2)}(t)|\,dx\right]^q\,dt\right]^{\frac1q}\\
    &\quad+ cr\sum_{j=0}^{i}2^{j(1-2s)}\dashint_{Q^s_{2^jr}(z_2)}|\nabla u-(\nabla u)_{Q^s_{2^jr}(z_2)}|\,dz\\
    &\quad+ c\left[\dashint_{I^s_{r}(t_2)}2^{-2siq}\mathrm{Tail}(u_2-(\nabla u_2)_{Q^s_{2^ir}(z_2)}\cdot (y-x_2)-(u_2)_{B_{2^ir}(x_2)}(t);B_{2^ir}(x_2))^q\,dt\right]^{\frac1q}\\
    &\eqqcolon J_1+J_2+J_3,
\end{align*}
where $q>1$ is determined in \eqref{defn.qq} and we choose the positive integer $i$ such that 
\begin{equation}\label{cond.i.last}
    1/100\leq2^{i}r<1/50.
\end{equation}
By Poincar\'e's inequality, we have
\begin{align*}
    J_1 &\leq c\left(\dashint_{I^s_{r}(t_2)}\left(\sum_{j=0}^{i}2^{-(2s-1)j}r\dashint_{B_{2^{j}r}(x_2)}|\nabla u_2-(\nabla u_2)_{Q^s_{2^jr}(z_2)}|\,dx\right)^q\,dt\right)^{\frac1q}\\
    &\leq c\sum_{j=0}^{i}2^{-(2s-1)j}r\|\nabla u_2\|_{L^\infty(Q^s_{4/5})} \leq crM.
\end{align*}
Similarly, we have $J_2\leq crM$. Using \eqref{defn.qq}, \eqref{cond.i.last} and \eqref{2ineq}, we estimate $J_3$ as 
\begin{align*}
    J_3
    &\leq cr\left(\int_{I^s_{4/5}}\mathrm{Tail}(u_2-(\nabla u_2)_{Q^s_{2^ir}(z_2)}\cdot (y-x_2)-(u_2)_{B_{2^ir}(x_2)}(t);B_{2^ir}(x_2))^q\,dt\right)^{\frac1q}\leq crM.
\end{align*}

Combining all the estimates of $J_i$ for each $i=1,2,3$ leads to the estimate
\begin{align}\label{ineq2.vmo}
    \overline{E}(u_2;Q^s_r(z_2))\leq cM.
\end{align}
Combining \eqref{ineq1.vmo} and \eqref{ineq2.vmo} now yields \eqref{ineq0.vmo}. 
Next, we use \eqref{exc.ineq} and Young's inequality to obtain
\begin{align*}
    E(u_2,\nabla;Q^s_{\rho r}(z_2))&\leq c\rho^{\alpha_1}E(u_2,\nabla;Q^s_{ r}(z_2))\\
    &\quad+c\rho^{-(n+2s+1+\alpha_1\theta/(1-\theta))}\frac{\left(|\mu_R|+|f_1|+|f_2|\right)(Q^s_r(z_2))}{r^{n+1}},
\end{align*}
where the constants $\alpha_1$ and $\theta$ are determined in Lemma \ref{lem.exc.dec.fin}. In view of \eqref{ineq0.vmo}, \eqref{f1.ineq}, \eqref{f2.ineq} and \eqref{ineq.supu1}, we have 
\begin{equation}\label{ineq3.vmo}
\begin{aligned}
    & E(u_2,\nabla;Q^s_{\rho r}(z_2)) \\&\leq c\rho^{\alpha_1}M+c\rho^{-N}\frac{|\mu|(Q^s_{rR}(z_2+z_0))}{r^{n+1}}\\
    &\quad+\frac{c}{r}\rho^{-N}\left[\int_{I^s_r(t_2)}\mathrm{Tail}(u_R;B_{5/2})\,dt+\int_{I^s_r(t_2)}\mathrm{Tail}(u_1;B_{6/5})\,dt\right]\\
    &\leq c\rho^{\alpha_1}M+c\rho^{-N}r^{2s-1}M+c\rho^{-N}\frac{|\mu|(Q^s_{rR}(z_2+z_0))}{(rR)^{n+1}}\\
    &\quad+\frac{c}{rR\rho^N}\int_{I^s_{Rr}(t_2+t_0)}\mathrm{Tail}(u-(u)_{Q^s_{R/2}(z_2+z_0)};B_{R/2}(x_2+x_0))\,dt,
\end{aligned}
\end{equation}
where we write $N\coloneqq n+2s+1+\alpha_1\theta/(1-\theta)$. For any $\epsilon>0$, we choose $\rho$ sufficiently small so that 
\begin{align*}
    \rho^{\alpha_1}\leq \epsilon(4cM).
\end{align*}
We next choose $r$ sufficiently small so that 
\begin{align*}
    r^{2s-1}\leq \epsilon\rho^{N}/(4cM )
\end{align*}
and
\begin{equation}\label{sec.r}
\begin{aligned}
    &\frac{c}{\rho^{N}}\frac{|\mu|(Q^s_{rR}(z_2+z_0))}{(rR)^{n+1}}\\
    &\quad+\frac{c}{\rho^N}\frac{1}{Rr}\int_{I^s_{Rr}(t_2+t_0)}\mathrm{Tail}(u-(u)_{Q^s_{R/2}(z_2+z_0)};B_{R/2}(x_2+x_0))\,dt\leq \frac{\epsilon}4.
\end{aligned}
\end{equation}
Here we have used the assumption \eqref{ass2.vmo} to find some $r$ satisfying \eqref{sec.r}. Thus, for any given $\epsilon>0$, there are constants $\rho$ and $r$ which are independent of the point $z_2$ such that
\begin{align*}
    E(u_2,\nabla; Q^s_{\rho r}(z_2))\leq \epsilon.
\end{align*}
This implies \eqref{fir.vmo} holds and that $\nabla u_2$ is VMO-regular in $Q^s_{1/25}$. Thus $\nabla u$ is also VMO-regular in $Q^s_{R/125}(z_0)$, which completes the proof of Theorem \ref{thm.vmo}. 

It remains to prove Theorem \ref{thm.conti}. We now assume \eqref{ass.coni}. Then \eqref{ass1.vmo} and \eqref{ass2.vmo} are true, which implies \eqref{fir.vmo}. Therefore, after a few calculations together with \eqref{ava.ff} and \eqref{ineq3.vmo}, we obtain 
\begin{align*}
    &|(\nabla u_2)_{Q^s_r(z_2)}-\nabla u_2(z_2)|\\
    &\leq cE(u_2,\nabla ;Q^s_r(z_2))+cI^{|\mu_R|+|f_1|+|f_2|}_{2s-1,s}(z_2;r)\\
    &\leq cE(u_2,\nabla ;Q^s_r(z_2))+cr^{2s-1}M+cI^{|\mu|}_{2s-1,s}(z_0+z_2;rR)\\
&\quad+c\int_{0}^{rR}\frac{1}{\varrho^2}\int_{I^s_\varrho(t_0+t_2)}\mathrm{Tail}(u-(u)_{Q^s_{R/2}(z_0+z_2)};B_{R/2}(x_0+x_2))\,dt\,d\varrho.
\end{align*}
Thus, using \eqref{fir.vmo} and \eqref{ass.coni}, we conclude that $\nabla u_2$ is continuous on $Q^s_{1/25}$, which implies that $\nabla u$ is continuous on $Q^s_{R/125}(z_0)$ and therefore also in $\Omega_T$. In view of a standard covering argument, the proof is complete.
\end{proof}

\begin{proof}[Proof of Corollary \ref{cor:Lorentz}] 
    Fix $Q^s_{2R}(z_0)\Subset\Omega_T$ and $z_1 \in Q^s_{R}(z_0)$. Then by Theorem \ref{thm.pt}, 
    \begin{align*}
        |\nabla u(z_1)|\leq cE(u/R;Q^s_{2R}(z_0))+cI^{|\mu|}_{2s-1,s}(z_1;R)+cI^{|\nu|}_{2s-1,s}(z_1;R),
    \end{align*}
    where we denote 
    \begin{equation*}
        \nu(x,t)=\int_{t-R^{2s}}^{t}\int_{\bbR^n\setminus B_{2R}(x_0)}\frac{|u(y,\tau)-(u)_{Q^s_{2R}(z_0)}|}{|y-x_0|^{n+2s}}\,dy\,d\tau.
    \end{equation*}
    Therefore, if $u\in L^{p,q}(0,T;L^1_{2s}(\bbR^n))$, then also
    \begin{equation*}
        \nu\in L^{p,q}(Q^s_{2R}(z_0)).
    \end{equation*}
    Setting $N\coloneqq N_{\mathrm{par},s}= n+2s$, we deduce from \cite[Equations (1.19), (1.20), (6.11), (6.12)]{DM2} that for any $f \in \mathcal{M}(\mathbb{R}^{n+1})$, we have
    \begin{align}\label{ineq.rie1}
        \|I^f_{\beta,s}\|_{L^{\frac{Np}{N-\beta p},q}(\bbR^{n+1})}\leq c\|f\|_{L^{p,q}(\bbR^{n+1})}
    \end{align}
    whenever $p>1$ and $\beta q<N$, and 
    \begin{align}\label{ineq.rie2}
        \|I^f_{\beta,s}\|_{L^{\frac{N}{N-\beta},\infty}(\bbR^{n+1})}\leq c\|f\|_{\mathcal{M}(\bbR^{n+1})}
    \end{align}
    whenever $\beta<N$. By the same reasoning as in \cite[Equation (6.13)]{DM2} together with \eqref{ineq.rie1}, \eqref{ineq.rie2}, we obtain
    \begin{align*}
        \mu,\nu\in \mathcal{M}(Q^s_{2R}(z_0))\Longrightarrow \nabla u\in L^{\frac{n+2s}{n+1},\infty}(Q^s_R(z_0))
    \end{align*}
    and 
    \begin{align*}
        \mu,\nu\in L^{p,q}(Q^s_{2R}(z_0))\Longrightarrow \nabla u\in L^{\frac{p(n+2s)}{n+2s-(2s-1)p},q}(Q^s_R(z_0))
    \end{align*}
    by taking $\beta=2s-1$. This completes the proof.
\end{proof}

\begin{proof}[Proof of Corollary \ref{cor:L}]
We first define 
\begin{equation*}
    \nu(I_r^s(t_0))\coloneqq \int_{I_r^s(t_0)}R^{-2s}\mathrm{Tail}(u-(u)_{Q_R^s(z_0)};B_R(x_0))\,dt,
\end{equation*}
which yields a measure defined in one dimension. Then we observe
\begin{align*}
    \int_{0}^{R}\frac{\nu(I_r^s(t_0))}{r^{2}}\,dr=\int_{0}^{R}\left(\int_{I_r^s(t_0)}R^{-2s}\mathrm{Tail}(u-(u)_{Q_R^s(z_0)};B_R(x_0))\,dt\right)\frac{dr}{r^2}.
\end{align*}
 By following the same lines as in the proof of Lemma \ref{lem.lorentz}, we have 
 \begin{equation}\label{lor.inq}
     \int_{0}^{R}\frac{\nu(I_r^s(t_0))}{r^{2}}\,dr \leq c\|\nu\|_{L^{\frac{2s}{2s-1},1}(I_{R}^s(t_0))}
 \end{equation}
 for some constant $c=c(n,s)$. Therefore, using Theorem \ref{thm.pt} and Lemma \ref{lem.lorentz} together with \eqref{lor.inq}, we obtain the desired result.
\end{proof}



\printbibliography

\end{document}